# OPÉRATEURS D'ENTRELACEMENT ET ALGÈBRES DE HECKE AVEC PARAMÈTRES D'UN GROUPE RÉDUCTIF $p$-ADIQUE - LE CAS DES GROUPES CLASSIQUES

Volker Heiermann

Abstract. For $G$ a symplectic or orthogonal $p$-adic group (not necessarily split), or an inner form of a general linear $p$-adic group, we compute the endomorphism algebras of some induced projective generators à la Bernstein of the category of smooth representations of $G$ and show that these algebras are isomorphic to the semi-direct product of a Hecke algebra with parameters by a finite group algebra. Our strategy and parts of our intermediate results apply to a general reductive connected $p$-adic group.

RÉSUMÉ: Pour $G$ un groupe symplectique ou orthognal $p$-adique (déployé ou non) ou une forme intérieure d'un groupe linéaire $p$-adique, nous calculons les algèbres d'endomorphismes de certains générateurs projectifs induits à la Bernstein dans la catégorie des représentations lisses de $G$, et nous montrons que ces algèbres sont isomorphes au produit semi-direct de l'algèbre d'un groupe fini avec une algèbre de Hecke avec paramètres. Notre stratégie et une bonne partie des résultats intermédiaires s'appliquent à un groupe réductif connexe $p$-adique arbitraire.

*À la mémoire de Jacqueline*

Fixons un corps local non archimédien $F$. Soit $G$ (le groupe des points d') un groupe réductif connexe défini sur $F$. Il a été montré par J. Bernstein [BD] que la catégorie $Rep(G)$ des représentations lisses complexes de $G$ est le produit direct de certaines sous-catégories pleines, notées ci-dessous $Rep_G({}^W\mathcal{O})$. Celles-ci sont généralement conjecturées d'être isomorphes à des catégories de modules à droite sur des algèbres de Hecke avec paramètres (ou plus généralement sur des algèbres de Hecke *étendues,* i.e. les produits d'une telle algèbre avec une algèbre de groupe fini, éventuellement tordue par un 2-cocycle). L'intérêt de ceci vient du fait que la théorie des représentations des algèbres de Hecke avec paramètres a été largement étudiée par G. Lusztig et d'autres. (Citons comme exemple [L2], [K] et [OS].)

Une approche pour prouver cette équivalence de catégorie est la "théorie des types" de C. Bushnell et Ph. Kutzko [BK]. D'un autre côté, J. Bernstein a explicité en 1992 dans un cours à l'université de Harvard certaines représentations $\Pi_{\mathcal{O}}$ de $G$ qui possèdent la propriété que $Rep_G({}^W\mathcal{O})$ est isomorphe à la catégorie des modules à droite sur l'algèbre d'entrelacement $\mathrm{End}_G(\Pi_{\mathcal{O}})$.

Le but de ce travail est de calculer l'algèbre $\mathrm{End}_G(\Pi_{\mathcal{O}})$ explicitement pour montrer qu'il s'agit bien d'une algèbre de Hecke avec paramètres étendue. Pour cela, je me restreins au cas où $G$ est un groupe symplectique ou orthogonal (déployé ou non) ou encore où $G$ est une forme intérieure de $\mathrm{GL}_n(F)$. Remarquons que, sauf dans le cas d'une forme intérieure de $\mathrm{GL}_n(F)$, cette équivalence de catégorie n'est à la connaissance de l'auteur pas encore accessible par la théorie des types.

En fait, les hypothèses dont j'ai besoin devraient également être vérifiées pour d'autres formes de ces groupes classiques. Par ailleurs, l'approche et une bonne partie des résultats intermédiaires sont générales. Certaines preuves ont d'ailleurs été motivées par l'obtention d'un résultat général et pourraient bien être simplifiées sous les hypothèses imposées dans ce papier.

Remarquons que les paramètres de ces algèbres de Hecke avec paramètres obtenues ici ont été explicités dans [H3] (voir aussi **7.5**).

Soyons maintenant plus précis: soient $P = MU$ un sous-groupe parabolique de $G$ et $(\sigma, E)$ une représentation irréductible cuspidale de $M$. Notons $\mathfrak{X}^{\mathrm{nr}}(M)$ le groupe des caractères non ramifiés de $M$, $W$ le groupe de Weyl de $G$ défini relatif à un tore déployé maximal contenu dans $M$ et $\mathcal{O}$ l'ensemble des classes d'équivalence de représentations de la forme $\sigma \otimes \chi$, $\chi \in \mathfrak{X}^{\mathrm{nr}}(M)$. Alors, ${}^W\mathcal{O}$ est l'orbite de $\mathcal{O}$ pour l'action par $W$, ${}^W\mathcal{O} = \{ {}^w\mathcal{O} \mid w \in W\}$.

Désignons par $i_P^G$ le foncteur de l'induction parabolique normalisé et par $r_P^G$ son adjoint à gauche. La catégorie $Rep_G({}^W\mathcal{O})$ est alors la sous-catégorie pleine de la catégorie $Rep(G)$ des représentations lisses complexes de $G$ dont les objets sont les représentations $\pi$ qui vérifient la propriété suivante: l'ensemble des classes d'équivalence des sous-quotients irréductibles de $r_{P'}^G\pi$ est contenu dans ${}^W\mathcal{O}$ si $P'$ est associé à $P$, et il ne contient aucun sous-quotient irréductible cuspidal sinon.

Notons $M^1$ l'intersection des noyaux des caractères non ramifiés de $M$ et $(\sigma_1, E_1)$ une composante irréductible de la restriction de $\sigma$ à $M^1$. Désignons par $ind_{M^1}^M$ le foncteur de l'induction compacte. Alors, il a été montré par Bernstein [Ro, 1.6] que la catégorie $Rep_G({}^W\mathcal{O})$ est isomorphe à la catégorie des $\mathrm{End}_G(i_P^G\,\mathrm{ind}_{M^1}^M E_1)$-modules à droite. Remarquons que $ind_{M^1}^M E_1$ ne dépend pas du choix de $(\sigma_1, E_1)$. Cet isomorphisme de catégories est par ailleurs compatible avec l'induction parabolique et le foncteur de Jacquet [Ro, 5.3].

L'objet d'étude de ce travail est alors l'algèbre $\mathrm{End}_G(i_P^G\,\mathrm{ind}_{M^1}^M E_1)$. Je me restreinds pour cela au cas où $G$ est un groupe symplectique ou orthogonal (déployé ou non) ou encore où $G$ est une forme intérieure de $\mathrm{GL}_n(F)$. Dans cette situation, je donne une base de cette algèbre en tant que module sur un anneau de fonctions régulières sur $\mathcal{O}$ (cf. proposition **7.4**) pour en déduire que cette algèbre est isomorphe au produit semi-direct d'une algèbre de Hecke avec paramètres (au sens de Lusztig [L1]) avec l'algèbre d'un groupe fini (cf. théorème **7.7**). (Le groupe fini est une sorte de "$R$-groupe".) Nous finissons notre article avec quelques remarques disant que l'isomorphisme de catégories qui est déduit de notre isomorphisme d'algèbres est compatible avec l'induction parabolique et le foncteur de Jacquet (cf. **7.9**).

Le déroulement antérieur de l'article est comme suit: le paragraphe 1 sert à fixer quelques définitions et notations et à établir ou rappeler des résultats basiques sur les opérateurs d'entrelacement standard et sur les représentations lisses d'un groupe $p$-adique. Le paragraphe 2 est un préliminaire au paragraphe 3 où nous déduisons des opérateurs d'entrelacement standard de Harish-Chandra des générateurs de l'algèbre d'entrelacement de la "grande" représentation $i_P^G ind_{M^1}^M E_{|M^1}$. Ceux-ci ont toutefois des pôles. En particulier, ils n'engendrent l'algèbre que sur le corps des fonctions rationnelles de $\mathcal{O}$. Au paragraphe 4, nous nous spécialisons à l'algèbre d'entrelacement de $i_P^G ind_{M^1}^M E_1$, et nous explicitons des générateurs qui ont toutefois toujours des pôles. Au paragraphe 5, nous définissons des opérateurs $T_w$, sans pôles et motivés par certaines relations que doivent vérifier les générateurs d'une algèbre de Hecke avec paramètres, pour montrer ensuite que ces opérateurs engendrent l'algèbre $\mathrm{End}_G$ $(i_P^G\,\mathrm{ind}_{M^1}^M\,E_1)$ (éventuellement ensemble avec certains des opérateurs d'entre-lacement ci-dessus reliés au "R-groupe"). Remarquons que la preuve que les opéra-teurs $T_w$ engendrent l'algèbre $\mathrm{End}_G(i_P^G\,\mathrm{ind}_{M^1}^M\,E_1)$ est cruciale dans ce travail. A cette occasion, une propriété des opérateurs d'entrelacement prouvée dans [H1] sera utile. Le paragraphe 6 sert à définir la donnée radicielle qui sera attachée à la description de $\mathrm{End}_G(i_P^G\,\mathrm{ind}_{M^1}^M\,E_1)$ comme algèbre de Hecke avec paramètres. Au paragraphe 7, la définition d'algèbre de Hecke avec paramètres de Lusztig sera rappelée avant d'établir les dernières propriétés indispensables de ces opérateurs $T_w$.

L'auteur remercie J.-L. Waldspurger pour lui avoir indiqué sur l'exemple de $\mathrm{SL}_2$ la possibilité de retrouver des algèbres de Hecke avec paramètres à partir d'opérateurs d'entrelacement. Il a également profité de discussions avec P. Schneider et E.-W. Zink sur différents aspects de cet article, et il remercie G. Henniart et J.-L. Waldspurger pour des remarques sur une version préliminaire de ce papier. Les critiques de différents rapporteurs ont finalement aidé à améliorer la version finale.

**1.** Nous gardons les hypothèses et notations de l'introduction. En outre, nous fixons un sous-groupe parabolique minimal $P_0$ contenu dans $P$, un sous-groupe de Levi $M_0$ de $P_0$ contenu dans $M$, $P_0 = M_0U_0$, et un tore $A_0$ déployé maximal (sur $F$) de $M_0$. Le sous-groupe parabolique $P$ est donc un sous-groupe parabolique standard (i.e. $P \supseteq P_0$), et son sous-groupe de Levi $M$ est un Levi standard (i.e.

$M \supseteq M_0$). Le tore déployé maximal contenu dans le centre de $M$ sera désigné par $A_M$. Nous noterons $W = W^G(= N_G(A_0)/Z_G(A_0))$ le groupe de Weyl de $G$ défini relatif à $A_0$, et on le munit de la longueur fixée par le choix de $P_0$. On notera cet ordre $l^G$ ou simplement $l$. Nous désignerons par $K$ un sous-groupe compact maximal de $G$ qui est en bonne position par rapport à $A_0$ par la théorie de Bruhat-Tits. Par ailleurs, lorsque $H$ est un groupe, nous écrirons $X(H)$ pour le groupe $\mathrm{Hom}(H, \mathbb{C}^\times)$.

L'ensemble des racines non triviales de $A_M$ dans l'algèbre de Lie de $G$ sera désigné par $\Sigma(A_M)$, le sous-ensemble des racines qui agissent dans l'algèbre de Lie de $U$, par $\Sigma(P)$, et l'ensemble des racines réduites par $\Sigma_{red}(P)$. On a une bijection $\alpha \mapsto M_\alpha$ entre $\Sigma_{red}(P)$ et l'ensemble des sous-groupes de Levi de $G$ qui contiennent $M$ et qui sont minimaux pour cette propriété.

Par ailleurs, on désignera par $a_M$ l'algèbre de Lie réelle de $A_M$, par $a_M^*$ son dual, par $a_{M,\mathbb{C}}^*$ son complexifié, par $|\cdot|_F$ le module de $F$, par $q$ le cardinal du corps résiduel de $F$ et par $H_M$ l'application $M \to a_M$ qui vérifie $q^{-\langle H_M(m),\alpha\rangle} = |\alpha(m)|_F$ pour tout caractère rationnel $\alpha$ de $M$ et tout $m \in M$. Le noyau de $H_M$ est égal à $M^1$. Si $s$ est un nombre complexe, on note $\chi_{\alpha\otimes s}$ le caractère non ramifié $m \mapsto |\alpha(m)|^s$ de $M$. Cette application se prolonge en un homomorphisme de groupes $a_{M,\mathbb{C}}^* \to \mathfrak{X}^{\mathrm{nr}}(M)$ qui est surjectif.

**1.1** Nous noterons $W(M)$ l'ensemble des représentants dans $W$ des éléments du groupe quotient $\{w \in W | w^{-1}Mw = M\}/W^M$ qui sont de longueur minimale dans leur classe à droite modulo $W^M$. Observons que $W(M)$ est un sous-groupe de $W$. (C'est aussi le sous-groupe de $W$, formé des éléments qui fixent $M$ et $P_0 \cap M$.) Nous désignons par $W(M, \mathcal{O})$ le sous-groupe de $W(M)$ formé des éléments qui stabilisent $\mathcal{O}$.

Pour $w \in W(M)$, posons $l_M(w) := |\Sigma_{red}(P) \cap \Sigma_{red}(w\overline{P}w^{-1})|$.

**Proposition:** *Pour que deux éléments $w$ et $w'$ de $W(M)$ vérifient $l_M(ww') = l_M(w') + l_M(w)$, il faut et il suffit que $l(ww') = l(w) + l(w')$.*

*Preuve:* Remarquons d'abord que, pour tout $w \in W(M)$, $w\Sigma_{red}(M \cap P_0) = \Sigma_{red}(M \cap P_0)$ et que $l(w) = |\Sigma_{red}(P_0) \cap \Sigma_{red}(w\overline{P_0})|$. Il s'ensuit tout d'abord que $\Sigma^M(A_0) \cap \Sigma_{red}(P_0) \cap \Sigma_{red}(w\overline{P_0}) = \emptyset$. Donc, un élément $\alpha$ de $\Sigma(A_0)$ est dans $\Sigma_{red}(P_0) \cap \Sigma_{red}(w\overline{P_0})$, si et seulement si $\alpha_{|A_M} \in \Sigma(P) \cap \Sigma(w\overline{P})$. Par suite, pour $w, w' \in W(M)$, $l_M(ww') = l_M(w) + l_M(w')$ équivaut à

$$|\Sigma_{red}(P_0) \cap \Sigma_{red}(ww'\overline{P_0})| = |\Sigma_{red}(P_0) \cap \Sigma_{red}(w\overline{P_0})| + |\Sigma_{red}(P_0) \cap \Sigma_{red}(w'\overline{P_0})|,$$

d'où la proposition par l'expression pour $l(w)$ rappelée au début. □

Pour tout $w \in W$, notons $P(w) = M(w)U(w)$ le sous-groupe parabolique standard minimal tel que $P, wP \subseteq P(w)$. (Si $M_\alpha$ est un sous-groupe de Levi standard et $w$ l'élément non trivial de $W^{M_\alpha}(M)$, alors $M(w) = M_\alpha$.) On peut choisir des représentants $\overline{w}$ des éléments de $W$ dans $G$ tels que $\overline{w} \in K \cap M(w)$, ce que l'on fera désormais. On identifiera dans la suite les éléments de $W$ avec leurs représentants dans $K$. Il sera toujours clair d'après le contexte, si $w$ désigne un élément de $W$ ou de $G$, en faisant les conventions suivantes: pour $w, w' \in W$, le symbole $ww'$, considéré comme élément de $K$, désigne le produit du représentant de $w$ dans $K$ avec celui de $w'$ et **non pas** le représentant de $ww'$. Par ailleurs, le symbole $w^{-1}$ correspondra à l'inverse (du représentant) de $w$ dans le groupe $G$ qui est bien un élément de $K \cap M_w$.

Remarquons que notre construction d'opérateurs de $\mathrm{End}_G(i_P^G E_B)$ sera essentiellement indépendante du choix de l'ensemble des représentants de $W$. On fera toutefois dans **1.15** quelques restrictions supplémentaires sur ces représentants.

Lorsque $P_1 = M_1 U_1$ est un sous-groupe parabolique semi-standard de $G$ (i.e. $M_1$ contient $A_0$), $(\pi, V)$ une représentation lisse de $M_1$ et $w \in W$, on désignera par $\lambda(w)$ l'isomorphisme $i_{P_1}^G V \to i_{wP_1}^G wV$, $v \mapsto v(w^{-1}\cdot)$, qui dépend du choix du représentant de $w$ dans $G$. On peut également définir $\lambda(g)$ pour tout élément $g$ de $G$. Si $g \in M_1$, alors $\lambda(g)$ est un isomorphisme $i_{P_1}^G V \to i_{P_1}^G gV$ qui est induit par fonctorialité de l'isomorphisme $\pi(g^{-1}) : V \to gV$. On écrira également $i_{P_1}^G(\pi(g^{-1}))$ (à ne pas confondre avec $i_P^G \pi(g^{-1})$).

**1.2** Rappelons que $\mathcal{O}$ est le quotient d'un tore complexe par un groupe fini. La fonction $\mu$ de Harish-Chandra est une fonction rationnelle définie sur $\mathcal{O}$. Nous écrirons $\mu^H$, si elle est définie par rapport à un sous-groupe réductif $H$ de $G$ qui contient $M$. Nous omettrons cet exposant, si $G = H$.

**Théorème:** ([Si1, 5.4.2.2 et 5.4.2.3]) (Harish-Chandra) *Soit $\alpha \in \Sigma(P)$ et soit $\sigma_1$ une représentation irréductible cuspidale de $M$.*

*a) Si $\mu^{M_\alpha}(\sigma_1) = 0$, alors il existe un unique élément non trivial $s_\alpha$ dans $W^{M_\alpha}$ $(M)$ tel que $s_\alpha(P \cap M_\alpha) = \overline{P} \cap M_\alpha$ et $s_\alpha \sigma_1 \simeq \sigma_1$.*

*b) Supposons qu'il existe un unique élément non trivial $s_\alpha$ dans $W^{M_\alpha}(M)$ tel que $s_\alpha(P \cap M_\alpha) = \overline{P} \cap M_\alpha$ et $s_\alpha \sigma_1 \simeq \sigma_1$. Alors, pour que $\mu^{M_\alpha}(\sigma_1) \neq 0$, il faut et il suffit que la représentation $i_{P \cap M_\alpha}^{M_\alpha} \sigma_1$ soit réductible. La représentation $i_{P \cap M_\alpha}^{M_\alpha} \sigma_1$ est alors somme directe de deux représentations irréductibles non isomorphes.*

**1.3** Notons $a_M^{M_\alpha}$ le sous-espace de $a_M$ qui est orthogonal à $a_{M_\alpha}^*$.

**Proposition:** *L'ensemble $\Sigma_{\mathcal{O},\mu} := \{\alpha \in \Sigma_{red}(A_M) |\ \mu^{M_\alpha}$ a un zéro sur $\mathcal{O}\}$ est un système de racines. Pour $\alpha \in \Sigma_{\mathcal{O},\mu}$, désignons par $s_\alpha$ l'unique élément de $W^{M_\alpha}$ $(M, \mathcal{O})$ qui conjugue $P \cap M_\alpha$ et $\overline{P} \cap M_\alpha$. Le groupe de Weyl $W_{\mathcal{O}}$ de $\Sigma_{\mathcal{O},\mu}$ s'identifie au sous-groupe de $W(M, \mathcal{O})$ engendré par les réflexions $s_\alpha$. Pour tout $\alpha \in \Sigma_{\mathcal{O},\mu}$, notons $\alpha^\vee$ l'unique élément de $a_M^{M_\alpha}$ qui vérifie $\langle \alpha, \alpha^\vee \rangle = 2$. Alors $\Sigma_{\mathcal{O},\mu}^\vee := \{\alpha^\vee | \alpha \in \Sigma_{\mathcal{O},\mu}\}$ est l'ensemble des coracines de $\Sigma_{\mathcal{O},\mu}$, la dualité étant celle*

*entre* $a_M$ *et* $a_M^*$.

*L'ensemble* $\Sigma(P)\cap\Sigma_{\mathcal{O},\mu}$ *est l'ensemble des racines positives pour un certain ordre sur* $\Sigma_{\mathcal{O},\mu}$.

*Preuve:* La preuve de la première partie de la proposition est analogue à celle de la proposition **4.2** dans [H2] (voir également [Si3]), après avoir vérifié que $\Sigma_{\mathcal{O}}$ est stable pour l'action de $W(M,\mathcal{O})$. Ceci résulte de la $W(M,\mathcal{O})$-invariance de la fonction $\mu$ de Harish-Chandra [W, V.2.1].

Soit $\alpha\in\Sigma_{\mathcal{O},\mu}$. Le sous-espace vectoriel de $a_M^*$ engendré par $\alpha$ est $a_M^{M_\alpha *}$. Celui-ci est l'orthogonal de $a_{M_\alpha}$. Par suite, $\alpha^\vee$ doit appartenir à $a_M^{M_\alpha}$. Comme $\langle\alpha,\alpha^\vee\rangle=2$, c'est donc bien la coracine.

La dernière assertion vient du fait que $\Sigma(P)$ définit un ordre sur $\Sigma_{red}(A_M)$, donc sur $\Sigma_{\mathcal{O},\mu}$. □

**1.4** *Définition:* On pose $\Sigma_{\mathcal{O},\mu}(P)=\Sigma(P)\cap\Sigma_{\mathcal{O},\mu}$, et on note $\Delta_{\mathcal{O},\mu}$ la base de $\Sigma_{\mathcal{O},\mu}$ déterminée par l'ordre pour lequel l'ensemble des racines positives est $\Sigma_{\mathcal{O},\mu}(P)$. La longueur sur $W_{\mathcal{O}}$ sera désignée par $l_{\mathcal{O}}$.

On fixe par ailleurs une représentation unitaire $\sigma$ dont la classe d'équivalence appartient à $\mathcal{O}$ telle que $\mu^{M_\alpha}(\sigma)=0$ pour tout $\alpha\in\Delta_{\mathcal{O},\mu}$. (Ceci est possible, puisque $\Delta_{\mathcal{O},\mu}$ est une base et que $\mu^{M_\alpha}$ est invariante pour l'action par $\mathfrak{X}^{\mathrm{nr}}(M_\alpha)$.)

**1.5** Pour $\alpha\in\Sigma_{\mathcal{O},\mu}$, fixons un élément $h_\alpha$ de $M\cap M_\alpha^1$ tel que $H_M(h_\alpha)$ soit un multiple de $\alpha^\vee$ par un nombre réel $>0$ et que, pour tout $\chi\in\mathfrak{X}^{\mathrm{nr}}(M)$, $\chi(h_\alpha)=1$ équivaut à $\chi\in\mathfrak{X}^{\mathrm{nr}}(M_\alpha)$. Ceci est bien défini, puisque le noyau de $H_M$ est égal à $M^1$ et que $M\cap M_\alpha^1/M^1$ est un $\mathbb{Z}$-module libre de rang 1 (cf. [H2, 1.2]). On voit par ailleurs que $H_M(h_\alpha)$ est unique.

Notons $\mathrm{Stab}(\mathcal{O})$ le groupe des caractères non ramifiés $\chi$ de $M$, vérifiant $\sigma\otimes\chi\simeq\sigma$ (il ne dépend que de $\mathcal{O}$), et $t_\alpha$ le plus petit entier $\geq 1$ tel que $\chi(h_\alpha^{t_\alpha})=1$ équivaut à $\chi\in\mathfrak{X}^{\mathrm{nr}}(M_\alpha)\,\mathrm{Stab}(\mathcal{O})$. Notons $b_{h_\alpha}$ l'application $\mathfrak{X}^{\mathrm{nr}}(M)\rightarrow\mathbb{C}$, $\chi\mapsto\chi(h_\alpha)$, et posons $Y_\alpha=b_{h_\alpha}$ et $X_\alpha=Y_\alpha^{t_\alpha}$. (La raison pour cette notation double deviendra claire dans la section 2.)

**Lemme:** *Pour* $w\in W(M)$, *on a* ${}^wY_\alpha=Y_{w\alpha}$. *En particulier,* ${}^{s_\alpha}Y_\alpha=Y_\alpha^{-1}$.

*Preuve:* Posons $H_M(h_\alpha)=m_\alpha\alpha^\vee$ avec $m_\alpha>0$. Remarquons que $Y_\alpha=b_{h_\alpha}$ et donc ${}^wY_\alpha={}^wb_{h_\alpha}=b_{wh_\alpha w^{-1}}$. Comme $H_M(wh_\alpha w^{-1})=m_\alpha(w\alpha^\vee)$, il reste à prouver que $m_\alpha=m_{w\alpha}$. Comme, pour $\lambda\in a_M^*$,

$$\chi_{w\lambda}(wh_\alpha w^{-1})={}^w\chi_\lambda(wh_\alpha w^{-1})=\chi_\lambda(h_\alpha),$$

on a $\chi_{w\lambda}(wh_\alpha w^{-1})=1$ si et seulement, si $\chi_\lambda\in\mathfrak{X}^{\mathrm{nr}}(M_\alpha)$. Mais, ceci équivaut à $\chi_{w\lambda}={}^w\chi_\lambda\in\mathfrak{X}^{\mathrm{nr}}(M_{w\alpha})$. Donc, $H_M(h_{w\alpha})=H_M(wh_\alpha w^{-1})$, et, par suite, $m_\alpha=m_{w\alpha}$. □

**1.6** Le résultat suivant a été montré par A. Silberger [Si2, 1.6] (voir également la remarque dans la preuve de la proposition **4.1** dans [H2]). (Remarquons que l'énoncé du théorème **1.6** dans [Si2] est relatif à une fonction notée $U(\sigma, \cdot)$, mais celle-ci est effectivement la fonction $\mu$ (cf. dernier paragraphe sur la page 575 de [Si2].)

**Proposition:** *Soit $\alpha \in \Sigma_{red}(P)$ et $s = s_\alpha$. Si $\mu^{M_\alpha}$ n'est pas constante, alors $\alpha \in \Sigma_{\mathcal{O},\mu}$, et on peut trouver une constante $c'_s > 0$ et des nombres réels $a_s \geq 0$ et $b_s \geq 0$ tels que l'on ait l'identité de fonctions rationnelles*

$$\mu^{M_\alpha}(\sigma \otimes \cdot)$$
$$= c'_s \frac{(1 - X_\alpha(\cdot))(1 - X_\alpha^{-1}(\cdot))}{(1 - X_\alpha(\cdot)q^{-a_s})(1 - X_\alpha^{-1}(\cdot)q^{-a_s})} \frac{(1 + X_\alpha(\cdot))(1 + X_\alpha^{-1}(\cdot))}{(1 + X_\alpha(\cdot)q^{-b_s})(1 + X_\alpha^{-1}(\cdot)q^{-b_s})}.$$

**1.7 Remarque:** Comme $\Delta_{\mathcal{O},\mu}$ est linéairement indépendant, on peut choisir $\sigma$ tel que, pour tout $\alpha \in \Delta_{\mathcal{O},\mu}$, les nombres réels $a_{s_\alpha}$ et $b_{s_\alpha}$ de l'énoncé de la proposition ci-dessus vérifient $a_{s_\alpha} \geq b_{s_\alpha}$, ce que l'on supposera désormais. En particulier, $a_{s_\alpha} > 0$.

**1.8** Lorsque $P' = MU'$ est un autre sous-groupe parabolique, notons $J_{P|P'}$ l'opérateur d'entrelacement défini dans [W, IV]. C'est un opérateur rationnel

$$\mathfrak{X}^{\mathrm{nr}}(M) \to \mathrm{Hom}_{\mathbb{C}}(i_{P'\cap K}^K E, i_{P\cap K}^K E), \chi \mapsto J_{P|P'}(\sigma \otimes \chi)$$

qui induit en tout point régulier $\chi$ un homomorphisme entre les représentations $i_{P'}^G(\sigma \otimes \chi)$ et $i_P^G(\sigma \otimes \chi)$.

**Lemme:** *Soit $\alpha \in \Delta_{\mathcal{O},\mu}$, $s = s_\alpha$, et supposons que $M_\alpha$ soit un sous-groupe de Levi standard de $G$. Les coordonnées de l'opérateur $J_{P|sP}$ sont des fonctions rationnelles en $Y_\alpha$. Les pôles de $J_{P|sP}$ sont précisément les zéros de $\mu^{M_\alpha}$. Tout pôle est d'ordre 1 et son résidu est bijectif. Par ailleurs, $J_{P|sP}J_{sP|P}$ est égal à $(\mu^{M_\alpha})^{-1}$ à multiplication par une constante non nulle près.*

*Preuve:* La propriété de rationnalité résulte de [W, IV.1.1], en remarquant que $J_{P|sP}$ est invariant par $\mathfrak{X}^{\mathrm{nr}}(M_\alpha)$. Grâce à la propriété de fonctorialité que l'opérateur d'entrelacement vérifie relative à l'induction parabolique, on peut se ramener au cas où $P$ est un sous-groupe parabolique maximal. Alors, $sP = \overline{P}$, et $J_{P|sP}J_{sP|P}$ est égal à $\mu^{-1}$ multiplié par une constante non nulle.

Soit $\sigma_1$ un pôle de $J_{P|\overline{P}}$. Alors, par [Si1, 5.4.2.1], $\sigma_1$ est une représentation unitaire. Comme, d'après [W, V.2.3], l'opérateur $\mu J_{P|\overline{P}}$ est régulier en $\sigma_1$, on a bien $\mu(\sigma_1) = 0$.

Réciproquement, supposons $\mu(\sigma_1) = 0$. Alors, comme $\mu^{-1} = J_{P|\overline{P}} J_{\overline{P}|P}$, au moins un des opérateurs $J_{P|\overline{P}}$ ou $J_{\overline{P}|P}$ doit avoir un pôle en $\sigma_1$. Or, un tel pôle est au plus d'ordre 1 (cf. [W, IV.1.2]). Comme les pôles de $\mu^{-1}$ sont d'ordre 2, $J_{P|\overline{P}}$ doit avoir un pôle d'ordre 1.

Finalement, $\mu(\sigma) = 0$ implique par le théorème de Harish-Chandra **1.2** que $s\sigma \simeq \sigma$ et que $i_P^G\sigma$ et $i_{\overline{P}}^G\sigma$ sont irréductibles. Comme le résidu de $J_{P|\overline{P}}$ en $\sigma$ est un opérateur d'entrelacement $i_{\overline{P}}^G\sigma \to i_P^G\sigma$, celui-ci doit être bijectif. □

**1.9** Soient $w, w' \in W$ et notons $P_{ww'}$, $P_{w'}$ les sous-groupes paraboliques standards de $G$ de sous-groupes de Levi $ww'Mw'^{-1}w^{-1}$ et $w'Mw'^{-1}$ respectivement. Pour simplifier les notations, désignons, pour $w' \in W$ et $P'$ un sous-groupe parabolique semi-standard, par $w'^{-1}P'$ le sous-groupe parabolique semi-standard égal à $w'^{-1}\overline{P_{w'}}w'$. Nous écrivons dans la suite $\mu_{P,w,w'}$ ou plus simplement $\mu_{w,w'}$ pour la fonction rationnelle sur $\mathfrak{X}^{\mathrm{nr}}(M)$ telle que $\mu_{w,w'}(\chi) = \prod_\alpha \mu^{M_\alpha}(\sigma\otimes\chi)$, le produit portant sur $\Sigma(P)\cap\Sigma(w'^{-1}\overline{P_{w'}})\cap\Sigma(w'^{-1}w^{-1}P_{ww'})$. (Remarquons que $P_{w'} = P_{ww'} = P$, si $w, w' \in W(M, \mathcal{O})$.)

**Proposition:** *Soient $w, w' \in W$. Alors, en tant que fonction rationnelle en $\chi \in \mathfrak{X}^{\mathrm{nr}}(M)$,*

$$\begin{aligned}&\lambda(ww')J_{w'^{-1}w^{-1}P_{ww'}|P}(\sigma\otimes\chi)\\ =&\mu_{w,w'}(\chi)\lambda(w)J_{w^{-1}P_{ww'}|P_{w'}}(w'\sigma\otimes w'\chi)\lambda(w')J_{w'^{-1}P_{w'}|P}(\sigma\otimes\chi).\end{aligned}$$

*Preuve:* Ceci est une conséquence immédiate des règles de composition pour les opérateurs d'entrelacement [W, IV.3]. □

**1.10 Proposition:** *Soit $P' = MU'$ un autre sous-groupe parabolique de Levi $M$. Si $\Sigma(P)\cap\Sigma(\overline{P'})$ a une intersection vide avec $\Sigma_{\mathcal{O},\mu}(P)$, l'opérateur $J_{P'|P}$ est bien défini et bijectif en tout point de $\mathcal{O}$.*

*Preuve:* Comme $J_{P'|P}$ se décompose en des opérateurs d'entrelacement élémentaires qui proviennent d'opérateurs d'entrelacement relatifs à des $M_\alpha$, $\alpha \notin \Sigma_{\mathcal{O},\mu}(P)$, on est ramené au cas où $P$ est un sous-groupe parabolique maximal de $G$ et $P' = \overline{P}$ avec $\mu$ constante. D'après [Si1, 5.4.2.1], les pôles de $J_{\overline{P}|P}$ sont des représentations qui deviennent unitaires après torsion par un caractères non ramifiés de $G$. Comme $\mu$ est constante et que l'opérateur $\mu J_{\overline{P}|P}$ est régulier en toute représentation unitaire de $\mathcal{O}$ [W, V.2.3], $J_{\overline{P}|P}$ doit lui-même être régulier sur tout $\mathcal{O}$, l'opérateur $J_{\overline{P}|P}$ étant polynomiale sur $\mathfrak{X}^{\mathrm{nr}}(G)$. (Cela se voit directement sur la définition de $J_{\overline{P}|P}$.) Il est bijectif en tout point de $\mathcal{O}$, puisque $\mu$ ne s'annulle pas. □

**1.11** *Définition:* On pose $R(\mathcal{O}) = \{w \in W(M,\mathcal{O}) | w\alpha \in \Sigma(P)$ pour tout $\alpha \in \Sigma(P)\cap\Sigma_{\mathcal{O},\mu}\}$.

**1.12 Proposition:** *Le groupe $R(\mathcal{O})$ est un sous-groupe de $W(M,\mathcal{O})$. On a $W(M,\mathcal{O}) = R(\mathcal{O}) \ltimes W_{\mathcal{O}}$.*

*Preuve:* Posons $R = R(\mathcal{O})$. Comme $\Sigma_{\mathcal{O},\mu}$ est $W(M,\mathcal{O})$-invariant (cf. preuve de **1.3**), $R$ est par définition un sous-groupe de $W(M,\mathcal{O})$. Le groupe $W_{\mathcal{O}}$ est engendré par les symétries $s_\alpha$, $\alpha \in \Delta_{\mathcal{O}}$. Comme $ws_\alpha w^{-1} = s_{w\alpha}$ et que $w\alpha \in \Sigma_{\mathcal{O},\mu}$, $W_{\mathcal{O}}$ est un sous-groupe distingué de $W(M,\mathcal{O})$. Comme $W_{\mathcal{O}}$ permute les chambres de Weyl dans $\Sigma_{\mathcal{O},\mu}$, on a $W_{\mathcal{O}} \cap R = \{1\}$. Il reste à montrer que $W(M,\mathcal{O}) = RW_{\mathcal{O}}$.

Or, soit $w \in W(M,\mathcal{O})$. Alors, $w(\Sigma(P) \cap \Sigma_{\mathcal{O},\mu}) = \Sigma(wPw^{-1}) \cap \Sigma_{\mathcal{O},\mu}$. Ceci est l'ensemble des racines positives dans $\Sigma_{\mathcal{O},\mu}$ pour un certain ordre sur $\Sigma_{\mathcal{O},\mu}$. Comme $W_{\mathcal{O}}$ permute transitivement les différents ordres sur $\Sigma_{\mathcal{O}}$, il existe $w' \in W_{\mathcal{O}}$ tel que $w'w(\Sigma(P) \cap \Sigma_{\mathcal{O},\mu}) = \Sigma(P) \cap \Sigma_{\mathcal{O},\mu}$, i.e. $w'w \in R$. □

**1.13 Proposition:** *Supposons que $G$ soit un groupe symplectique ou orthogonal, et notons $d = rg_F(G)$.*

*a) En remplaçant $\sigma$ par un autre élément de $\mathcal{O}$ et en conjuguant $(M,\sigma)$ par un élément de $G$, on peut supposer $M = \underline{M}(F)$ avec*

$$\underline{M} = \mathrm{GL}_{k_1} \times \cdots \times \mathrm{GL}_{k_1} \times \mathrm{GL}_{k_2} \times \cdots \times \mathrm{GL}_{k_2} \times \cdots \times \mathrm{GL}_{k_r} \times \cdots \times \mathrm{GL}_{k_r} \times \underline{H}_k,$$

*où $\underline{H}_k$ désigne un groupe semi-simple de rang absolu $k$ du même type que $G$, et*

$$\sigma = \sigma_1 \otimes \cdots \otimes \sigma_1 \otimes \sigma_2 \otimes \cdots \otimes \sigma_2 \otimes \cdots \otimes \sigma_r \otimes \cdots \otimes \sigma_r \otimes \tau,$$

*les classes inertielles des $\sigma_i$ étant deux à deux distinctes, ainsi que $\sigma_i \simeq \sigma_i^\vee$ si $\sigma_i$ et $\sigma_i^\vee$ sont dans une même orbite inertielle.*

*b) Notons $d_i$ le nombre de facteurs égaux à $\sigma_i$ et identifions $A_{\underline{M}}$ à $\mathbb{T} = \mathbb{G}_m^{d_1} \times \mathbb{G}_m^{d_2} \times \cdots \times \mathbb{G}_m^{d_r}$. Notons $\alpha_{i,j}$ le caractère rationnel de $A_{\underline{M}}$ (identifié à $\mathbb{T}$) qui envoie un élément $x = (x_{1,1}, \dots, x_{1,d_1}, x_{2,1}, \dots, x_{2,d_2}, \dots, x_{r,1}, \dots, x_{r,d_r})$ sur $x_{i,j}x_{i,j+1}^{-1}$, si $j < d_i$, et sur $x_{i,d_i}$, si $j = d_i$. Le système de racines $\Sigma_{\mathcal{O},\mu}$ est la somme directe de $r$ composantes irréductibles ou vides $\Sigma_{\mathcal{O},\mu,i}$, $i = 1, \dots, r$, définies de la manière suivante:*

*Supposons d'abord ou $k \neq 0$ ou $G$ de système de racine de type $B_d$.*

*(i) Si la fonction $s \mapsto \mu(\sigma_i|\det_{k_i}|^s \otimes \tau)$ (définie relativement à $\mathrm{GL}_{k_i} \times \underline{H}_k$ et $\underline{H}_{k_i+k}$) a un pôle sur $\mathbb{C}$, alors une base de $\Sigma_{\mathcal{O},\mu,i}$ est donnée par $\{\alpha_{i,1}, \dots, \alpha_{i,d_i}\}$, et ce système est de type $B_{d_i}$.*

*(ii) Si la fonction $s \mapsto \mu(\sigma_i|\det_{k_i}|^s \otimes \tau)$ (définie relativement à $\mathrm{GL}_{k_i} \times \underline{H}_k$ et $\underline{H}_{k_i+k}$) est régulière sur $\mathbb{C}$ et $\sigma_i \simeq \sigma_i^\vee$, alors une base de $\Sigma_{\mathcal{O},\mu,i}$ est donnée par $\{\alpha_{i,1}, \dots, \alpha_{i,d_i-1}, \alpha_{i,d_i-1} + 2\alpha_{i,d_i}\}$, et ce système est de type $D_{d_i}$.*

*(iii) Sinon, une base de $\Sigma_{\mathcal{O},\mu,i}$ est donnée par $\{\alpha_{i,1}, \dots, \alpha_{i,d_i-1}\}$, et ce système est de type $A_{d_i-1}$.*

*Supposons maintenant $k=0$ et ou $G$ de système de racines de type $C_d$ ou $G$ de système de racines de type $D_d$ et $k_i\geq 2$. Alors, dans le cas (i) ci-dessus, une base de $\Sigma_{\mathcal{O},\mu,i}$ est donnée par $\{\alpha_{i,1},\dots,2\alpha_{i,d_i}\}$, et le système est de type $C_{d_i}$. Dans les autres cas, la situation reste inchangée.*

*Supposons finalement $k=0$, $G$ de système de racines de type $D_d$ et $k_i=1$. Alors, si $\sigma_i\simeq\sigma_i^\vee$, une base de $\Sigma_{\mathcal{O},\mu,i}$ est donnée par $\{\alpha_{i,1},\dots,\alpha_{i,d_i-1},\alpha_{i,d_i-1}+2\alpha_{i,d_i}\}$ et ce système est de type $D_{d_i}$. Sinon, le système est de type $A_{d_i-1}$ et une base est donnée par $\{\alpha_{i,1},\dots,\alpha_{i,d_i-1}\}$.*

*Preuve:* Les sous-groupes de Levi des groupes classiques sont décrits dans [G, Kap. 5]. Par permutation des facteurs, on met $\sigma$ dans la forme voulue. Supposons qu'il existe un caractère non ramifié $\chi_i$ de $\mathrm{GL}_{k_i}$ tel que $\sigma_i^\vee\simeq\sigma_i\otimes\chi_i$. On peut écrire $\chi_i=|\det_{k_i}|_F^{s_i}$, où $s_i$ est un nombre complexe. Posons $\chi_i^{1/2}=|\det_{k_i}|^{s_i/2}$. Alors $(\sigma_i\otimes\chi_i^{1/2})^\vee\simeq\sigma_i^\vee\otimes\chi_i^{-1/2}\chi_i=\sigma_i\otimes\chi_i^{1/2}$.

Identifions $A_{\underline{M}}$ au tore $\mathbb{T}$, écrivons $x=(x_{1,1},\dots,x_{1,d_1},x_{2,1},\dots,x_{2,d_2},\dots,x_{r,1},\dots,x_{r,d_r})$ pour l'élément général de $\mathbb{T}$ et examinons les différents cas de figure.

Supposons d'abord ou $k\neq 0$ ou $G$ de système de racines de type $B_d$. Les racines réduites dans $\Sigma(A_M)$ sont alors $x\mapsto x_{i,j}^{\pm1}x_{i',j'}^{\pm1}$, $(i,j)\neq(i',j')$, ainsi que $x\mapsto x_{i,j}^{\pm1}$. (Comme $k\neq 0$, il existe bien dans le cas d'un système de racines de type $C_d$ ou $D_d$ des racines de restriction à $A_M$ égale à $x\mapsto x_{i,j}^{\pm1}$.) Si $\alpha\in\Sigma(A_M)$ correspond à une racine $x\mapsto x_{i,j}x_{i',j'}^{-1}$, $(i,j)\neq(i',j')$, alors la fonction $\mu^{M_\alpha}$ est égale à celle définie à partir de la représentation $\sigma_i\otimes\sigma_{i'}$ du sous-groupe de Levi $\mathrm{GL}_{k_i}\times\mathrm{GL}_{k_{i'}}$ de $\mathrm{GL}_{k_i+k_{i'}}$. Si $\alpha\in\Sigma(A_M)$ correspond à une racine $x\mapsto x_{i,j}x_{i',j'}$, $(i,j)\neq(i',j')$, $i\leq i'$, alors la fonction $\mu^{M_\alpha}$ est égale à celle définie à partir de la représentation $\sigma_i\otimes\sigma_{i'}^\vee$ du sous-groupe de Levi $\mathrm{GL}_{k_i}\times\mathrm{GL}_{k_{i'}}$ de $\mathrm{GL}_{k_i+k_{i'}}$. Dans les autres cas, la fonction $\mu^{M_\alpha}$ est celle définie par la représentation $\sigma_i\otimes\tau$ du sous-groupe de Levi $\mathrm{GL}_{k_i}\times H_m$ de $H_{m+k_i}$. Selon les différents cas présentés dans (i), (ii) et (iii), on déduit alors du théorème de Harish-Chandra **1.2** et des résultats de Bernstein-Zelevinsky [BZ2, 4.2] le résultat indiqué.

Pour $k=0$ et $G$ de système de racines de type $C_d$, les racines réduites dans $\Sigma(A_M)$ sont $x\mapsto x_{i,j}^{\pm1}x_{i',j'}^{\pm1}$, $(i,j)\neq(i',j')$, ainsi que $x\mapsto x_{i,j}^{\pm2}$. On conclut alors comme ci-dessus.

Pour $k=0$ et $G$ de système de racines de type $D_d$ finalement, les racines réduites dans $\Sigma(A_M)$ sont $x\mapsto x_{i,j}^{\pm1}x_{i',j'}^{\pm1}$, $(i,j)\neq(i',j')$, ainsi que, si $k_i\geq 2$, $x\mapsto x_{i,j}^{\pm2}$. On conclut alors comme ci-dessus, en tenant compte du fait que la fonction $\mu$ est bien connue pour les groupes déployés de rang 1. □

**1.14 Proposition:** *Supposons que $G$ soit le groupe des points rationnels d'une forme intérieure de $\mathrm{GL}_n$, i.e. $G=\mathrm{GL}_n(D)$, où $D$ est une algèbre à division de centre $F$.*

*En remplaçant $\sigma$ par un autre élément de $\mathcal{O}$ et en conjuguant $(M,\sigma)$ par un*

*élément de $G$, on peut supposer $M$ égal à*

$$\mathrm{GL}_{k_1}(D)\times\cdots\times\mathrm{GL}_{k_1}(D)\times\mathrm{GL}_{k_2}(D)\times\cdots\times\mathrm{GL}_{k_2}(D)\times\cdots\times\mathrm{GL}_{k_r}(D)\times\cdots\times\mathrm{GL}_{k_r}(D)$$

*et*

$$\sigma=\sigma_1\otimes\cdots\otimes\sigma_1\otimes\sigma_2\otimes\cdots\otimes\sigma_2\otimes\cdots\otimes\sigma_r\otimes\cdots\otimes\sigma_r,$$

*les classes inertielles des $\sigma_i$ étant deux à deux distinctes.*

*Notons $d_i$ le nombre de facteurs égaux à $\sigma_i$ et identifions $A_{\underline{M}}$ à $\mathbb{T}=\mathbb{G}_m^{d_1}\times\mathbb{G}_m^{d_2}\times\cdots\times\mathbb{G}_m^{d_r}$. Notons $\alpha_{i,j}$ le caractère rationnel de $A_{\underline{M}}$ (identifié à $\mathbb{T}$) qui envoie un élément $x=(x_{1,1},\dots,x_{1,d_1},x_{2,1},\dots,\ x_{2,d_2},\ \dots,x_{r,1},\dots\ ,x_{r,d_r})$ sur $x_{i,j}x_{i,j+1}^{-1}$. Le système des racines $\Sigma_{\mathcal{O},\mu}$ est la somme directe de $r$ composantes irréductibles ou vides $\Sigma_{\mathcal{O},\mu,i}$, $i=1,\dots,r$, de type $A_{d_i-1}$ et de base $\{\alpha_{i,1},\dots,\alpha_{i,d_i-1}\}$ respectivement.*

*Preuve:* Pour $G$ déployé, ceci résulte des travaux de Bernstein-Zelevinsky [BZ2, 4.2] avec le théorème de Harish-Chandra **1.2**. Dans le cas général, on le déduit des travaux de Bernstein-Zelevinsky avec la formule des traces [DKV, T]. □

**1.15 Proposition:** *(i) Si $G$ est un groupe linéaire ou le groupe multiplicatif d'une algèbre à division, alors $R(\mathcal{O})=1$.*

*(ii) Sinon, $R(\mathcal{O})\neq 1$, si et seulement si les conclusions du (ii) de la proposition* **1.13** *sont vérifiées pour au moins un indice $i$ avec, si $G$ est de type $D_n$, de plus $k_i$ pair ou bien $H\neq 1$ et $\tau$ invariant par l'automorphisme extérieur. Le groupe $R(\mathcal{O})$ est alors le produit direct des groupes $R(\mathcal{O})_i$ indexés par les $i$ pour lesquels les conclusions ci-dessus sont vérifiées avec $R(\mathcal{O})_i=\{1,s_{\alpha_{i,d_i}}\}$ ($s_{\alpha_{i,d_i}}$ échangeant $\alpha_{i,d_{i-1}}$ et $\alpha_{i,d_{i-1}}+2\alpha_{i,d_i}$ et vérifiant $s_{\alpha_{i,d_i}}\sigma_i\simeq\sigma_i^\vee$).*

*Preuve:* Il est facile de voir que les conditions du théorème **1.2** sont réunies exactement dans les cas décrits dans l'énoncé, en prenant en compte [BJ, 3.4] pour les groupes de type $D_n$. □

La proposition **1.15** ainsi que la proposition **1.12** nous conduisent à ajouter les deux hypothèses suivantes sur le choix des représentants des éléments de $W(M,\mathcal{O})$: concernant les représentants des éléments dans $R(\mathcal{O})$, on suppose que, si $r$ est un produit d'éléments $s_{\alpha_i}$ dans $R(\mathcal{O})_i$ deux à deux distincts dans $R(\mathcal{O})$, alors il en est ainsi pour les représentants. Concernant un élément arbitraire $w$ de $W(M,\mathcal{O})$, on suppose que, si $w=w_{\mathcal{O}}r$ avec $w_{\mathcal{O}}\in W_{\mathcal{O}}$ et $r\in R(\mathcal{O})$, alors le représentant de $w$ est le produit de celui de $w_{\mathcal{O}}$ avec celui de $r$.

**1.16** La proposition suivante est un cas particulier d'un résultat bien connu pour les groupes finis [He, Lemme 7.5]. Remarquons qu'elle n'est probablement pas valable pour un groupe réductif arbitraire sur $F$, quoique la preuve reste valable dans d'autres cas que celui d'un groupe orthogonal ou symplectique ou d'une forme intérieure d'un groupe linéaire général considérés ici.

**Proposition:** *Soit $(\sigma_1, E_1)$ une composante irréductible de $\sigma_{|M^1}$. Alors $\sigma_1$ se prolonge en une représentation $\sigma_2$ de $M^\sigma = \{m \in M| \ ^m\sigma_1 \simeq \sigma_1\}$ telle que $\sigma_{|M^\sigma} = \bigoplus_{m\in M/M^\sigma} \ ^m\sigma_2$. En particulier, $\sigma = \mathrm{Ind}_{M^\sigma}^M \sigma_2$ et $^m\sigma_2 \not\simeq \sigma_2$ si $m \notin M^\sigma$. Par ailleurs, $\sigma_{|M^1}$ est somme directe de représentations irréductibles deux à deux non isomorphes.*

*Preuve:* Nous allons d'abord supposer $A_M$ de rang 1. Remarquons d'abord que, si $H$ est un sous-groupe ouvert fermé d'indice fini d'un groupe localement profini $G$, alors $\mathrm{ind}_H^G = \mathrm{Ind}_H^G$ et, comme dans le cas des groupes finis, une représentation $(\sigma, E)$ est induite par une représentation $(\sigma_2, E_1)$, si et seulement si $V = \bigoplus_{g\in G/H} \sigma(g)E_1$.

Observons maintenant que $A_M \subseteq M^\sigma$. Le quotient $M/M^\sigma$ est donc un groupe commutatif fini. Comme $\sigma$ est cuspidale, on peut écrire $E = E_{\sigma_1} \oplus E'$ de sorte que $E_{\sigma_1}$ soit somme directe de représentations isomorphes à $\sigma_1$ et qu'aucun sous-quotient de $E'$ ne soit isomorphe à $\sigma_1$ [BZ1, 2.44]. Le groupe $M^\sigma$ agit donc dans $E_{\sigma_1}$, et on a $E = \bigoplus_{m\in M/M^\sigma} \sigma(m)E_{\sigma_1}$. Par suite, par la remarque ci-dessus, $\mathrm{Ind}_{M^\sigma}^M E_{\sigma_1} = E$. Comme $E$ est irréductible, la représentation de $M^\sigma$ dans $E_{\sigma_1}$ doit être irréductible. On va maintenant montrer qu'en fait $E_{\sigma_1} = E_1$.

Voyons d'abord que $\sigma_1$ se prolonge en une représentation irréductible $\sigma_2$ de $M^\sigma$. D'abord, $\sigma_1$ se prolonge en une représentation irréductible $\widetilde{\sigma}_1$ de $A_M M$, donnée par $am \mapsto \chi_\sigma(a)\sigma_1(m)$, $\chi_\sigma$ désignant le caractère central de $\sigma$. Comme le quotient $M^\sigma/A_M M^1$ est par hypothèse un groupe cyclique fini, il est alors facile de définir un prolongement à $M^\sigma$. (Étant donné un élément $m'$ de $M^\sigma$ dont l'image engendre le groupe quotient cyclique $M^\sigma/A_M M^1$, on peut toujours trouver un isomorphisme $A_M M^1$-équivariant $\varphi_{m'} : E_1 \to {m'}^{-1}E_1$ tel que ${m'}^i m \mapsto \varphi_{m'}^i \widetilde{\sigma}_1(m)$ pour $m \in A_M M^1$ donne le prolongement voulu.)

Montrons maintenant que $\mathrm{Ind}_{A_M M^1}^{M^\sigma} \widetilde{\sigma}_1 = \bigoplus_{\chi\in X(M^\sigma/A_M M^1)} \sigma_2 \otimes \chi$. Les représentations $\sigma_2 \otimes \chi$, $\chi \in X(M^\sigma/A_M M^1)$ sont deux à deux non isomorphes, puisque tout isomorphisme $\sigma_2 \otimes \chi \mapsto \sigma_2 \otimes \chi'$ induit un automorphisme de $\widetilde{\sigma}_1$ et, par le lemme de Schur, il serait donc scalaire, ce qui est impossible. Par la réciprocité de Frobenius, les sous-représentations irréductibles de $\mathrm{Ind}_{A_M M^1}^{M^\sigma} \widetilde{\sigma}_1$ sont les représentations irréductibles de $M^\sigma$ de restriction à $A_M M^1$ égale à $\widetilde{\sigma}_1$ et leur multiplicité est 1. Or, les arguments de la preuve de [BZ1, 3.29] montrent que toute représentation irréductible de $M^\sigma$ de restriction à $A_M M^1$ égale à $\widetilde{\sigma}_1$ est isomorphe à $\sigma_2 \otimes \chi$ pour un $\chi \in X(M^\sigma/A_M M^1)$. Comme $Ind_{A_M M^1}^{M^\sigma} \widetilde{\sigma}_1$ est unitaire, puisque $\widetilde{\sigma}_1$ l'est, on a bien la décomposition indiquée.

Il s'ensuit que $E_{\sigma_1} = E_1$: la représentation $\rho$ de $M^\sigma$ dans $E_{\sigma_1}$ est irréductible, de restriction égale à un multiple de $\sigma_1$. La réciprocité de Frobenius donne donc $\mathrm{Hom}_M(\rho, \mathrm{Ind}_{A_M M^1}^{M_\sigma} \widetilde{\sigma}_1) \neq 0$, ce qui implique que $\rho$ est isomorphe à $\sigma_2 \otimes \chi$ pour un $\chi \in X(M^\sigma/A_M M^1)$. Par suite, $E_{\sigma_1} = E_1$ et on peut supposer dans la suite que $\rho = \sigma_2$.

Il en résulte, en remarquant que $M^\sigma$ est distingué dans $M$, que $\sigma = \mathrm{Ind}_{M^\sigma}^M \sigma_2$ et $\sigma_{|M^\sigma} = \bigoplus_{m \in M/M^\sigma} {}^m\sigma_2$. Comme $\sigma$ est irréductible, les ${}^m\sigma_2$ sont deux à deux non isomorphes. Comme par ailleurs $({}^m\sigma_2)_{|M^1} = {}^m\sigma_1$ avec $m$ parcourant un système de représentants de $M/M^\sigma$, $\sigma_{|M^1}$ est somme directe de représentations irréductibles deux à deux non isomorphes.

Supprimons maintenant l'hypothèse que $A_M$ est de rang 1. Comme $G$ est un groupe symplectique ou orthogonal, ou le groupe multiplicatif d'une algèbre simple, $M$ est de la forme $\mathrm{GL}_{m_1}(D) \times \cdots \times \mathrm{GL}_{m_r}(D) \times H$, où $D$ est une algèbre à division de centre $F$ et où $H$ est un groupe semi-simple du même type que $G$ ou le groupe trivial. On en déduit que $M^1$ est de la forme $\mathrm{GL}_{m_1}(D)^1 \times \cdots \times \mathrm{GL}_{m_r}(D)^1 \times H$, que $\sigma$ est isomorphe à un produit $\sigma_1 \otimes \cdots \otimes \sigma_r \otimes \tau$ avec $\sigma_i$, $\tau$ des représentations cuspidales, et que $M^\sigma$ est de la forme $\mathrm{GL}_{m_1}(D)^{\sigma_1} \times \cdots \times \mathrm{GL}_{m_r}(D)^{\sigma_r} \times H$. On est donc ramené au cas où $M = \mathrm{GL}_m(D)$. Mais, alors $A_M$ est de rang 1, et l'assertion est vraie par ce qui précédait. □

**1.17** Lorsque $\chi$ est un caractère non ramifié de $M$, notons $E_\chi$ l'espace $E$ muni de la représentation $\sigma \otimes \chi$.

**Corollaire:** (avec les notations de la proposition **1.16**) *On a* $\mathrm{Stab}(\mathcal{O}) = X(M/M^\sigma)$ *et, pour tout* $\chi \in X(M/M^\sigma)$*, l'application* $\phi_{\sigma,\chi} : E \to E_\chi$ *qui envoie un élément* $e \in \sigma(m)E_1$ *sur* $\chi(m)e$ *est un isomorphisme.*

*Preuve:* Soit $\chi \in \mathfrak{X}^{\mathrm{nr}}(M)$ tel que $\sigma \otimes \chi \simeq \sigma$. Alors, il existe $m \in M$ tel que $\sigma_2 \otimes \chi \simeq {}^m\sigma_2$. Par suite, $\sigma_1 \simeq (\sigma_2 \otimes \chi)_{|M^1} \simeq ({}^m\sigma_2)_{|M^1} = {}^m\sigma_1$, i.e. on a $m \in M^\sigma$ et donc $\sigma_2 \otimes \chi \simeq \sigma_2$. Soit $\varphi_\chi$ un automorphisme de $E_1$ qui entrelace $\sigma_2 \otimes \chi$ et $\sigma_2$. Alors il entrelace $\sigma_1$ avec lui-même, i.e. $\varphi_\chi$ est la multiplication par un scalaire. Il s'ensuit que $\chi_{|M^\sigma} \equiv 1$. Donc $Stab(\mathcal{O}) \subseteq X(M/M^\sigma)$. Inversement, c'est une vérification immédiate que pour tout $\chi \in X(M/M^\sigma)$ les applications $\phi_{\sigma,\chi}$ définissent des isomorphismes entre $\sigma$ et $\sigma \otimes \chi$. □

**2.** Notons $B = B_M$ l'anneau des polynômes sur la variété affine complexe $\mathfrak{X}^{\mathrm{nr}}(M)$. C'est un anneau intègre, et on désignera par $K(B)$ son corps des fractions.

Écrivons $E_B$ (resp. $E_{K(B)}$) pour l'espace vectoriel complexe $E \otimes_{\mathbb{C}} B$ (resp. $E_B \otimes_B K(B) \simeq E \otimes_{\mathbb{C}} K(B)$). Munissons-le de l'action de $M$ donnée par $\sigma_B : M \to Aut_{\mathbb{C}}(E_{K(B)})$, $\sigma_B(m)(e \otimes b) := \sigma(m)e \otimes bb_m$. Ici, $b_m : \mathfrak{X}^{\mathrm{nr}}(M) \to \mathbb{C}$ est défini par $b_m(\chi) = \chi(m)$.

Par restriction à $K$, l'espace $i_P^G E_B$ (resp. $i_P^G E_{K(B)}$) est isomorphe à $i_{K\cap P}^K E \otimes_{\mathbb{C}} B$ (resp. $i_{K\cap P}^K E \otimes_{\mathbb{C}} K(B)$). On identifie ainsi $i_{K\cap P}^K E$ à son image dans $i_P^G E_B$ (resp. $i_P^G E_{K(B)}$) par l'homomorphisme canonique $v \mapsto v \otimes 1$.

**2.1** Remarquons d'abord que $E_B$ est canoniquement isomorphe à $\mathrm{ind}_{M^1}^M E$. En effet, notons $\mathcal{R}(M/M^1)$ un système de représentants de $M/M^1$. On obtient un isomorphisme canonique entre $B$ et $\mathrm{ind}_{M^1}^M \mathbb{C} = \mathbb{C}[M/M^1]$, en envoyant un élément de la forme $b_m$, $m \in \mathcal{R}(M/M^1)$, sur l'application $M \to \mathbb{C}$ qui vaut 1 sur $m^{-1}M^1$ et 0 en dehors de cet ensemble. Identifions ces deux espaces à l'aide de cet isomorphisme. On a un isomorphisme canonique $E\otimes \mathrm{ind}_{M^1}^M \mathbb{C} \to \mathrm{ind}_{M^1}^M E_{|M^1}$, en envoyant $e\otimes b$ sur l'application $v_{e\otimes b}$ dans $\mathrm{ind}_{M^1}^M E_{|M^1}$ qui envoie $m$ sur $b(m)\sigma(m)e$. L'isomorphisme réciproque est donné par $\mathrm{ind}_{M^1}^M E_{|M^1} \to E \otimes \mathrm{ind}_{M^1}^M \mathbb{C}$, $v \mapsto \sum_{m\in\mathcal{R}(M/M^1)} \sigma(m^{-1})$ $v(m) \otimes b_{m^{-1}}$. (Ici, $b_{m^{-1}}$ est considéré comme élément de $\mathrm{ind}_{M^1}^M \mathbb{C}$ au moyen de l'isomorphisme ci-dessus.)

**2.2** Pour $w \in W^G$, notons $wE_B$ (resp. $wE_{K(B)}$) l'espace $E_B$ (resp. $E_{K(B)}$) muni de la représentation $w\sigma_B$ de $wMw^{-1}$. (On a donc $w\sigma_B(wmw^{-1}) = \sigma_B(m)$.) Notons $\tau_w$ l'automorphisme de $\mathbb{C}$-espace vectoriel de $E\otimes_{\mathbb{C}} K(B)$ qui envoie $e\otimes b$ sur $e\otimes\ {}^wb$, ${}^wb(\chi) := b(w^{-1}\chi)$. (On a donc ${}^wb_m = b_{wmw^{-1}}$ pour $m \in M$.) Il définit un isomorphisme $M$-équivariant entre $wE_{K(B)}$ et $(wE)_{K(B)}$ et l'on notera encore $\tau_w$ l'isomorphisme $G$-équivariant $i_P^G w(E_{K(B)}) \to i_P^G\ (wE)_{K(B)}$ qui en est déduit par fonctorialité.

**2.3** Pour $b \in K(B)$, notons $b_\chi$ l'élément de $K(B)$ donné par $b_\chi(\chi') = b(\chi\chi')$. On a $(b_m)_\chi = \chi(m)b_m$. Le résultat suivant est une vérification directe:

**Lemme:** *Soit $\chi \in \mathfrak{X}^{\mathrm{nr}}(M)$ et notons $E_\chi$ l'espace $E$ muni de la représentation $\sigma\otimes\chi$ de $M$. L'application $\rho_\chi : E\otimes_{\mathbb{C}} B \to E_\chi \otimes_{\mathbb{C}} B$, $e\otimes b \mapsto e\otimes b_\chi$, définit un isomorphisme entre ces représentations lisses de $M$. On a $\rho_\chi^{-1} = \rho_{\chi^{-1}}$ et $\tau_w \circ \rho_\chi = \rho_{w\chi}\tau_w$.*

**2.4** Nous allons définir pour tout $w \in W_{\mathcal{O}}$ un isomorphisme $\rho_w : i_P^G wE \to i_P^G E$ qui envoie $i_P^G(wE)_B$ sur $i_P^G E_B$. Pour cela, posons

$$\rho_w = [((\prod_{\alpha\in\Sigma_{\mathcal{O},\mu}(P)\cap w^{-1}\Sigma_{\mathcal{O},\mu}(\overline{P})} (Y_\alpha(\chi)-1))\lambda(w)J_{w^{-1}P|P}(\sigma\otimes\chi))_{|\chi=1}]^{-1},$$

si cette expression est bien définie. On écrira $\rho_{\sigma,w}$, si on veut souligner la dépendance de $\sigma$. Il est par ailleurs clair que les opérateurs $\rho_{\sigma,w}$ et $\tau_w$ commutent.

**Proposition:** *L'application $\rho_w$ provient d'un isomorphisme $wE \simeq E$, et il définit pour tout $w \in W_{\mathcal{O}}$ un isomorphisme $i_P^G wE \to i_P^G E$. Il vérifie les propriétés de compatibilité suivantes: pour tout $w' \in W$ avec $w'Mw'^{-1}$ égal au sous-groupe de Levi standard d'un sous-groupe parabolique standard $P_{w'}$ de $G$, $\rho_w$ induit un iso-* *morphisme* $i^G_{P_{w'}}(w'wE)_{K(B)} \to i^G_{P_{w'}}(w'E)_{K(B)}$, *et on a la formule de commutation*

$$\lambda(w')J_{w'^{-1}P_{w'}|P}(\sigma\otimes\cdot)\rho_w = \rho_w\lambda(w')J_{w'^{-1}P_{w'}|P}(w\sigma\otimes\cdot).$$

*Par ailleurs, pour* $w'\in W_{\mathcal{O}}$, *on a, avec* $w''$ *égal au représentant de* $ww'$ *dans* $K$,

$$\rho_w\rho_{w'} = (\mu_{w,w'}\prod_\alpha (Y_\alpha-1)^{-1}(Y_\alpha^{-1}-1)^{-1})_{|\chi=1}\lambda(w''(ww')^{-1})\rho_{ww'},$$

*le produit portant sur l'ensemble* $\Sigma_{\mathcal{O},\mu}(P)\cap {w'}^{-1}\Sigma_{\mathcal{O},\mu}(\overline{P})\cap {w'}^{-1}w^{-1}\Sigma_{\mathcal{O},\mu}(P)$, *ainsi que, pour* $r\in R(\mathcal{O})$ *et* $s_\alpha$ *une réflexion simple dans* $W_{\mathcal{O}}$,

$$\rho_{r\sigma,s_{r\alpha}} = \lambda(rs_\alpha r^{-1}s_{r\alpha}^{-1})\rho_{\sigma,s_\alpha}.$$

*Preuve:* Supposons d'abord $w=s_\alpha$ avec $\alpha\in\Delta_{\mathcal{O},\mu}$. Alors, $\mu^{M_\alpha}(\sigma)=0$, et il résulte du théorème de Harish-Chandra **1.2** que $i^{M_\alpha}_{P\cap M_\alpha}\sigma$ est irréductible et que $s_\alpha\sigma\simeq\sigma$.

Si $M_\alpha$ est standard, on déduit de **1.8** que l'opérateur $J_{s_\alpha P|P}$ considéré comme fonction rationnelle en $Y_\alpha$ a un pôle d'ordre 1 en $Y_\alpha=1$ et que l'opérateur $Res_{Y_\alpha=1}(\lambda(s_\alpha)J_{s_\alpha P|P}(\sigma\otimes(\cdot)) : i^G_P\sigma \to i^G_P s_\alpha\sigma$ est bijectif. Par suite, $\rho_{s_\alpha}$ est bien défini et bijectif. Il est induit par fonctorialité par un isomorphisme entre les représentations irréductibles $i^{M_\alpha}_{P\cap M_\alpha}s_\alpha\sigma$ et $i^{M_\alpha}_{P\cap M_\alpha}\sigma$. Or, un tel isomorphisme est, par le lemme de Schur, uniquement déterminé à une constante près. En conséquence, $\rho_{s_\alpha}$ est lui-même induit par fonctorialité par un isomorphisme $s_\alpha\sigma\to\sigma$. Les deux propriétés de compatibilité en résultent par la propriété de fonctorialité pour les opérateurs d'entrelacement.

Si $M_\alpha$ n'est pas standard, alors ou $\alpha$ est contenue dans une composante irréductible de type $D_{d_i}$ de $\Sigma_{\mathcal{O},\mu}$, et il existe par **1.15** et par **1.13** un élément $r\in R(\mathcal{O})$ tel que $\alpha=r\beta$ avec $\beta\in\Delta_{\mathcal{O},\mu}$, $M_\beta$ standard et $r\sigma\simeq\sigma$, ou bien cette composante irréductible est de type $B_{d_i}$ ou $C_{d_i}$ et $\alpha$ est la racine courte (resp. longue).

Dans le premier cas, suivant **1.9**, on décompose

$$\begin{aligned}&\lambda(s_\alpha)J_{s_\alpha P|P}(\sigma\otimes\chi) \qquad\qquad (\#)\\ =&\mu_{r,s_\beta r^{-1}}(\chi)\mu_{s_\beta,r^{-1}}(\chi)\lambda(s_\alpha r s_\beta^{-1}r^{-1})\lambda(r)J_{r^{-1}P|P}(s_\beta r^{-1}\sigma\otimes s_\beta r^{-1}\chi)\times\\ &\times\lambda(s_\beta)J_{s_\beta P|P}(r^{-1}\sigma\otimes r^{-1}\chi)\lambda(r^{-1})J_{rP|P}(\sigma\otimes\chi).\end{aligned}$$

On a $\Sigma(P)\cap\Sigma(r\overline{P})\cap\Sigma(rs_\beta P)=\emptyset$, puisque $r\beta\in\Sigma(P)$ et $\Sigma(P)\cap\Sigma(s_\beta\overline{P})=\{\beta\}$. La fonction $\mu_{s_\beta,r}$ est donc constante. De même, l'ensemble $\Sigma(P)\cap\Sigma(rs_\beta\overline{P})\cap\Sigma(s_\alpha P)$ a une intersection vide avec $\Sigma_{\mathcal{O},\mu}$, ce qui montre que $\mu_{r,s_\beta r^{-1}}$ est constante. Comme $r\sigma\simeq\sigma$ et que $J_{rP|P}(\sigma\otimes\chi)$ est d'après **1.10** bien défini et bijectif pour tout $\chi$, on voit, en appliquant ce qui précède à $\rho_{r^{-1}\sigma,s_\beta}$ et en utilisant l'égalité ${}^rX_\beta=X_\alpha$ (cf. **1.5**), que l'opérateur $\rho_{s_\alpha}$ est bien défini et bijectif. Comme $\rho_{s_\beta}$ est induit par fonctorialité d'un isomorphisme $s_\beta E \to E$, on en déduit de même pour $\rho_{s_\alpha}$, après avoir appliqué les propriétés de compatibilité de $\rho_{s_\beta}$ à la formule de décomposition pour $\rho_{s_\alpha}^{-1}$ ci-dessus. Les propriétés de compatibilité pour $\rho_{s_\alpha}$ en résultent.

Si $M_\alpha$ n'est pas standard et $\alpha$ appartient à une composante irréductible de type $B_{d_i}$ ou $C_{d_i}$, alors notons $w$ l'élément de longueur minimale dans $W$ qui envoie un élément $x$ dans $A_M$ sur l'élément de $A_0$ déduit de $x$, en échangeant les composantes $(x_{i,1},\dots,x_{i,d_i})$ et $(x_{r,1},\dots,x_{r,d_r})$. Le groupe $wMw^{-1}$ est alors un sous-groupe de Levi standard de $G$ et la racine $\beta = w\alpha$ de $A_{wMw^{-1}}$ est dans $\Delta_{w\,\mathcal{O},\mu}$ avec $(wMw^{-1})_\beta$ standard. Notons $P_w$ le sous-groupe parabolique standard de $G$ de sous-groupe de Levi $wMw^{-1}$. Les racines dans $\Sigma_{red}(P)\cap\Sigma_{red}(w^{-1}\overline{P_w})$ sont de la forme $x\mapsto x_{i,j_1}x_{i_2,j_2}^{-1}$ avec $i<i_2\leq r$ et $x\mapsto x_{i_1,j_1}x_{r,j_2}^{-1}$ avec $i\leq i_1<r$. L'intersection avec $\Sigma_{\mathcal{O},\mu}$ est donc vide. Il en est de même pour celle de $w(\Sigma_{red}(P)\cap \Sigma_{red}(w^{-1}\overline{P_w}))$ avec $\Sigma_{s_\beta w\,\mathcal{O},\mu}$. En conséquence, les opérateurs $\lambda(w)J_{w^{-1}P_w|P}(\sigma)$ et $\lambda(w^{-1})J_{wP|P_w}(s_\beta w\sigma)$ sont d'après **1.10** bien définis et bijectifs. Les fonctions $\mu_{w^{-1},s_\beta w}$ et $\mu_{s_\beta,w}$ sont constantes. Par la formule de décomposition **1.9**, on trouve donc, en tant que fonction rationnelle en $\chi$,

$$\begin{aligned}&\lambda(s_\alpha)J_{s_\alpha P|P}(\sigma\otimes\chi)\\ =&\mu_{w^{-1},s_\beta w}\mu_{s_\beta,w}\lambda(s_\alpha w^{-1}s_\beta^{-1}w)\times\\ &\times\lambda(w^{-1})J_{wP|P_w}(s_\beta w(\sigma\otimes\chi))\lambda(s_\beta)J_{s_\beta P_w|P_w}(w(\sigma\otimes\chi))\lambda(w)J_{w^{-1}P_w|P}(\sigma\otimes\chi).\end{aligned}$$

Comme $Y_\beta({}^{w}\chi) = {}^{w^{-1}}Y_\beta(\chi) = Y_\alpha(\chi)$, il en résulte

$$\rho_{\sigma,s_\alpha} = \mu_{w^{-1},s_\beta w}\mu_{s_\beta,w}\lambda(w^{-1})J_{wP|P_w}(s_\beta w\sigma)\lambda(s_\beta)\rho_{w\sigma,s_\beta}\lambda(w)J_{w^{-1}P_w|P}(\sigma),$$

i.e., suite à ce qui a été préalablement établi pour $\rho_{w\sigma,s_\beta}$, $\rho_{\sigma,s_\alpha}^{-1}$ est bien défini et bijectif. En fait, on déduit des propriétés de compatibilités pour $\rho_{w\sigma,s_\beta}$ que $\rho_{\sigma,s_\alpha}$ est un produit de $\rho_{w\sigma,s_\beta}$ par un scalaire non nul. Comme $s_\alpha E = w^{-1}s_\beta(wE)$, les propriétés de compatibilité pour $\rho_{\sigma,s_\alpha}$ sont immédiates.

La dernière assertion de la proposition résulte directement de la décomposition (#) pour $r\in R(\mathcal{O})$, après l'avoir multipliée avec $X_\alpha - 1$, évaluée en 1 et après avoir utilisé l'égalité ${}^{r}X_\beta = X_\alpha$.

Supposons maintenant avoir montré l'égalité

$$(*)\qquad \rho_{w'}^{-1}\rho_w^{-1} = (\mu_{w,w'}^{-1}\prod_\alpha (Y_\alpha-1)(Y_\alpha^{-1}-1))_{|\chi=1}\lambda(ww'w''^{-1})\rho_{ww'}^{-1}.$$

avec $w''=ww'$, lorsque $w,w'$ sont des éléments dans $W_{\mathcal{O}}$ tels que $\rho_w$ et $\rho_{w'}$ soient bien définis, bijectifs et tels que $\rho_{w'}$ vérifie les deux propriétés de compatibilité, le produit portant sur l'ensemble $\Sigma_{\mathcal{O},\mu}(P)\cap w'^{-1}\Sigma_{\mathcal{O},\mu}(\overline{P})\cap w'^{-1}w^{-1}\,\Sigma_{\mathcal{O},\mu}(P)$.

Déduisons-en d'abord que $\rho_w$ est bien défini et bijectif, lorsque $w \in W_{\mathcal{O}}$, et que les propriétés de compatibilité sont bien vérifiées. Pour cela, effectuons une récurrence sur la longueur de $w$. Le cas $l_{\mathcal{O}}(w) = 1$ a déjà été établi. Soit $w'' = ws_\alpha \in W_{\mathcal{O}}$ avec $\alpha \in \Delta_{\mathcal{O},\mu}$, $l_{\mathcal{O}}(w'') = l_{\mathcal{O}}(w) - 1$. Par hypothèse de récurrence, les opérateurs $\rho_w$ et $\rho_{s_\alpha}$ sont bien définis, bijectifs, et ils vérifient les propriétés de compatibilité. On peut donc appliquer $(*)$, et on trouve

$$\rho_{w''}^{-1} = \lambda(w'' s_\alpha^{-1} w^{-1}) \rho_{s_\alpha}^{-1} \rho_w^{-1},$$

l'ensemble $\Sigma_{\mathcal{O},\mu}(P) \cap s_\alpha \Sigma_{\mathcal{O},\mu}(\overline{P}) \cap w^{-1} s_\alpha^{-1} \Sigma_{\mathcal{O},\mu}(P)$ étant vide grâce à l'hypothèse $l(ws_\alpha) = l(w) + l(s_\alpha)$. Comme $w'' s_\alpha^{-1} w^{-1} \in M$, il s'ensuit que $\rho_{w''}^{-1}$ est bijectif, que $\rho_{w''}$ est bien défini et qu'il vérifie les propriétés de compatibilité.

Il reste à montrer l'égalité $(*)$. Soient $w, w' \in W_{\mathcal{O}}$ vérifiant les hypothèses indiquées. On a

$$\begin{aligned}
&\rho_{w'}^{-1} \rho_w^{-1} \\
=&\rho_{w'}^{-1} ((\prod_\alpha (Y_\alpha - 1)) \lambda(w) J_{w^{-1}P|P}(\sigma \otimes \chi))_{|\chi=1} \\
=&((\prod_\alpha (Y_\alpha - 1)) \lambda(w) J_{w^{-1}P|P}(w'\sigma \otimes \chi))_{|\chi=1} \rho_{w'}^{-1} \\
=&((\prod_\alpha ({}^{w'^{-1}} Y_\alpha - 1))(\prod_\beta (Y_\beta - 1)) \lambda(w) J_{w^{-1}P|P}(w'\sigma \otimes w'\chi) \times \\
&\times \lambda(w') J_{w'^{-1}P|P}(\sigma \otimes \chi))_{|\chi=1} \\
=&((\prod_\alpha ({}^{w'^{-1}} Y_\alpha - 1))(\prod_\beta (Y_\beta - 1)) \mu_{w,w'}^{-1} \lambda(ww') J_{w'^{-1}w^{-1}P|P}(\sigma \otimes \chi))_{|\chi=1},
\end{aligned}$$

le produit sur $\alpha$ portant sur l'ensemble $\Sigma_{\mathcal{O},\mu}(P) \cap w^{-1} \Sigma_{\mathcal{O},\mu}(\overline{P})$ et celui sur $\beta$ sur l'ensemble $\Sigma_{\mathcal{O},\mu}(P) \cap w'^{-1} \Sigma_{\mathcal{O},\mu}(\overline{P})$. Remarquons l'égalité ensembliste

$$\begin{aligned}
&(w'^{-1} \Sigma_{\mathcal{O},\mu}(P) \cap w'^{-1} w^{-1} \Sigma_{\mathcal{O},\mu}(\overline{P})) \cup (\Sigma_{\mathcal{O},\mu}(P) \cap w'^{-1} \Sigma_{\mathcal{O},\mu}(\overline{P})) \\
=&(\Sigma_{\mathcal{O},\mu}(P) \cap w'^{-1} w^{-1} \Sigma_{\mathcal{O},\mu}(\overline{P})) \\
&\cup \pm(\Sigma_{\mathcal{O},\mu}(P) \cap w'^{-1} \Sigma_{\mathcal{O},\mu}(\overline{P}) \cap w'^{-1} w^{-1} \Sigma_{\mathcal{O},\mu}(P)),
\end{aligned}$$

les réunions étant disjointes. Comme ${}^{w'^{-1}} Y_\alpha = Y_{w'^{-1}\alpha}$ par le lemme **1.5**, on en déduit bien l'égalité $(*)$. □

**2.5** On va maintenant fixer un isomorphisme $\rho_w : wE \to E$ d'abord pour $w \in R(\mathcal{O})$ et ensuite pour $w \in W(M, \mathcal{O})$.

Rappelons que $R(\mathcal{O})$ est un produit de sous-groupes $R(\mathcal{O})_i$ d'ordre 2 (cf. **1.15** (ii)). Notons $s_{\alpha_i}$ l'élément non trivial de $R(\mathcal{O})_i$. Choisissons pour tout $i$ un isomorphisme $\rho_{s_{\alpha_i}} : s_{\alpha_i} E \to E$ tel que $\rho_{s_{\alpha_i}}^2 \sigma(s_{\alpha_i}^{-2}) = id_E$. (Ceci est toujours possible, quitte à multiplier $\rho_{s_{\alpha_i}}$ par un nombre complexe non nul. L'isomorphisme $\rho_{s_{\alpha_i}}$ est uniquement déterminé au produit par $\pm 1$ près.)

Changer le représentant de $s_{\alpha_i}$ à droite par un élément de $m \in M$ revient alors à composer $\rho_{s_{\alpha_i}}$ à droite avec $\sigma(m)$, comme on le vérifie facilement.

Il est clair que les $\rho_{s_{\alpha_i}}$ commutent entre eux, puisque les $s_{\alpha_i}$ agissent sur des composantes différentes de $\sigma$. Pour $r \in R(\mathcal{O})$, $r = s_{\alpha_{i_1}} \cdots s_{\alpha_{i_l}}$, on pose $\rho_r = \rho_{s_{\alpha_{i_1}}} \cdots \rho_{s_{\alpha_{i_l}}}$. Par ce qui précède, il est clair que ceci est indépendant de la décomposition choisie pour $r$. Multiplier le représentant de $r$ à droite par un élément $m \in M$ revient à composer $\rho_r$ à droite avec $\sigma(m)$. En effet, cela revient à multiplier chaque élément $s_{\alpha_i}$ dans la décomposition de $r$ à droite par un certain élément $m_i$ de $M_{s_{\alpha_i}}$ par notre convention dans **1.5** sur le choix des représentants des éléments dans $R(\mathcal{O})$.

Finalement, si $w \in W(M, \mathcal{O})$, on peut par **1.12** écrire $w = rw_{\mathcal{O}}$ avec $r \in R(\mathcal{O})$ et $w_{\mathcal{O}} \in W_{\mathcal{O}}$ déterminés de façon unique par $w$, et on pose alors $\rho_w = \rho_r \rho_{w_{\mathcal{O}}}$.

L'isomorphisme $i_P^G wE \to i_P^G E$ déduit de $\rho_w$ par fonctorialité sera (par abus de notation) noté encore $\rho_w$. Multiplier le représentant de $w$ à droite par un élément $m$ de $M$ revient à composer $\rho_w$ à droite par $i_P^G(\sigma(m))$.

**2.6** *Définition:* Soit $\chi \in \mathrm{Stab}(\mathcal{O})$ et rappelons l'isomorphisme $\phi_{\sigma,\chi} : E \to E_\chi$ défini dans **1.17**. On notera $\phi_\chi$ l'automorphisme de $E_B$ qui envoie $e \otimes b$ sur $(\phi_{\sigma,\chi} e) \otimes b_{\chi^{-1}}$. On écrira encore $\phi_\chi$ pour les automorphismes de $i_P^G E_B$ et de $i_P^G E_{K(B)}$ déduits de $\phi_\chi$ par fonctorialité.

**2.7** Le lemme suivant sera utile dans la suite:

**Lemme:** *Soient $\alpha$, $\alpha' \in \Delta_{\mathcal{O}}$ des racines simples, $s = s_\alpha$ et $s' = s_{\alpha'}$. Soient $m_s \in M_\alpha^1 \cap M$ et $m_{s'} \in M_{\alpha'}^1 \cap M$. Alors*

$$b_{m_s}\ {}^{s}b_{m_{s'}}\ {}^{ss'}b_{m_s}\ {}^{ss's}b_{m_{s'}} \cdots = b_{m_{s'}}\ {}^{s'}b_{m_s}\ {}^{s's}b_{m_{s'}}\ {}^{s'ss'}b_{m_s} \cdots,$$

*le nombre de facteurs de chaque côté étant égal à l'ordre de $ss'$.*

*Preuve:* C'est un calcul simple, en considérant les différents systèmes de racines de rang 2: en effet, par **1.5**, il existe des entiers $k$ et $k'$ tels que $b_{m_s} = Y_\alpha^k$ et $b_{m_{s'}} = Y_{\alpha'}^{k'}$. Par ailleurs, on a ${}^{s}b_{m_s} = Y_\alpha^{-k}$, ${}^{s'}b_{m_{s'}} = Y_{\alpha'}^{-k'}$, ${}^{s'}b_{m_s} = Y_{s'(\alpha)}^k$ et ${}^{s}b_{m_{s'}} = Y_{s(\alpha')}^{k'}$. Traitons à titre d'exemple le cas où $\alpha$ et $\alpha'$ engendrent un système de racines de type $A_2$. L'ordre de $ss'$ est alors 3, et on a $s(\alpha') = \alpha + \alpha' = s'(\alpha)$. On en déduit ${}^{s'}b_{m_s} = Y_\alpha^k Y_{\alpha'}^k$, ${}^{s}b_{m_{s'}} = Y_\alpha^{k'} Y_{\alpha'}^{k'}$, ${}^{ss'}b_{m_s} = Y_{\alpha'}^k$ et ${}^{s's}b_{m_{s'}} = Y_\alpha^{k'}$, d'où

$$b_{m_s}\ {}^{s}b_{m_{s'}}\ {}^{ss'}b_{m_s} = Y_\alpha^k Y_\alpha^{k'} Y_{\alpha'}^{k'} Y_{\alpha'}^k = Y_\alpha^{k'} Y_\alpha^k Y_{\alpha'}^k Y_\alpha^{k'} = b_{m_{s'}}\ {}^{s'}b_{m_s}\ {}^{s's}b_{m_{s'}}.$$

□

**3.** Considérons $\mathrm{End}_G(i_P^G E_B)$ comme $B$-module à gauche. On va construire pour tout $w \in W(M, \mathcal{O})$ un opérateur d'entrelacement $A_w \in \mathrm{Hom}_G(i_P^G E_B, i_P^G E_{K(B)})$ qui se prolonge de façon canonique en un élément de $\mathrm{End}_G(i_P^G E_{K(B)})$.

**3.1** Pour $\chi \in \mathfrak{X}^{\mathrm{nr}}(M)$, notons $B_\chi$ l'idéal des polynômes $b \in B$ avec $b(\chi) = 0$ et $sp_\chi$ les applications de spécialisation $E \otimes B \to E$ et $E \otimes K(B) \to E$, $e \otimes b \mapsto b(\chi)e$, ainsi que les applications induites $i_P^G E_B \to i_P^G E$ et $i_P^G E_{K(B)} \to i_P^G E$. (Dans le cas de $E_{K(B)}$, l'application $sp_\chi$ n'est définie que sur le sous-espace des éléments qui sont réguliers en $\chi$.) Soit $w \in W(M)$. D'après [W, IV.1.1], il existe $b_w \in B$ et un homomorphisme de $G - B$-modules

$$J_{B,w} : i_P^G E_B \to i_P^G (wE_B)$$

tel que, pour tout $\chi \in \mathfrak{X}^{\mathrm{nr}}(M)$, $v \in i_P^G E_B$,

$$b_w(\chi)(\lambda(w) J_{w^{-1}P|P}(\sigma \otimes \chi) sp_\chi v) = sp_\chi J_{B,w} v.$$

L'homomorphisme $b_w^{-1} J_{B,w}$ définit un élément de $\mathrm{Hom}_G(i_P^G E_B, i_P^G (wE_{K(B)}))$ que l'on notera $J_{K(B),w}$. Soit maintenant $w \in W(M, \mathcal{O})$. Alors on note $A_w$ l'élément de $\mathrm{Hom}_G(i_P^G E_B, i_P^G E_{K(B)})$ obtenu, en composant $J_{K(B),w}$ à gauche avec $\rho_w \circ \tau_w$. Il se prolonge de façon canonique en un élément de l'algèbre $\mathrm{End}_G(i_P^G E_{K(B)})$.

**Proposition:** *Soit $w \in W(M, \mathcal{O})$. Multiplier le représentant de $w$ à droite par un élément $m \in M$ revient à multiplier l'opérateur $A_w$ à gauche par ${}^w b_m^{-1}$. En particulier, $A_w$ ne dépend pas du choix d'un représentant dans $K$ de $w$.*

*Preuve:* Par la définition ci-dessus, celle dans **2.5** et la remarque dans **1.1**, on a $\rho_{wm} = \rho_w i_P^G(\sigma(m))$ et $J_{K(B),wm} = i_P^G(\sigma_B(m^{-1})) J_{K(B),w}$, d'où le résultat. □

**3.2 Lemme:** *Soient $w \in W(M, \mathcal{O})$, $v \in i_P^G E_B$ et $\chi \in \mathfrak{X}^{\mathrm{nr}}(M)$ tels que $\sigma \otimes \chi$ soit un point régulier pour $J_{w^{-1}P|P}$.*

*Alors*

$$sp_\chi A_w v = \rho_w \lambda(w) J_{w^{-1}P|P}(\sigma \otimes w^{-1}\chi) sp_{w^{-1}\chi} v.$$

*Preuve:* On a

$$sp_\chi A_w v = sp_\chi \rho_w \tau_w J_{K(B),w} v = \rho_w sp_\chi \tau_w J_{K(B),w} v.$$

Mais, l'application canonique $B \to B/B_\chi$ composée à droite avec $\tau_w$ a comme noyau $B_{w^{-1}\chi}$. Cette application composée est donc égale à $sp_{w^{-1}\chi}$. L'expression ci-dessus est donc égale à $\rho_w sp_{w^{-1}\chi} J_{K(B),w} v$. Ceci est bien égal à l'expression de l'énoncé. □

**3.3 Proposition:** *Soient $w, w' \in W_{\mathcal{O}}$. Alors,*

$$A_w A_{w'} = {}^{ww'}\mu_{w,w'}^{-1}(\mu_{w,w'} \prod_{\alpha} (Y_\alpha - 1)^{-1}(Y_\alpha^{-1} - 1)^{-1})_{|\chi=1} A_{ww'},$$

*le produit portant sur l'ensemble $\Sigma_{\mathcal{O},\mu}(P) \cap w'^{-1}\Sigma_{\mathcal{O},\mu}(\overline{P}) \cap w'^{-1}w^{-1}\Sigma_{\mathcal{O},\mu}(P)$.*

*Preuve:* Il suffit de montrer que, pour tout $\chi \in \mathfrak{X}^{\mathrm{nr}}(M)$, les deux opérateurs ont la même spécialisation en $\chi$. Posons $w'' = ww'$. D'après le lemme **3.2** et la propriété de commutation de $\rho_w$ remarquée en **2.4**, la spécialisation en $\chi$ de l'opérateur à gauche est égale à

$$\begin{aligned}
&\rho_w \lambda(w) J_{w^{-1}P|P}(\sigma \otimes w^{-1}\chi) \rho_{w'} \lambda(w') J_{w'^{-1}P|P}(\sigma \otimes w'^{-1} w^{-1}\chi)) sp_{w'^{-1}w^{-1}\chi} \\
&= \rho_w \rho_{w'} \lambda(w) J_{w^{-1}P|P}(w'\sigma \otimes w^{-1}\chi) \lambda(w') J_{w'^{-1}P|P}(\sigma \otimes w'^{-1} w^{-1}\chi) sp_{w'^{-1}w^{-1}\chi} \\
&= \rho_w \rho_{w'} \lambda(w) \lambda(w')\ {}^{ww'}\mu_{w,w'}^{-1}(\chi) J_{w'^{-1}w^{-1}P|P}(\sigma \otimes w'^{-1} w^{-1}\chi) sp_{w'^{-1}w^{-1}\chi} \\
&= \rho_w \rho_{w'} sp_{w^{-1}\chi} \tau_{w'} \lambda(ww'w''^{-1})\ {}^{ww'}\mu_{w,w'}^{-1}(\chi) J_{K(B),ww'} \\
&= \rho_w sp_{w^{-1}\chi} \rho_{w'} \tau_{w'} \lambda(ww'w''^{-1})\ {}^{ww'}\mu_{w,w'}^{-1}(\chi) J_{K(B),ww'} \\
&= sp_\chi \rho_w \rho_{w'} \tau_{ww'} \lambda(ww'w''^{-1})\ {}^{ww'}\mu_{w,w'}^{-1}(\chi) J_{K(B),ww'}.
\end{aligned}$$

Il reste alors à remplacer $\rho_w \rho_{w'}$ par l'expression donnée dans la proposition **2.4**, d'où la proposition. □

**3.4 Corollaire:** *Soient $w, w' \in W_{\mathcal{O}}$ tels que $l_{\mathcal{O}}(ww') = l_{\mathcal{O}}(w) + l_{\mathcal{O}}(w')$. Alors $A_w A_{w'} = A_{ww'}$. Par ailleurs, si $s = s_\alpha$ est une réflexion simple dans $W_{\mathcal{O}}$, alors il existe un nombre complexe $c_s'' \neq 0$ tel que ${A_s}^2 = c_s''(\mu^{M_\alpha})^{-1}$.*

*Preuve:* Si $l_{\mathcal{O}}(ww') = l_{\mathcal{O}}(w) + l_{\mathcal{O}}(w')$, alors $\Sigma_{\mathcal{O},\mu}(P) \cap w'^{-1}\Sigma_{\mathcal{O},\mu}(\overline{P}) \cap w'^{-1}w^{-1}$ $\Sigma_{\mathcal{O},\mu}(P)$ est l'ensemble vide, et la fonction rationnelle $\mu_{w,w'}$ est constante. La première assertion est donc une conséquence immédiate de **3.3**. La deuxième assertion résulte de **3.3** et du fait que $\mu_{s,s}$ est égal à $\mu^{M_\alpha}$ multiplié par une constante et que $A_1 = id$. □

**3.5 Proposition:** *(i) Pour $w \in W_{\mathcal{O}}$ et $r \in R(\mathcal{O})$, on a*

$$A_r A_w = A_{rw} = A_{rwr^{-1}} A_r.$$

*(ii) Pour tous $r, r' \in R(\mathcal{O})$, on a $A_r A_{r'} = A_{rr'}$.*

*(iii) Pour tout $b \in B$ et tout $w \in W(M, \mathcal{O})$, $A_w b = {}^w b A_w$.*

*Preuve:* (i) Comme $r\Delta_{\mathcal{O},\mu} = \Delta_{\mathcal{O},\mu}$ pour $r \in R(\mathcal{O})$, il suffit, suite au corollaire **3.4**, de considérer le cas où $w$ est un élément simple $s_\alpha$ de $W_{\mathcal{O}}$. En appliquant le même raisonnement comme dans la preuve de la proposition **3.3**, on trouve pour $\chi$ un élément générique de $\mathfrak{X}^{\mathrm{nr}}(M)$ que

$$sp_\chi(A_r A_{s_\alpha}) = \rho_r \rho_{s_\alpha} \lambda(r)\lambda(s_\alpha) J_{s_\alpha^{-1} r^{-1} P|P}(\sigma \otimes s_\alpha^{-1} r^{-1}\chi) sp_{s_\alpha^{-1} r^{-1}\chi} = sp_\chi(A_{rs_\alpha})$$

et

$$sp_\chi(A_{rs_\alpha r^{-1}} A_r) = \rho_{s_{r\alpha}} \rho_r \lambda(s_{r\alpha})\lambda(r) J_{r^{-1} s_{r\alpha}^{-1} P|P}(\sigma \otimes r^{-1} s_{r\alpha}^{-1}\chi) sp_{r^{-1} s_{r\alpha}^{-1}\chi}.$$

Ces deux expressions sont égales, si et seulement si

$$\rho_r \rho_{s_\alpha} \lambda(r)\lambda(s_\alpha) = \rho_{s_{r\alpha}} \rho_r \lambda(s_{r\alpha})\lambda(r),$$

ce qui équivaut à

$$\rho_r^{-1} \rho_{s_{r\alpha}}^{-1} \rho_r \rho_{s_\alpha} = \lambda(s_{r\alpha} r s_\alpha^{-1} r^{-1}).$$

Or, comme le côté gauche de cette dernière égalité est égal à $\rho_{r\sigma, s_{r\alpha}}^{-1} \rho_{\sigma, s_\alpha}$, ceci résulte de la dernière assertion de la proposition **2.4**

(ii) Il suffit de considérer le cas où $r = s_{\alpha_i}$ est l'élément non nul d'un $R(\mathcal{O})_i$. On va d'abord considére le cas où la projection de $r'$ sur $R(\mathcal{O})_i$ est nulle. Par le même raisonnement que dans la preuve de la proposition **3.3**, on trouve pour $\chi$ générique

$$sp_\chi(A_r A_{r'}) = \rho_r \rho_{r'} \lambda(r)\lambda(r') J_{r'^{-1} r^{-1} P|P}(\sigma \otimes {r'}^{-1} r^{-1}\chi) sp_{r'^{-1} r^{-1}\chi}.$$

Comme $\rho_r \rho_{r'} = \rho_{rr'}$ si la projection de $r'$ sur $R(\mathcal{O})_i$ est nulle, cette expression est égale à $sp_\chi(A_{rr'})$. On voit également que les opérateurs $A_{s_{\alpha_i}}$ et $A_{r'}$ commutent.

Supposons maintenant la projection de $A_{r'}$ sur $R(\mathcal{O})_i$ non nulle. Il suffit alors, d'après la commutativité, de considérer le cas $r = r'$. Pour $\chi$ générique dans $\mathfrak{X}^{\mathrm{nr}}(M)$, on trouve, par un raisonnement analogue à celui ci-dessus,

$$sp_\chi(A_r^2) = \rho_r \rho_r \lambda(r^2) J_{P|P}(\sigma \otimes r^2\chi) sp_{r^2\chi} = sp_\chi,$$

puisque $r^2 \in M \cap K$ et que $\rho_r$ a été choisi tel que $\rho_r^2 \lambda(r^2) = \mathrm{id}$.

(iii) C'est immédiat. □

**3.6 Proposition:** *L'espace* $\mathrm{Hom}_M(E_B, E_{K(B)})$ *est isomorphe à* $\oplus_{\chi \in \mathrm{Stab}(\mathcal{O})} K(B)\phi_\chi$ *en tant que* $K(B)$*-espace vectoriel.*

*Preuve:* On a $E_B = \mathrm{ind}_{M^1}^M E_{|M^1}$. Il faut donc déterminer $\mathrm{Hom}_M(\mathrm{ind}_{M^1}^M E, E_{K(B)})$. Observons que $(E_{|M^1})^\vee = (E^\vee)_{|M^1}$, puisque $M^1$ est un sous-groupe ouvert de $M$. Comme $M^1$ et $M$ sont unimodulaires, la réciprocité de Frobenius relative à l'induction compacte [BZ1, 2.29] donne donc $\mathrm{Hom}_M(\mathrm{ind}_{M^1}^M E_{|M^1}, E_{K(B)}) \simeq \mathrm{Hom}_{M^1}(E_{|M^1}, (E_{K(B)})_{|M^1})$. On définit un isomorphisme $\beta : \mathrm{Hom}_M(\mathrm{ind}_{M^1}^M E_{|M^1}, E_{K(B)}) \to \mathrm{Hom}_{M^1}(E_{|M^1}, (E_{K(B)})_{|M^1})$, $\varphi \mapsto \beta(\varphi)$, en posant $\beta(\varphi)(e) = \varphi(v_e)$, où $v_e$ désigne l'élément de $\mathrm{ind}_{M^1}^M E$ de valeur 0 pour $m \in M - M^1$ et de valeur $e$ en 1. Notons $E_{\overline{m}}$, $\overline{m} \in M/M^\sigma$, les différentes composantes irréductibles de $E_{|M^1}$. Elles sont deux à deux non isomorphes d'après la proposition **1.16**. Il en résulte que $\mathrm{Hom}_{M^1}(E_{|M^1}, (E_{K(B)})_{|M^1}) = \bigoplus_{\overline{m}\in M/M^\sigma} \mathrm{Hom}_{M^1}(E_{\overline{m}}, E_{\overline{m}} \otimes K(B))$. Choisissons une $\mathbb{C}$-base $f_i$, $i \in I$, de $K(B)$. Alors, en tant que $M^1$-module, $E_{\overline{m}} \otimes K(B) = \bigoplus_{i\in I} E_{\overline{m}} \otimes f_i$. On déduit alors du lemme de Schur que les éléments de $\mathrm{Hom}_{M^1}(E_{\overline{m}}, E_{\overline{m}} \otimes K(B))$ sont de la forme $e \mapsto e \otimes f$ avec $f \in K(B)$. Par suite, $\mathrm{Hom}_{M^1}(E_{|M^1}, (E_{K(B)})_{|M^1})$ est un $K(B)$-module isomorphe à $K(B)^{M/M^\sigma}$.

Pour prouver la proposition, il suffit donc grâce à la $K(B)$-linéarité de $\beta$ de montrer que les $\beta(\phi_\chi)$ forment une base du $K(B)$-module $\mathrm{Hom}_{M^1}(E_{|M^1}, (E_{K(B)})_{|M^1})$. Identifions $\mathrm{ind}_{M^1}^M E_{|M^1}$ et $E \otimes B$. On a $\beta(\phi_\chi)(e) = \phi_\chi(e \otimes 1) = \phi_{\sigma,\chi}(e) \otimes 1$. La projection de $\beta(\phi_\chi)$ sur $\mathrm{Hom}_{M^1}(E_{\overline{m}}, E_{\overline{m}} \otimes K(B))$ est donc la multiplication par $\chi(m)$. On déduit alors de la dualité des groupes $\mathrm{Stab}(\mathcal{O})$ et $M/M^\sigma$ et de l'indépendance linéaire des caractères que les $\beta(\phi_\chi)$ forment en effet une base du $K(B)$-module $\mathrm{Hom}_{M^1}(E_{|M^1}, (E_{K(B)})_{|M^1})$. □

**3.7 Proposition:** *Les opérateurs d'entrelacement $\phi_\chi A_w$, où $\chi \in \mathrm{Stab}(\mathcal{O})$ et $w \in W(M, \mathcal{O})$, sont $K(B)$-linéairement indépendants dans $\mathrm{Hom}_G(i_P^G E_B, i_P^G E_{K(B)})$ muni de la structure de $K(B)$-espace vectoriel induite par celle de $i_P^G E_{K(B)}$.*

*Preuve:* Supposons par absurde les opérateurs $\phi_\chi A_w$ $K(B)$-linéairement dépendants. Il existe alors des éléments $b_{w,\chi}$ de $B$ tels que $\sum_{w\in W(M,\mathcal{O}),\chi\in\mathrm{Stab}(\mathcal{O})} b_{w,\chi}\phi_\chi A_w = 0$, et on peut supposer le nombre $l$ d'éléments $b_{w,\chi} \neq 0$ minimal. Comme aucun des $\phi_\chi A_w$ n'est nul, il existe au moins deux indices $(w_1, \chi_1)$ et $(w_2, \chi_2)$ avec $b_{w_1,\chi_1} \neq 0 \neq b_{w_2,\chi_2}$. On peut trouver $b \in B$ avec $({}^{w_1}b)_{\chi_1^{-1}} = b$ et $({}^{w_2}b)_{\chi_2^{-1}} \neq b$. On a $\sum_{w\in W(M,\mathcal{O})} bb_{w,\chi}\phi_\chi A_w = 0$ et $0 = \sum_{w\in W(M,\mathcal{O})} b_{w,\chi}\ \phi_\chi A_w b = \sum_{w\in W(M,\mathcal{O})} b_{w,\chi}\ ({}^w b)_{\chi^{-1}}\phi_\chi A_w$, où on a appliqué la règle de commutation facile $A_w b = {}^w b A_w$. Soustrayant ces deux équations, on obtient une combinaison linéaire des $\phi_\chi A_w$ qui est nulle et dont moins de $l$ coefficients sont non nuls. Ceci donne une contradiction. □

**3.8 Théorème:** *En tant que $K(B)$-espaces vectoriels, on a*

$$\mathrm{Hom}_G(i_P^G E_B, i_P^G E_{K(B)}) = \bigoplus_{w\in W(M,\mathcal{O}),\chi\in\mathrm{Stab}(\mathcal{O})} K(B)\phi_\chi A_w.$$

*Preuve:* Par la réciprocité de Frobénius, on a

$$\mathrm{Hom}_G(i_P^G E_B, i_P^G E_{K(B)}) = \mathrm{Hom}_M(r_P^G i_P^G E_B, E_{K(B)}).$$

Par le lemme géométrique [BZ2, 2.12], $r_P^G i_P^G E_B$ admet une filtration par des sous-espaces $\mathcal{F}_w$, $w \in W^M \backslash W^G / W^M$, dont les sous-quotients sont isomorphes à $wE_B$. On en déduit une filtration du $K(B)$-module $\mathrm{Hom}_M(r_P^G i_P^G E_B, E_{K(B)})$ dont les sous-quotients sont, en tenant compte de **3.6**, isomorphes à $\mathrm{Hom}_M(wE_B, E_{K(B)})$.

Comme $\mathrm{Hom}_M(wE_B, E_{K(B)}) \neq 0$, si et seulement si $w \in W(M, \mathcal{O})$ et alors, par la proposition **3.6**, $\mathrm{Hom}_M\ (wE_B, E_{K(B)})$ est de dimension $|\,\mathrm{Stab}(\mathcal{O})|$, on en déduit que la dimension de $\mathrm{Hom}_G(i_P^G E_B, i_P^G E_{K(B)})$ est égale à $|W(M, \mathcal{O})|\,|\,\mathrm{Stab}(\mathcal{O})|$.

Comme les opérateurs d'entrelacement $\phi_\chi A_w$, $w \in W(M, \mathcal{O})$, $\chi \in \mathrm{Stab}(\mathcal{O})$ sont linéairement indépendants par la proposition **3.7**, le théorème en résulte. □

**4.** Maintenant, on fixe une composante irréductible $E_1$ de $E_{|M^1}$, et on note $\sigma_1$ la représentation de $M^1$ dans cet espace. On va s'intéresser à l'algèbre $\mathrm{End}_G(i_P^G(\mathrm{ind}_{M^1}^M E_1))$. Remarquons que ni cette algèbre ni la représentation $\mathrm{ind}_{M^1}^M E_1$ ne dépendent du choix de $(\sigma_1, E_1)$. Notons $B_\mathcal{O}$ la sous-algèbre de $B$ formée des polynômes qui sont invariants par translation par des éléments de $\mathrm{Stab}(\mathcal{O})$. D'après le corollaire **1.17**, c'est l'anneau des fonctions régulières de la variété affine quotient $X(M/M^1)/X(M/M^\sigma)$. Elle s'identifie donc à $\mathbb{C}[M^\sigma/M^1]$. En particulier, c'est un anneau factoriel, puisque $M^\sigma/M^1$ est un $\mathbb{Z}$-module libre de type fini et de même rang que $M/M^1$. Rappelons que l'on a noté $\mathcal{R}(M/M^1)$ un système de représentants de $M/M^1$.

**4.1 Lemme:** *L'isomorphisme canonique* $\mathrm{ind}_{M^1}^M E_{|M^1} \to E \otimes B$ *de* **2.1** *envoie* $\mathrm{ind}_{M^1}^M E_1$ *sur l'ensemble* $\sum_{m \in \mathcal{R}(M/M^1)} \sigma(m) E_1 \otimes b_m$. *En particulier, cet espace est un* $B_\mathcal{O}$*-sous-module de* $E \otimes B$.

*Preuve:* C'est immédiat. □

On identifiera dans la suite $\mathrm{ind}_{M^1}^M E_1$ et $\sum_{m \in \mathcal{R}(M/M^1)} \sigma(m) E_1 \otimes b_m$ au moyen de l'isomorphisme dans **2.1**. On écrira $E_{B_\mathcal{O}}$ pour $\mathrm{ind}_{M^1}^M E_1$, bien que cet espace ne soit pas isomorphe à $E \otimes B_\mathcal{O}$.

**4.2** Notons $K(B_\mathcal{O})$ le corps des fractions de $B_\mathcal{O}$ et posons $(E_{B_\mathcal{O}})_{K(B_\mathcal{O})} = E_{B_\mathcal{O}} \otimes_{B_\mathcal{O}} K(B_\mathcal{O})$.

**Lemme:** *On a*

$$\mathrm{Hom}_M(E_{B_\mathcal{O}}, E_{K(B_\mathcal{O})}) \simeq K(B_\mathcal{O}).$$

*Preuve:* Comme dans la preuve de la proposition **3.6**, on déduit de la réciprocité de Frobénius que

$$\begin{aligned}
&\mathrm{Hom}_M(ind_{M^1}^M E_1, (ind_{M^1}^M E_1)_{K(B_{\mathcal{O}})})\\
=&\mathrm{Hom}_{M^1}(E_1, (\mathrm{ind}_{M^1}^M E_1 \otimes_{B_{\mathcal{O}}} K(B_{\mathcal{O}}))_{|M^1})\\
=&\mathrm{Hom}_{M^1}(E_1, (\bigoplus_{m\in\mathcal{R}(M/M^\sigma)} \sigma(m)E_1 \otimes b_m) \otimes_{B_{\mathcal{O}}} K(B_{\mathcal{O}}))\\
=&\mathrm{Hom}_{M^1}(E_1, E_1 \otimes K(B_{\mathcal{O}}))\\
=&K(B_{\mathcal{O}}).
\end{aligned}$$

□

**4.3 Proposition:** *On a*

$$\bigcap_{\eta\in\mathrm{Stab}(\mathcal{O})} \ker(\phi_\eta - \mathrm{id}) = (E_{B_{\mathcal{O}}})_{K(B_{\mathcal{O}})}.$$

*En particulier, les automorphismes $\phi_\eta$ sont triviaux sur $E_{B_{\mathcal{O}}}$.*

*Preuve:* Fixons des systèmes de représentants $\mathcal{R} = \mathcal{R}(M/M^1)$ de $M/M^1$ et $\mathcal{R}' = \mathcal{R}(M/M^\sigma)$ de $M/M^\sigma$ respectivement. Tout élément de $(E_{B_{\mathcal{O}}})_{K(B_{\mathcal{O}})}$ s'écrit sous la forme $b^{-1}\sum_{m\in\mathcal{R}}\sigma(m)e_m \otimes b_m$ avec $e_m \in E_1$ et $b \in B_{\mathcal{O}}$. On vérifie tout de suite qu'un tel élément est dans le noyau de $\phi_\eta - \mathrm{id}$ pour tout $\eta$.

Inversement, soit $v \in E_{K(B)}$. Il peut s'écrire sous la forme $b^{-1}\sum_{(m,m')\in\mathcal{R}\times\mathcal{R}'} \sigma(m')e_{m,m'} \otimes b_m$ avec $e_{m,m'} \in E_1$ et $b \in B$. Soit $\eta \in \mathrm{Stab}(\mathcal{O})$. Alors,

$$\phi_\eta(v) = b_{\eta^{-1}}^{-1} \sum_{(m,m')\in\mathcal{R}\times\mathcal{R}'} \sigma(m')\eta(m')\eta(m)^{-1}e_{m,m'} \otimes b_m.$$

Donc, $\phi_\eta(v) = v$ implique que

$$b_{\eta^{-1}} \sum_{(m,m')\in\mathcal{R}\times\mathcal{R}'} \sigma(m')e_{m,m'} \otimes b_m = b \sum_{(m,m')\in\mathcal{R}\times\mathcal{R}'} \sigma(m')\eta(m'm^{-1})e_{m,m'} \otimes b_m.$$

Par suite, pour tout $m' \in \mathcal{R}'$,

$$b_{\eta^{-1}} \sum_{m\in\mathcal{R}} e_{m,m'} \otimes b_m = b \sum_{m\in\mathcal{R}} \eta(m'm^{-1})e_{m,m'} \otimes b_m. \qquad (*)$$

On va d'abord montrer que toute composante primaire de $b$ divise $b_{\eta^{-1}}$. Comme $b$ et $b_{\eta^{-1}}$ sont des monômes de Laurent de même degré, il s'ensuivra qu'ils ne diffèrent que par une constante, et il sera alors clair que $b^{-1}b_{\eta^{-1}} = \eta(m'm^{-1})$ si $e_{m,m'} \neq 0$.

Soit $b_0$ un facteur irréductible de $b$. On peut supposer que $b_0$ ne divise pas $\sum_{(m,m')\in\mathcal{R}\times\mathcal{R}'}\sigma(m')e_{m,m'}\otimes b_m$ dans $E_B$. Il existe donc $m'\in\mathcal{R}'$ tel que $b_0$ ne divise pas $\sum_{m\in\mathcal{R}}e_{m,m'}\otimes b_m$. Fixons une base $(e_i)_{i\in I}$ de $E_1$ et écrivons $e_{m,m'}=\sum_i c_i(m,m')e_i$ avec $c_i(m,m')\in\mathbb{C}$. Alors l'égalité (*) implique que

$$b_{\eta^{-1}}\sum_{m\in\mathcal{R}}\sum_{i\in I}c_i(m,m')e_i\otimes b_m=b\sum_{m\in\mathcal{R}}\sum_{i\in I}c_i(m,m')\eta(m'm^{-1})e_i\otimes b_m.$$

On en déduit $b_{\eta^{-1}}\sum_{m\in\mathcal{R}}c_i(m,m')b_m=b\sum_{m\in\mathcal{R}}c_i(m,m')\eta(m'm^{-1})b_m$ pour tout $i$. Par ailleurs, il existe au moins un $i$ tel que $b_0$ ne divise pas $\sum_{m\in\mathcal{R}}c_i(m,m')b_m$ dans $B$. Par suite, la composante $b_0$-primaire de $b$ divise $b_{\eta^{-1}}$.

On a donc bien $b^{-1}b_{\eta^{-1}}=\eta(m'm^{-1})$, si $e_{m,m'}\neq 0$. Supposons qu'il existe $(m_0,m_0')\in\mathcal{R}\times\mathcal{R}'$ tel que $e_{m_0,m_0'}\neq 0$. Posons $h=m_0'm_0^{-1}$. Alors, $\frac{b(\chi\eta^{-1})}{b(\chi)}$ est constant en $\chi$ de valeur $\eta(h)$.

Écrivons $b=\sum_{m\in M/M^1}c_mb_m$ avec $c_m\in\mathbb{C}$. Alors $b_{\eta^{-1}}=\sum_{m\in M/M^1}c_m\eta(m)^{-1}$ $b_m$. Comme $b_{\eta^{-1}}=\eta(h)b$, on trouve

$$\sum_{m\in M/M^1}c_m\eta(m)^{-1}b_m=\sum_{m\in M/M^1}\eta(h)c_mb_m.$$

Considérant $M/M^1$ comme groupe de caractères de $\mathfrak{X}^{\mathrm{nr}}(M)$, on déduit de l'indépendance linéaire des caractères que $\eta(h)c_m=\eta^{-1}(m)c_m$ pour tout $m$ et tout $\eta\in\mathrm{Stab}(\mathcal{O})$. Donc, $c_m\neq 0$ implique $m\in h^{-1}M^\sigma$, d'où $b=\sum_{m\in M^\sigma/M^1}c_mb_{h^{-1}m}=b_{h^{-1}}b'$ avec $b':=\sum_{m\in M^\sigma/M^1}c_mb_m\in B_{\mathcal{O}}$.

Comme $e_{m,m'}\neq 0$ implique par ce qui précédait $\eta(m'm^{-1})=\eta(h)$ pour tout $\eta\in\mathrm{Stab}(\mathcal{O})$, on a alors $m'm^{-1}\in hM^\sigma$. On en déduit $v=b'^{-1}b_h\sum_{m\in\mathcal{R}}\sigma(mh)e_m\otimes b_m$ avec $e_m\in E_1$. Or, ceci vaut $b'^{-1}\sum_{m\in\mathcal{R}}\sigma(mh)e_m\otimes b_{mh}$, et, par suite, $v$ a la forme indiquée. □

**4.4 Lemme:** *Supposons $M$ Levi maximal de $G$ et $W_{\mathcal{O}}\neq 1$. Soit $s$ l'unique élément $\neq 1$ de $W_{\mathcal{O}}$ et $(\sigma^1,E^1)$ une composante irréductible de $\sigma_{|M\cap G^1}$. Alors $s\sigma^1\simeq\sigma^1$.*

*Preuve:* Comme $K$ et $U$ sont contenus dans $G^1$, l'espace $i_{P\cap K}^KE^1$ est bien défini comme sous-espace $G^1$-invariant de $i_{P\cap K}^KE$. Posons $\Sigma_{\mathcal{O}}=\{\pm\alpha\}$. Comme $\overline{U}\subseteq G^1$, il résulte directement de la définition des opérateurs d'entrelacement que les opérateurs $J_{\overline{P}|P}(\sigma\otimes\chi_{\lambda\alpha})$, $\lambda\in\mathbb{C}$, envoient l'espace $i_{P\cap K}^KE^1$ dans l'espace $i_{\overline{P}\cap K}^KE^1$, si $\lambda$ est régulier. L'opérateur $J_{P|\overline{P}}$ possède la propriété analogue. Comme le composé $J_{P|\overline{P}}(\sigma\otimes\chi_{\lambda\alpha})J_{\overline{P}|P}(\sigma\otimes\chi_{\lambda\alpha})$ est scalaire, qu'il possède un pôle d'ordre 2 en $\lambda=0$, et que les opérateurs $J_{P|\overline{P}}(\sigma\otimes\chi_{\lambda\alpha})$ et $J_{\overline{P}|P}(\sigma\otimes\chi_{\lambda\alpha})$ possèdent au plus des pôles d'ordre 1, il en résulte que la restriction de $J_{\overline{P}|P}(\sigma\otimes\chi_{\lambda\alpha})$ à $i_{P\cap K}^KE^1$ admet un pôle d'ordre 1 en $\lambda = 0$. On va en déduire que $s\sigma^1 \simeq \sigma^1$. Plus précisément, on va montrer que $J_{\overline{P}|P}(\sigma \otimes \chi_{\lambda\alpha})$ est régulier en $\lambda = 0$, si $s\sigma^1 \not\simeq \sigma^1$. Si $G^1 = G$, ceci est prouvé dans [W, IV.2.2], puisque $\sigma$ est cuspidal unitaire et que $s$ normalise $M$. Comme en général $G^1 \neq G$ et que ce n'est pas nécessairement un groupe réductif, mais seulement un groupe localement profini, il faut généraliser cette preuve, ce qui demande quelques préparations préliminaires.

Pour simplifier, posons $P^1 = (G^1 \cap M)U$. Remarquons d'abord que dans la théorie des représentations lisses des groupes localement profinis [BZ1, 2.], la représentation de $G^1$ dans $i_{P\cap K}^K E^1$ est la représentation induite par la représentation $(\delta_P)_{|G^1\cap M}^{1/2}\sigma^1$ prolongée trivialement à $(G^1 \cap M)U$. On la notera $i_{P^1}^{G^1}\sigma^1$. Posons $V^1 = i_{P\cap K}^K E^1$. Comme $\overline{U} \subseteq G^1$, l'espace vectoriel $V^1(\overline{U})$ engendré par les éléments de la forme $v^1 - (i_{(G^1\cap M)U}^{G^1}\sigma^1)(\overline{u})v^1$, $v^1 \in V^1$, $\overline{u} \in \overline{U}$, est contenu dans $V^1$. Par ailleurs, il est $G^1 \cap M$-invariant. Il s'ensuit que le foncteur de Jacquet $r_{\overline{P}}^G$ envoie $V^1$ sur un sous-espace $G^1 \cap M$-invariant de $r_{\overline{P}}^G i_{P\cap K}^K E$ que l'on munira de l'action de $G^1 \cap M$ donnée par la représentation $r_{\overline{P}}^G i_P^G \sigma$. On écrira $r_{\overline{P}^1}^{G^1} V^1$ et $r_{\overline{P}^1}^{G^1} i_{P^1}^{G^1}\sigma^1$. On vérifie la réciprocité de Frobenius

$$\mathrm{Hom}_{G^1}(i_{P^1}^{G^1}\sigma^1, i_{\overline{P}^1}^{G^1}\sigma^1) = \mathrm{Hom}_{G^1\cap M}(r_{\overline{P}^1}^{G^1} i_{P^1}^{G^1}\sigma^1, \sigma^1), \qquad \Phi \mapsto r_{\overline{P}^1}^{G^1}\Phi,$$

où $(r_{\overline{P}^1}^{G^1}\Phi)(r_{\overline{P}^1}^{G^1}v^1) := (\Phi(v^1))(1)$. L'anneau $B^1 = \mathbb{C}[M \cap G^1/M^1]$ s'identifie à l'anneau des fonctions régulières sur $\mathfrak{X}^{\mathrm{nr}}(M \cap G^1) := \mathfrak{X}^{\mathrm{nr}}(M)_{|M\cap G^1}$. Posons $E_{B^1}^1 = E^1 \otimes B^1$ et notons $\sigma_{B^1}^1$ la représentation de $M \cap G^1$ dans $E_{B^1}^1$ donnée par $\sigma_{B^1}^1(m^1)(e^1 \otimes b^1) = \sigma^1(m^1)e^1 \otimes b^1 b_{m_1}$. Les résultats ci-dessus s'appliquent également à la représentation $V_{B^1}^1 := i_{P^1}^{G^1} E_{B^1}^1$, considérée comme sous-espace $G^1$-invariant de $i_P^G$ $E_B$. Notons $V^2$ et $V_{B^1}^2$ les sous-espaces vectoriels de $V^1$ (resp. $V_{B^1}^1$) formés des éléments à support dans $P^1\overline{P}^1$. L'application $p_1 : V^2 \to E^1$, $v^2 \mapsto \int_{\overline{U}} v^2(\overline{u})d\overline{u}$, induit un isomorphisme $r_{\overline{P}^1}^{G^1} V^2 \to E^1$ de représentations de $M \cap G^1$. L'application $p_s : V^1 \to sE^1$, $v^1 \mapsto v^1(s)$, définit un isomorphisme $r_{\overline{P}^1}^{G^1} V^1 / r_{\overline{P}^1}^{G^1} V^2 \to sE^1$. On a la situation analogue pour $V_{B^1}^1$ et $V_{B^1}^2$.

On va maintenant procéder à la preuve que $J_{\overline{P}|P}(\sigma \otimes \chi_{\lambda\alpha})$ est régulier en $\lambda = 0$, si $s\sigma^1 \not\simeq \sigma^1$: soit $n$ le plus petit entier tel que $\lambda^n J_{\overline{P}|P}(\sigma \otimes \chi_{\lambda\alpha})$ soit régulier en $\lambda = 0$. Notons cet opérateur $J(\lambda)$ et considérons l'opérateur $r_{\overline{P}^1}^{G^1} J(0)$ dans $\mathrm{Hom}_{M\cap G^1}(r_{\overline{P}^1}^{G^1} i_{P^1}^{G^1}\sigma^1, \sigma^1)$, obtenu par réciprocité de Frobenius. Comme le quotient $r_{\overline{P}^1}^{G^1} V^1 / r_{\overline{P}^1}^{G^1} V^2$ est isomorphe à $s\sigma^1$ et que cette représentation n'est par hypothèse pas isomorphe à $\sigma^1$, la restriction de $r_{\overline{P}^1}^{G^1} J(0)$ à $r_{\overline{P}^1}^{G^1} V^2$ est non triviale. Il existe donc $v^2 \in V_{B^1}^2$ tel que

$$0 \neq (r_{\overline{P}^1}^{G^1} J(0))(r_{\overline{P}^1}^{G^1} sp_1 v^2) = (J(0)sp_1 v^2)(1).$$

Mais, l'application $(J_{\overline{P}|P}(\sigma\otimes\chi_{\lambda\alpha})v^2)(1)$ est régulier en tout $\lambda\in\mathbb{C}$ [W, p. 283], puisque $v^2$ est à support dans $P\overline{P}$. Par conséquent, $n=0$. □

**4.5 Lemme:** *Pour tout $w\in W(M,\mathcal{O})$, il existe $m_w$ dans $M$ tel que $b_{m_w}A_w$ laisse invariant l'espace $(i_P^G E_{B_\mathcal{O}})_{K(B_\mathcal{O})}$. Si $w=s_\alpha$ avec $\alpha\in\Delta_{\mathcal{O},\mu}$ ou $w\in R(\mathcal{O})_i$, on peut choisir $m_w$ dans $M\cap M_\alpha^1$.*

*Preuve:* Par définition, $A_w=\rho_w\tau_w b_w^{-1}J_{B,w}$ avec $b_w\in B^\times$. En fait, on peut choisir $b_w\in B_\mathcal{O}^\times$, puisque la valeur de $J_{w^{-1}P|P}$ en une représentation $\sigma'$ ne dépend en un certain sens que de la classe d'isomorphie de $\sigma'$ [W, IV.1]. L'opérateur $\lambda(w)^{-1}$ $J_{B,w}$ ne dépendant que de la classe d'isomorphie de $\sigma$ [W, IV.1], il commute avec les opérateurs $\phi_\eta$, $\eta\in\mathrm{Stab}(\mathcal{O})$. On déduit donc de la proposition **4.3** que $J_{B,w}$ envoie l'espace $i_P^G E_{B_\mathcal{O}}$ dans l'espace $(i_P^G wE_{B_\mathcal{O}})_{K(B_\mathcal{O})}$. Considérons la restriction de l'isomorphisme $\rho_w\tau_w: i_P^G wE_B\to i_P^G E_B$ à $i_P^G\ wE_{B_\mathcal{O}}$.

Soit $v:=\sum_{m\in\mathcal{R}(M/M^1)}\sigma(m)e_m\otimes b_m$ dans $E_{B_\mathcal{O}}$ ($e_m\in E_1$). On a

$$\tau_w(v)=\sum_{m\in\mathcal{R}(M/M^1)}\sigma(m)e_m\otimes b_{wmw^{-1}}=\sum_{m\in\mathcal{R}(M/M^1)}(w\sigma)(m)e_{w^{-1}mw}\otimes b_m$$

. L'isomorphisme $\rho_w$ envoie $(w\sigma)(m)e_{w^{-1}mw}$ sur un élément dans $\sigma(m)\rho_{\sigma,w}E_1$, d'où

$$\rho_w\tau_w v\in\sum_{m\in\mathcal{R}(M/M^1)}\sigma(m)\rho_{\sigma,w}E_1\otimes b_m.$$

L'espace $\rho_{\sigma,w}E_1$ est le sous-espace de $E_{|M^1}$ muni de la représentation $w^{-1}\sigma_1$. Par **1.16**, il existe $m_w\in M$ tel que $\sigma(m_w)\rho_{\sigma,w}E_1$ soit égal à $E_1$. Il en résulte que l'opérateur $A_w$ laisse invariant $i_P^G E_{K(B_\mathcal{O})}$, si on multiplie le représentant de $w$ dans $G$ à droite par $w^{-1}m_w w$. Or, ceci revient par **3.1** à multiplier $A_w$ par $b_{m_w}^{-1}$.

Soit maintenant $\alpha\in\Delta_{\mathcal{O},\mu}$ et $s=s_\alpha\in W_\mathcal{O}$. Notons $\sigma^1$ la composante irréductible de $\sigma_{|M\cap M_\alpha^1}$ dont la restriction à $M^1$ contient $\sigma_1$. Comme $s\sigma^1\simeq\sigma^1$ par le lemme **4.4**, il existe $m_s\in M\cap M_\alpha^1$ tel que $s\sigma_1\simeq m_s\sigma_1$. L'opérateur $b_{m_s}A_s$ laisse l'espace $i_P^G E_{K(B_\mathcal{O})}$ invariant, par ce qui précède.

Si $w\in R(\mathcal{O})_i$, alors l'effet de $w$ sur $\sigma$ porte sur la seule composante $\sigma_i$ de $\sigma$ avec $w\sigma_i\simeq\sigma_i$. On peut donc choisir $m_w$ dans le groupe $\mathrm{GL}_{k_i}$ sur lequel $\sigma_i$ est définie. Il appartient donc à une composante de $M_\alpha$ qui est dans les notations de la proposition **1.13** égale à $\underline{H}_{k_i+k}$. Comme cette composante est semi-simple, $m_w$ appartient bien à $M_\alpha^1$. □

**4.6** *Définition:* Si $s$ est une symétrie simple dans $W_\mathcal{O}$, fixons un élément $m_s$ de $M\cap M_\alpha^1$ tel que les conclusions du lemme **4.5** soient vérifiées et posons $J_s=b_{m_s}A_s$. C'est un élément de $\mathrm{Hom}_G(i_P^G E_{B_\mathcal{O}}, i_P^G E_{K(B_\mathcal{O})})$. Pour $r\in R(\mathcal{O})_i$, fixons $m_r\in R(\mathcal{O})$ tel que les conclusions du lemme **4.5** relatives à $A_r$ soient vérifiées et posons $J_r=$ $b_{m_r}A_r$. Pour $w \in W_{\mathcal{O}}$, $w = s_1 \cdots s_l$, définissons $J_w = J_{s_1} \cdots J_{s_l}$. (Il résultera du lemme ci-après que cette définition est bonne.) Pour $r \in R(\mathcal{O})$, on fait de même à partir des $J_{s_{\alpha_i}}$, $\alpha_i \in R(\mathcal{O})_i$. Si finalement $w \in W(M, \mathcal{O})$, $w = rw_{\mathcal{O}}$, avec $r \in R(\mathcal{O})$ et $w_{\mathcal{O}} \in W_{\mathcal{O}}$, on écrira $J_w = J_r J_{w_{\mathcal{O}}}$.

**4.7 Lemme:** *(i) Soit $w \in W_{\mathcal{O}}$ et soit $w = s_1 \cdots s_l$ une décomposition de $w$ en symétries simples. L'opérateur $J_{s_1} \cdots J_{s_l}$ ne dépend que de $w$ et non pas de la décomposition de $w$.*

*(ii) Si $s = s_\alpha$ est une réflexion simple dans $W_{\mathcal{O}}$, alors $J_s^2 = c_s''(\mu^{M_\alpha})^{-1}$, où $c_s''$ est le nombre complexe défini dans* **3.4**.

*(iii) Pour tous $r, r' \in R(\mathcal{O})$, on a $J_r J_{r'} = J_{rr'}$.*

*Preuve:* Pour (i), d'après [Sp, 8.3.3], il suffit de montrer l'assertion suivante: si $s$ et $s'$ sont des symétries simples dans $W_{\mathcal{O}}$ et si $m(s, s')$ désigne l'ordre de $ss'$, alors

$$J_s J_{s'} J_s \cdots = J_{s'} J_s J_{s'} \cdots ,$$

le nombre de facteurs de chaque côté étant égal à $m(s, s')$. Vu le résultat de composition établi pour les opérateurs $A_s$ dans **3.5** et la définition de $J_s$ et de $J_{s'}$, ceci revient à prouver que

$$b_{m_s} {}^s b_{m_{s'}} {}^{ss'} b_{m_s} {}^{ss's} b_{m_{s'}} \cdots = b_{m_{s'}} {}^{s'} b_{m_s} {}^{s's} b_{m_{s'}} {}^{s'ss'} b_{m_s} \cdots .$$

Or, ceci résulte du lemme **2.7**.

Pour (ii), la preuve de **3.4** reste valable, puisque $J_s^2$ est la restriction de $A_s^2$.

Par itération, il suffit de vérifier l'égalité (iii) lorsque $r = s_{\alpha_i}$ est l'élément non trivial d'un $R(\mathcal{O})_i$. Si la projection de $r'$ sur $R(\mathcal{O})_i$ est triviale, alors l'assertion résulte directement de la définition **4.6**. Sinon, il suffit de montrer que $J_r^2 = 1$. Par **3.5**, on a

$$J_r^2 = b_{m_r} A_r b_{m_r} A_r = b_{m_r} {}^r b_{m_r}.$$

Comme ${}^r b_{m_r} = b_{m_r}^{-1}$, on a bien $b_{m_r} {}^r b_{m_r} = 1$. □

**4.8 Lemme:** *Les opérateurs $J_w$, $w \in W(M, \mathcal{O})$, sont $K(B_{\mathcal{O}})$-linéairement indépendants dans* $\mathrm{Hom}_G(i_P^G E_{B_{\mathcal{O}}}, i_P^G E_{K(B_{\mathcal{O}})})$, *et, on a, pour $b$ dans $K(B_{\mathcal{O}})$, la relation de commutation $J_w b = {}^w b J_w$.*

*Preuve:* La preuve de **3.7** se généralise, après avoir remarqué que les opérateurs $J_w$ sont non nuls. Or, ceci résulte du fait que $sp_\chi$ envoie l'espace $E_{B_{\mathcal{O}}}$, considéré comme sous-espace de $E_B$, sur l'espace $E$ muni de la représentation $\sigma \otimes \chi$. Comme $sp_\chi \circ J_{w,K(B)} = \lambda(w) J_{w^{-1}P|P}(\sigma \otimes \chi) sp_\chi$ et que $J_{w^{-1}P|P}(\sigma \otimes \chi) \neq 0$ pour $\chi$ régulier, l'opérateur $J_w$ est bien non nul.

La relation de commutation est une conséquence immédiate de la définition de $J_w$ et de la relation de commutation pour les opérateurs $A_w$ (cf. **3.5**) □

**4.9 Théorème:** *En tant que $K(B_{\mathcal{O}})$-modules,*

$$\mathrm{Hom}_G(i_P^G E_{B_{\mathcal{O}}}, i_P^G E_{K(B_{\mathcal{O}})}) = \bigoplus_{w\in W(M,\mathcal{O})} K(B_{\mathcal{O}}) J_w.$$

*Preuve:* Par la réciprocité de Frobénius, on a

$$\mathrm{Hom}_G(i_P^G E_{B_{\mathcal{O}}}, i_P^G E_{K(B_{\mathcal{O}})}) = \mathrm{Hom}_M(r_P^G i_P^G E_{B_{\mathcal{O}}}, E_{K(B_{\mathcal{O}})})$$

Par le lemme géométrique [BZ2, 2.12], $r_P^G i_P^G E_{B_{\mathcal{O}}}$ admet une filtration par des sous-espaces $\mathcal{F}_w$, $w \in W^M\backslash W^G/W^M$, dont les sous-quotients sont isomorphes à $wE_{B_{\mathcal{O}}}$. On en déduit une filtration du $K(B_{\mathcal{O}})$-module $\mathrm{Hom}_M(r_P^G i_P^G E_{B_{\mathcal{O}}}, E_{K(B_{\mathcal{O}})})$ dont les sous-quotients sont isomorphes à $\mathrm{Hom}_M(wE_{B_{\mathcal{O}}}, E_{K(B_{\mathcal{O}})})$.

Compte tenu de l'indépendance $K(B_{\mathcal{O}})$-linéaire des opérateurs $J_w$, $w \in W(M, \mathcal{O})$, il reste à prouver que

$$\mathrm{Hom}_G(wE_{B_{\mathcal{O}}}, (E_{B_{\mathcal{O}}})_{K(B_{\mathcal{O}})}) \simeq \begin{cases} 0, & \text{si } w \notin W(M,\mathcal{O}); \\ K(B_{\mathcal{O}}), & \text{sinon.} \end{cases}$$

Comme $wE_{B_{\mathcal{O}}}$ est un quotient de $wE_B$ et que $(E_{B_{\mathcal{O}}})_{K(B_{\mathcal{O}})}$ est un sous-module de $E_{K(B)}$, la nullité de $\mathrm{Hom}_M(wE_B, E_{K(B)})$ pour $w \notin W(M,\mathcal{O})$ implique celle de $\mathrm{Hom}_M(wE_{B_{\mathcal{O}}}, (E_{B_{\mathcal{O}}})_{K(B_{\mathcal{O}})})$.

Supposons maintenant $w \in W(M,\mathcal{O})$. On a $wE_{B_{\mathcal{O}}} = \mathrm{ind}_{M^1}^M wE_1$. Comme $W(M,\mathcal{O})$ permute les composantes irréductibles de $E_{|M^1}$, on trouve que $\mathrm{ind}_{M^1}^M wE_1 = E_{B_{\mathcal{O}}}$, d'où, d'après **4.2**, $\mathrm{Hom}_M(wE_{B_{\mathcal{O}}}, E_{K(B_{\mathcal{O}})}) = K(B_{\mathcal{O}})$. □

**5.** Dans cette section, on va définir pour tout $w \in W_{\mathcal{O}}$ un opérateur $T_w \in \mathrm{End}_G(i_P^G E_{B_{\mathcal{O}}})$ et montrer que les opérateurs $J_r T_w$, $r \in R(\mathcal{O})$ et $w \in W_{\mathcal{O}}$, forment une base du $B_{\mathcal{O}}$-module $\mathrm{End}_G(i_P^G E_{B_{\mathcal{O}}})$.

**5.1 Lemme:** *Soit $\alpha \in \Delta_{\mathcal{O}}$, $s = s_\alpha$. Posons $c_s = {c'_s}^{-1} c''_s$ (où $c'_s$ et $c''_s$ sont les nombres complexes définis respectivement en **1.5** et en **3.4**).*

*(i) L'opérateur $sp_1(X_\alpha - 1)J_s$ est scalaire sur $i_{P\cap K}^K E$ et*

$$(sp_1(X_\alpha - 1)J_s)^2 = c_s \frac{(1-q^{-a_s})^2(1+q^{-b_s})^2}{4}.$$

*Par ailleurs, si $\chi \in \mathfrak{X}^{\mathrm{nr}}(M)$ vérifie $X_\alpha(\chi) = 1$, alors, pour $v \in i_{P\cap K}^K E_{B_{\mathcal{O}}}$,*

$$sp_\chi(X_\alpha - 1)J_s v = (sp_1(X_\alpha - 1)J_s) sp_\chi v.$$

*(ii) Supposons que $\mu^{M_\alpha}$ s'annule en $X_\alpha=-1$. Fixons $\chi_{-1}\in\mathfrak{X}^{\mathrm{nr}}(M)$ tel que $X_\alpha(\chi_{-1})=-1$. L'opérateur $sp_{\chi_{-1}}(X_\alpha+1)J_s$ est scalaire sur $i_{P\cap K}^K E$ et*

$$(sp_{\chi_{-1}}(X_\alpha+1)J_s)^2=c_s\frac{(1+q^{-a_s})^2(1-q^{-b_s})^2}{4}.$$

*Par ailleurs, si $\chi\in\mathfrak{X}^{\mathrm{nr}}(M)$ vérifie $X_\alpha(\chi)=-1$, alors, pour $v\in i_{P\cap K}^K E_{B_\mathcal{O}}$,*

$$sp_\chi(X_\alpha+1)J_sv=(sp_{\chi_{-1}}(X_\alpha-1)J_s)sp_\chi v.$$

*Preuve:* (i) D'après le lemme **1.8** et nos choix, la représentation $i_{P\cap M_\alpha}^{M_\alpha}\sigma$ est irréductible, et la restriction de $J_s$ à $i_{P\cap K}^K E$ a un pôle simple en $Y_\alpha=1$. L'opérateur $sp_1(X_\alpha-1)J_s$ est défini à partir d'un entrelacement de cette représentation avec elle-même. Cela doit donc être un scalaire. On a

$$(sp_1(X_\alpha-1)J_s)^2=sp_1((X_\alpha-1)J_s(X_\alpha-1)J_s)=sp_1((X_\alpha-1)(X_\alpha^{-1}-1)J_s^2).$$

Comme $(J_s)^2=c_s''(\mu^{M_\alpha})^{-1}$, on trouve bien l'expression dans l'énoncé, après avoir remplacé $\mu^{M_\alpha}$ par l'expression donnée dans la proposition **1.6**.

Soit maintenant $\chi\in\mathfrak{X}^{\mathrm{nr}}(M)$ tel que $X_\alpha(\chi)=1$. Ils existent alors $\eta\in\mathrm{Stab}(\mathcal{O})$ et $\chi_1\in\mathfrak{X}^{\mathrm{nr}}(M_\alpha)$ tels que $\chi=\eta\chi_1$. On a ${}^s\chi_1=\chi_1$ et $\chi_1(m_s)=1$. Remarquons par ailleurs qu'en tant qu'espace vectoriel, $i_{P\cap K}^K E_{B_\mathcal{O}}=\bigoplus_{m\in M/M^\sigma} i_{P\cap K}^K\sigma(m)E_1\otimes b_m$. Soit $m\in M$ et $v_m\in i_{P\cap K}^K\sigma(m)E_1$. On a, en tant que fonction rationnelle en $\chi'\in\mathfrak{X}^{\mathrm{nr}}(M)$,

$$\begin{aligned}
&sp_{\chi'\chi}(X_\alpha-1)J_sv_m\otimes b_m\\
=&(X_\alpha(\chi\chi')-1)(\chi\chi')(m_s)\rho_{\sigma,s}\lambda(s)J_{sP|P}(\sigma\otimes{}^s(\chi'\chi))\,{}^s(\chi'\chi)(m)v_m\\
=&\chi_1(m)(X_\alpha(\chi')-1)(\eta\chi')(m_s)\rho_{\sigma,s}\lambda(s)J_{sP|P}(\sigma\otimes{}^s(\chi'\eta))\,{}^s\chi'(m)\,{}^s\eta(m)v_m\\
=&\chi_1(m)(X_\alpha(\chi')-1)(\eta\chi')(m_s)\rho_{\sigma,s}\lambda(s)J_{sP|P}(\sigma\otimes{}^s(\chi'\eta))\,{}^s\chi'(m)\phi_{\sigma,s\eta}v_m\\
=&\chi_1(m)(X_\alpha(\chi')-1)(\eta\chi')(m_s)\rho_{\sigma,s}\phi_{\sigma,s\eta}\lambda(s)J_{sP|P}(\sigma\otimes{}^s\chi')\,{}^s\chi'(m)v_m.
\end{aligned}$$

Spécialisant en $\chi'=1$, on trouve

$$sp_\chi(X_\alpha-1)J_sv_m\otimes b_m=\chi_1(m)\eta(m_s)\rho_{\sigma,s}\phi_{\sigma,s\eta}\rho_{\sigma,s}^{-1}sp_1(X_\alpha-1)J_sv_m\otimes b_m.$$

On vient de voir que $sp_1(X_\alpha-1)J_sv_m\otimes b_m$ est un élément de $i_{P\cap K}^K\sigma(m)E_1$. Soit $e_1\in E_1$. Remarquons que $\rho_{\sigma,s}^{-1}e_1$ est un élément de $(s\sigma)(m_s^{-1})E_1$. Alors

$$\begin{aligned}
\rho_{\sigma,s}\phi_{\sigma,s\eta}\rho_{\sigma,s}^{-1}\sigma(m)e_1&=\rho_{\sigma,s}\phi_{\sigma,s\eta}(s\sigma)(m)\rho_{\sigma,s}^{-1}e_1\\
&=\rho_{\sigma,s}\eta(m_s^{-1})\eta(m)(s\sigma)(m)\rho_{\sigma,s}^{-1}e_1\\
&=\eta(m_s^{-1})\eta(m)\sigma(m)e_1,
\end{aligned}$$

i.e. $(\rho_{\sigma,s}\phi_{\sigma,s\eta}\rho_{\sigma,s}^{-1})_{|i_{P\cap K}^K\sigma(m)E_1} = \eta(m_s^{-1})\eta(m)$. La dernière assertion du (i) du lemme en résulte.

(ii) Observons d'abord que ${}^s\chi_{-1}\chi_{-1}^{-1}\in \mathrm{Stab}(\mathcal{O})$: comme $\mu^{M_\alpha}$ s'annule en $X_\alpha = -1$, on trouve par les résultats de Harish-Chandra que $\sigma\otimes\chi_{-1}\simeq s(\sigma\otimes\chi_{-1})\simeq s\sigma\otimes\ {}^s\chi_{-1}\simeq\sigma\otimes\ {}^s\chi_{-1}$, d'où $\sigma\simeq\sigma\otimes\ {}^s\chi_{-1}\chi_{-1}^{-1}$. Posons $\eta_{-1}=\ {}^s\chi_{-1}\chi_{-1}^{-1}$. Soit $v_m\in i_{P\cap K}^K\sigma(m)E_1$. On a alors, en tant que fonction rationnelle en $\chi'\in\mathfrak{X}^{\mathrm{nr}}(M)$,

$$\begin{aligned}
&sp_{\chi'\chi_{-1}}(X_\alpha+1)J_sv_m\otimes b_m\\
=&(X_\alpha(\chi'\chi_{-1})+1)(\chi'\chi_{-1})(m_s)\rho_{\sigma,s}\lambda(s)J_{sP|P}(\sigma\otimes\ {}^s(\chi_{-1}\chi'))^s(\chi_{-1}\chi')(m)v_m\\
=&(X_\alpha(\chi'\chi_{-1})+1)(\chi'\chi_{-1})(m_s)\rho_{\sigma,s}\lambda(s)J_{sP|P}(\sigma\otimes\eta_{-1}\chi_{-1}\ {}^s\chi')(\chi_{-1}\ {}^s\chi')(m)\times\\
&\times\phi_{\sigma,\eta_{-1}}v_m\\
=&(X_\alpha(\chi'\chi_{-1})+1)(\chi'\chi_{-1})(m_s)\rho_{\sigma,s}\phi_{\sigma,\eta_{-1}}\lambda(s)J_{sP|P}(\sigma\otimes\chi_{-1}\ {}^s\chi')(\chi_{-1}\ {}^s\chi')(m)v_m\\
=&\chi_{-1}(m)\chi_{-1}(m_s)\rho_{\sigma,s}\phi_{\sigma,\eta_{-1}}\rho_{\sigma\otimes\chi_{-1},s}^{-1}(-X_\alpha(\chi')+1)\chi'(m_s)\rho_{\sigma\otimes\chi_{-1},s}\lambda(s)\times\\
&\times J_{sP|P}(\sigma\otimes\chi_{-1}\ {}^s\chi')\chi'(m)v_m
\end{aligned}$$

Notons $J_{s,\chi_{-1}}$ l'opérateur dans $\mathrm{End}_G((E_{\chi_{-1}})_{B_\mathcal{O}})_{K(B_\mathcal{O})}$ défini comme $J_s$, en remplaçant $\sigma$ par $\sigma\otimes\chi_{-1}$. Alors, d'après (i), on a pour $\chi\in\mathfrak{X}^{\mathrm{nr}}(M)$ avec $X_\alpha(\chi)=1$,

$$\begin{aligned}
&sp_{\chi_{-1}\chi}((X_\alpha+1)J_s)v_m\otimes b_m\\
=&\chi_{-1}(m)\chi_{-1}(m_s)\rho_{\sigma,s}\phi_{\sigma,\eta_{-1}}\rho_{\sigma\otimes\chi_{-1},s}^{-1}sp_\chi(-X_\alpha+1)J_{s,\chi_{-1}}v_m\otimes b_m\\
=&\chi_{-1}(m_s)\rho_{\sigma,s}\phi_{\sigma,\eta_{-1}}\rho_{\sigma\otimes\chi_{-1},s}^{-1}(sp_1(-X_\alpha+1)J_{s,\chi_{-1}})sp_{\chi_{-1}\chi}v_m\otimes b_m.
\end{aligned}$$

L'opérateur $\rho_{\sigma,s}\phi_{\sigma,\eta_{-1}}\rho_{\sigma\otimes\chi_{-1},s}^{-1}$ vient par fonctorialité d'un automorphisme d'une représentation irréductible dans $E$, puisque $s(\sigma\otimes\chi_{-1})=s\sigma\otimes\chi_{-1}\eta_{-1}$. C'est donc un scalaire. En posant ci-dessus $\chi=1$ et en utilisant la partie (i) du lemme relative à l'opérateur $J_{s,\chi_{-1}}$, on trouve que l'opérateur $sp_{\chi_{-1}}(X_\alpha+1)J_s$ est scalaire sur $i_{P\cap K}^KE$. L'égalité ci-dessus implique alors également la dernière assertion de la partie (ii) du lemme.

Concernant la deuxième assertion, il suffit de remarquer que

$$(sp_{\chi_{-1}}(X_\alpha+1)J_s)^2=sp_{\chi_{-1}}(X_\alpha+1)(X_\alpha^{-1}+1)J_s^2,$$

et on conclut comme pour (i). □

**5.2** *Définition:* Soit $\alpha\in\Delta_\mathcal{O}$, $s=s_\alpha$. Fixons une racine carrée $c_s^{1/2}$ de $c_s$. D'après le lemme **5.1**, il existe $\epsilon_1,\epsilon_{-1}\in\{\pm1\}$ tels que

$$sp_1(X_\alpha-1)J_s=\epsilon_1c_s^{1/2}\frac{(1-q^{-a_s})(1+q^{-b_s})}{2}sp_1$$

et $$sp_{\chi_{-1}}(X_\alpha+1)J_s=\epsilon_{-1}c_s^{1/2}\frac{(1+q^{-a_s})(1-q^{-b_s})}{2}sp_{\chi_{-1}}.$$

Posons
$$R_s=\begin{cases}-\epsilon_1 q^{a_s+b_s}c_s^{-1/2}J_s, & \text{si } \epsilon_1 b_s\neq\epsilon_{-1}b_s;\\ -\epsilon_1 q^{a_s+b_s}c_s^{-1/2}X_\alpha J_s, & \text{si } \epsilon_1 b_s=\epsilon_{-1}b_s.\end{cases}$$

et $$T_s=R_s+(q^{a_s+b_s}-1)\frac{X_\alpha(X_\alpha-\frac{q^{b_s}-q^{a_s}}{q^{a_s+b_s}-1})}{X_\alpha^2-1}.$$

**5.3 Lemme:** *Un élément* $A\in \mathrm{Hom}_G(i_P^G E_{B_{\mathcal{O}}}, i_P^G E_{K(B_{\mathcal{O}})})$ *est dans* $\mathrm{End}_G(i_P^G E_{B_{\mathcal{O}}})$, *si et seulement si, pour tout* $v\in i_P^G E_{B_{\mathcal{O}}}$ *et pour tout* $v^\vee\in i_{P\cap K}^K E^\vee$, *l'application* $\mathfrak{X}^{\mathrm{nr}}(M)\to\mathbb{C}$, $\chi\mapsto\langle sp_\chi Av, v^\vee\rangle$, *est régulière.*

*Preuve:* La condition est nécessaire, puisque, si $A(v)\in i_P^G E_{B_{\mathcal{O}}}$, alors $sp_\chi(A(v))\in i_P^G E_\chi$ pour tout $\chi\in\mathfrak{X}^{\mathrm{nr}}(M)$.

La condition est suffisante: Soit $v\in i_P^G E_{B_{\mathcal{O}}}$. Par l'inclusion $i_P^G E_{B_{\mathcal{O}}}\subseteq i_{P\cap K}^K E\otimes K(B)$, on peut écrire $Av=\sum_{i\in I} v_i\otimes b_i$ avec $b_i\in K(B)$ et les $v_i$ $\mathbb{C}$-linéairement indépendants dans $i_{P\cap K}^K E$.

Soient $v_i^\vee\in i_{P\cap K}^K E^\vee$, $i\in I$, des éléments duaux aux $v_i$. Alors $\langle sp_\chi(A(v)), v_i^\vee\rangle = b_i$. Il en résulte donc que $b_i\in B$ pour tout $i\in I$, i.e. $Av\in i_P^G E_B\cap(i_P^G E_{B_{\mathcal{O}}})_{K(B_{\mathcal{O}})} = i_P^G E_{B_{\mathcal{O}}}$. □

**5.4 Proposition:** *Les opérateurs* $T_{s_\alpha}$, $\alpha\in\Delta_{\mathcal{O}}$, *définis dans* **5.2**, *appartiennent à* $\mathrm{End}_G(i_P^G E_{B_{\mathcal{O}}})$.

*Preuve:* Posons $s=s_\alpha$. Il faut montrer que $T_s$ est régulier en tout $\chi\in\mathfrak{X}^{\mathrm{nr}}(M)$. Écrivons $T_s=p_sJ_s+r_s$ avec $p_s\in B^\times$ et $r_s\in K(B_{\mathcal{O}})$. Soient $v\in i_{P\cap K}^K E_{\mathcal{O}}$ et $v^\vee\in E^\vee$. L'application $\chi\mapsto\langle sp_\chi T_s v, v^\vee\rangle$ est régulière en $\chi$, sauf peut-être si $X_\alpha(\chi)=1$ ou si $X_\alpha(\chi)=-1$ et $b_s\neq 0$.

Remarquons d'abord que
$$Res_{X_\alpha=\pm 1} r_s=\frac{\pm 1}{2}(q^{a_s+b_s}-1-(\pm 1)(q^{b_s}-q^{a_s})).$$

La régularité en $\chi\in\mathfrak{X}^{\mathrm{nr}}(M)$ tel que $X_\alpha(\chi)=1$ résulte alors de l'égalité $sp_\chi(X_\alpha-1)p_sJ_s=-q^{a_s+b_s}\frac{(1-q^{-a_s})(1+q^{-b_s})}{2}sp_\chi$, et celle en $\chi\in\mathfrak{X}^{\mathrm{nr}}(M)$ tel que $X_\alpha(\chi)=-1$ de l'égalité $sp_\chi(X_\alpha+1)p_sJ_s=q^{a_s+b_s}\frac{(1+q^{-a_s})(1-q^{-b_s})}{2}sp_\chi$, si $b_s\neq 0$. □

**5.5 Proposition:** *Soit* $\alpha\in\Delta_{\mathcal{O}}$, $s=s_\alpha$. *Alors* $(T_s+1)(T_s-q^{a_s+b_s})=0$.

*Preuve:* Soit $T_s$ un opérateur de la forme $p_sJ_s + r_s$ avec $p_s \in \mathbb{C}^\times$ et $r_s \in K(B)$. Alors

$$\begin{aligned}(T_s+1)(T_s-q^{a_s+b_s}) =& T_s^2 + (1-q^{a_s+b_s})T_s - q^{a_s+b_s}\\ =& (r_s + {}^sr_s + 1 - q^{a_s+b_s})p_sJ_s\\ &+ p_s^2J_s^2 + r_s^2 + (1-q^{a_s+b_s})r_s - q^{a_s+b_s}.\end{aligned}$$

Comme 1 et $J_s$ sont $K(B)$-linéairement indépendants, l'équation $(T_s + 1)(T_s - q^{a_s+b_s}) = 0$ équivaut à

$$r_s + {}^sr_s + 1 - q^{a_s+b_s} = 0 \qquad \text{et} \qquad p_s^2J_s^2 + r_s^2 + (1-q^{a_s+b_s})r_s - q^{a_s+b_s} = 0.$$

Pour $T_s$ défini comme ci-dessus dans le cas $\epsilon_1 b_s \neq \epsilon_{-1} b_s$, la vérification de ces identités se fait par un calcul élémentaire. Le cas $\epsilon_1 b_s = \epsilon_{-1} b_s$ en résulte, puisque $(X_\alpha J_s)\ (X_\alpha J_s) = J_s^2$. □

**5.6 Définition:** Fixons pour tout $w \in W_\mathcal{O}$ une décomposition réduite en symétries simples $w = s_1 \cdot \ldots \cdot s_l$ et définissons $T_w = T_{s_1} \cdot \ldots \cdot T_{s_l}$. (En particulier, $T_1 = id$.)

**5.7 Proposition:** *Soient $w, w' \in W_\mathcal{O}$.*

*(i) Si $l_\mathcal{O}(ww') = l_\mathcal{O}(w) + l_\mathcal{O}(w')$, alors $T_wT_{w'}$ est de la forme $b_{w,w'}J_{ww'} + \sum_{w''} f_{w,w',w''}\ J_{w''}$, avec $b_{w,w'} \in B_\mathcal{O}^\times$ et $f_{w,w',w''} \in K(B_\mathcal{O})$, $f_{w,w',w''} = 0$ lorsque $l_\mathcal{O}(w'') \geq l_\mathcal{O}(ww')$.*

*(ii) Si $l_\mathcal{O}(ww') < l_\mathcal{O}(w) + l_\mathcal{O}(w')$, alors $T_wT_{w'}$ est de la forme $\sum_{w''} f_{w,w',w''}J_{w''}$ avec $f_{w,w',w''} \in K(B_\mathcal{O})$, $f_{w,w',w''} = 0$ lorsque $l_\mathcal{O}(w'') \geq l_\mathcal{O}(w) + l_\mathcal{O}(w')$.*

*Preuve:* Ceci résulte directement de la définition de $T_w$ et de **4.7** (i), (ii). □

**5.8 Lemme:** *Soit $w' \in W_\mathcal{O}$ et notons $w_0$ l'élément le plus long dans $W_\mathcal{O}$. Supposons pour tout $w \in W_\mathcal{O}$ donné un opérateur $T_{w,w'}$ dans* $\mathrm{End}_G(i_P^G E_{B_\mathcal{O}})$ *de la forme $\sum_{w''} f_{w,w''}J_{w''}$, avec $f_{w,w''} \in K(B_\mathcal{O})$, $f_{w,w_0} = 0$ si $w \neq w'$ et $f_{w',w_0} \in B_\mathcal{O}^\times$. Soit $\chi \in \mathfrak{X}^{\mathrm{nr}}(M)$ et soient $c_{r,w}$, $r \in R(\mathcal{O})$, $w \in W_\mathcal{O}$, des nombres complexes tels que $\sum_{r,w} c_{r,w} sp_\chi J_r T_{w,w'} = 0$. Alors $c_{1,w'} = 0$.*

*Preuve:* Fixons un sous-groupe $H$ de $G$ vérifiant les propriétés indiquées dans [H1, 5.] relatives à $P$ et $\mathcal{O}$. Fixons un élément $v \neq 0$ de $(i_{P\cap K}^K E)^H$ de support contenu dans $(P \cap K)H$ et à valeurs dans $E_1$. Il s'identifie donc à l'élément de $i_P^G E_{B_\mathcal{O}}$ défini pour $m \in M$, $u \in U$ et $k \in K$ par $muk \mapsto \delta_P^{1/2}\sigma_B(m)v(k)$. Soit $v^\vee$ un élément de $(i_{\overline{P}\cap K}^K E^\vee)^H$ de support contenu dans $(\overline{P} \cap K)H$. On suppose que $v^\vee$ a été choisi tel que $\langle (v(1), v^\vee(1)\rangle \neq 0$. L'opérateur rationnel $J_{P|\overline{w_0P}}\lambda(w_0)$ est régulier sur $\mathcal{O}^\vee$ par **1.10** (où $w_0$ désigne l'élément de longueur maximale de $W_\mathcal{O}$). Posons $A = J_{P|\overline{w_0 P}}(\sigma \otimes \cdot)\lambda(w_0)$.

On déduit de [H1, 5.3] (où il faut inverser les rôles de $v$ et $v^\vee$) que, pour $w \in W(M, \mathcal{O})$, $\langle J_{P|wP}(\sigma \otimes \chi')\lambda(w)\rho_w v, A(\chi')\rho_{w_0}^\vee v^\vee\rangle$ est une fonction constante en $\chi'$ qui ne peut être $\neq 0$ que si $w = w_0$. Si $w = w_0$, alors sa valeur est de la forme $c_w\langle v(1), v^\vee(1)\rangle$, où $c_w$ est un nombre complexe $\neq 0$ qui ne dépend pas de $v$ et $v^\vee$.

Supposons $\sum_{r,w} c_{r,w} sp_\chi J_r T_{w,w'} = 0$. Alors,

$$\begin{aligned}
&\langle \sum_{r,w} c_{r,w} J_r T_{w,w'} v, A\rho_{w_0}^\vee v^\vee\rangle \\
=&\sum_{r,w}\sum_{w''}\langle c_{r,w} J_r f_{w,w''} J_{w''} v, A\rho_{w_0}^\vee v^\vee\rangle \\
=&\sum_{r,w}\sum_{w''} c_{r,w} f_{w,w''}\langle J_r J_{w''} v, A\rho_{w_0}^\vee v^\vee\rangle
\end{aligned}$$

Observons maintenant que, d'après le lemme **3.2** et la preuve de **3.3**, pour $w \in W_\mathcal{O}$, $r \in R(\mathcal{O})$ et $\chi' \in \mathfrak{X}^{\mathrm{nr}}(M)$ un point régulier, et compte tenu de la définition de $J_r J_w$, $sp_{\chi'} J_r J_w v$ est le produit d'un nombre complexe non nul avec

$$\begin{aligned}
&\rho_{\sigma,rw}\lambda(rw) J_{w^{-1}r^{-1}P|P}(\sigma \otimes (w^{-1}r^{-1}\chi')) sp_{w^{-1}r^{-1}\chi'} v \\
=&\rho_{\sigma,rw} J_{P|rwP}(rw\sigma \otimes \chi')\lambda(rw) v \\
=&J_{P|rwP}(\sigma \otimes \chi')\lambda(rw)\rho_{\sigma,rw} v,
\end{aligned}$$

les propriétés de commutations pour $\rho_{\sigma,rw}$ par rapport à $J_{P|rwP}(\cdot)\lambda(rw)$ résultant du fait que celui-ci est défini par fonctorialité à partir d'un isomorphisme entre des représentations dans $E$ (cf. **2.5**).

D'après ce qui précède, on voit que le produit $\langle J_{P|rwP}(\sigma\otimes\chi')\lambda(rw)\rho_{\sigma,rw} v, A(\chi')$ $\rho_{w_0}^\vee v^\vee\rangle$ est nul, sauf si $r = 1$ et $w = w_0$, et alors sa valeur est par choix de $v^\vee$ non nulle égale à $c_{w_0}\langle (v(1), v^\vee(1)\rangle =: c$. Il s'ensuit que

$$\begin{aligned}
0 =& \langle \sum_{r,w} sp_\chi c_{r,w} J_r T_{w,w'} v, A(\chi)\rho_{w_0}^\vee v^\vee\rangle \\
=& c_{1,w'} f_{w',w_0}(\chi) c.
\end{aligned}$$

Comme $f_{w',w_0} \in B^\times$, on a $f_{w',w_0}(\chi) \neq 0$, et il en résulte que $c_{1,w'} c = 0$, d'où $c_{1,w'} = 0$. □

**5.9 Proposition:** *Pour tout $\chi \in \mathfrak{X}^{\mathrm{nr}}(M)$, les opérateurs $sp_\chi J_r T_w$, $r \in R(\mathcal{O})$, $w \in W_\mathcal{O}$, sont linéairement indépendants.*

*Preuve:* Soit $\chi \in \mathfrak{X}^{\mathrm{nr}}(M)$. On va d'abord prouver par récurrence décroissante sur $l$, $0 \leq l \leq l_\mathcal{O}(w_0)$, que, lorsque $c_{r,w}$, $r \in R(\mathcal{O})$, $w \in W_\mathcal{O}$, sont des nombres complexes tels que $\sum_{r,w} c_{r,w} sp_\chi J_r T_w = 0$, alors $c_{1,w} = 0$ pour tout $w$ de longueur supérieure ou égale à $l$.

Il résulte du lemme précédent et de **5.7** (en posant $w' = 1$) que cette affirmation est vraie pour $l = l(w_0)$. Soit donc $0 \leq l < l(w_0)$ et supposons-la vraie pour $l+1$. Soient $c_{r,w}$ des nombres complexes tels que $\sum_{r,w} c_{r,w} sp_\chi J_r T_w = 0$. On a donc $c_{1,w} = 0$ lorsque $l(w) > l$. Soit $w_1 \in W_\mathcal{O}$ tel que $l(w_1) = l$. Posons $w' = w_1^{-1} w_0$ et $A := \sum_{r,w} c_{r,w} J_r T_w$. Par hypothèse, $sp_\chi A = 0$ et donc $sp_\chi A T_{w'} = 0$.

D'autre part,

$$AT_{w'} = \sum_{r,w} c_{r,w} J_r (T_w T_{w'}).$$

On va déduire du lemme **5.8** que $c_{1,w_1} = 0$. Si $l(w) > l$, alors $c_{1,w} = 0$ par hypothèse de récurrence. Si $l(w) = l$, alors $l(ww') \leq l(w) + l(w') = l(w_0)$, puisque par [B, VI, paragraphe 1, Cor. 3], $l(w') = l(w_0) - l$. Par la proposition **5.7**, $T_w T_{w'}$ est, pour $l(w) \leq l$ et $w \neq w_1$, combinaison $K(B_\mathcal{O})$-linéaire d'opérateurs $J_{w''}$, $l(w'') < l(w_0)$, alors que $T_{w_1} T_{w'}$ est de la forme $f_{w_1,w_0} J_{w_0} + \sum_{w'' \neq w_0} f_{w_1,w''} J_{w''}$ avec $f_{w_1,w_0} \in B_\mathcal{O}^\times$ et $f_{w_1,w''} \in K(B_\mathcal{O})$. Les hypothèses du lemme **5.8** sont donc bien vérifiées, et, par suite, $c_{1,w_1} = 0$.

Déduisons-en maintenant que $c_{r',w} = 0$ pour tout $r' \in R(\mathcal{O})$ et $w \in W_\mathcal{O}$. Soit $r' \in R(\mathcal{O})$. Remarquons que $T_w J_{r'^{-1}} = J_{r'^{-1}} T_{r' w r'^{-1}}$ par **3.5**.

Comme $0 = sp_\chi A J_{r'^{-1}}$ et, d'autre part,

$$AJ_{r'^{-1}} = \sum_{r,w} c_{r,w} J_r J_{r'^{-1}}\ T_{r' w r'^{-1}},$$

on trouve, en utilisant **4.7** (iii),

$$0 = sp_\chi A J_{r'^{-1}} = \sum_{r,w} c_{r, r'^{-1} w r'} sp_\chi J_{r r'^{-1}} T_w.$$

Par ce qui a été montré dans la première partie, il en résulte bien que $c_{r',w} = 0$ pour tout $w \in W_\mathcal{O}$. □

**5.10 Théorème:**

$$\mathrm{End}_G(i_P^G E_{B_\mathcal{O}}) = \bigoplus_{r \in R(\mathcal{O}), w \in W_\mathcal{O}} B_\mathcal{O} J_r T_w.$$

*Preuve:* Il est clair par définition et **5.4** que l'espace de droite est inclus dans l'espace de gauche.

D'autre part, par le théorème **4.9** et **5.7**, l'espace de gauche est inclus dans $\bigoplus_{r,w} K(B_\mathcal{O}) J_r T_w$. Il suffit donc de prouver qu'une combinaison linéaire de la forme $A := \sum_{r,w} b^{-1} b_{r,w} J_r T_w$ avec $b_{r,w}$ et $b$ dans $B_{\mathcal{O}}$ est dans $\mathrm{End}_G(i_P^G E_{B_{\mathcal{O}}})$, si et seulement si $b$ divise tous les $b_{r,w}$ dans $B_{\mathcal{O}}$.

On se ramène au cas où $b$ et $b_{r,w}$ sont premiers entre eux dans $B_{\mathcal{O}}$. Il faut alors prouver que $b$ est une unité. Supposons par l'absurde que $b$ n'est pas une unité. L'anneau $B_{\mathcal{O}}$ étant factoriel, il suffit de considérer le cas où $b$ est irréductible. Alors il existe au moins un élément $\chi \in \mathfrak{X}^{\mathrm{nr}}(M)$ tel que $b(\chi) = 0$, alors que $b_{r,w}(\chi) \neq 0$ pour au moins un couple $(r, w)$. Comme $0 = sp_\chi bA = \sum_{r,w} b_{r,w}(\chi) sp_\chi J_r\, T_w$, ceci contredit la proposition **5.9** qui dit que les opérateurs $sp_\chi J_r T_w$ sont $\mathbb{C}$-linéairement indépendants pour tout $\chi \in \mathfrak{X}^{\mathrm{nr}}(M)$. □

**5.11 Corollaire:** *Le centre $\mathcal{Z}_{\mathcal{O}}$ de $\mathrm{End}_G(i_P^G E_{B_{\mathcal{O}}})$ est formé des éléments $W(M, \mathcal{O})$-invariants dans $B_{\mathcal{O}}$.*

*Preuve:* En effet, supposons $\sum_{r,w} b_{r,w} J_r J_w$ dans le centre. Alors, par les relations de commutation (cf. **4.8**), on a pour tout $r \in R(\mathcal{O})$, $w \in W_{\mathcal{O}}$, $b \in B_{\mathcal{O}}$, $bb_{r,w} = {}^{rw}b b_{r,w}$. Par suite, par l'indépendance linéaire des opérateurs $J_r J_w$ (cf. **4.8**), $b_{r,w} = 0$ si $rw \neq 1$. Le centre de $\mathrm{End}_G(i_P^G E_{B_{\mathcal{O}}})$ est donc contenu dans $B_{\mathcal{O}}$, et il est alors clair que c'est $B_{\mathcal{O}}^{W(M,\mathcal{O})}$. □

**5.12** Notons $K(\mathcal{Z}_{\mathcal{O}})$ le corps de fractions de $\mathcal{Z}_{\mathcal{O}}$ et posons $\mathrm{End}_G(i_P^G E_{B_{\mathcal{O}}})_{K(\mathcal{Z}_{\mathcal{O}})} = \mathrm{End}_G(i_P^G E_{B_{\mathcal{O}}}) \otimes_{\mathcal{Z}_{\mathcal{O}}} K(\mathcal{Z}_{\mathcal{O}})$.

**Corollaire:** *L'algèbre $\mathrm{End}_G(i_P^G E_{B_{\mathcal{O}}})_{K(\mathcal{Z}_{\mathcal{O}})}$ est un $K(B_{\mathcal{O}})$-module qui est canoniquement isomorphe à $\mathrm{Hom}_G(i_P^G E_{B_{\mathcal{O}}}, i_P^G E_{K(B_{\mathcal{O}})})$.*

*Preuve:* Remarquons tout d'abord que $B_{\mathcal{O}} \otimes_{\mathcal{Z}_{\mathcal{O}}} K(\mathcal{Z}_{\mathcal{O}})$ est canoniquement isomorphe à $K(B_{\mathcal{O}})$ (cf. [L1]). Il en résulte par **5.10** que $\mathrm{End}_G(i_P^G E_{B_{\mathcal{O}}})_{K(\mathcal{Z}_{\mathcal{O}})}$ est isomorphe à $\bigoplus_{r,w} K(B_{\mathcal{O}}) J_r T_w$, d'où le corollaire par **4.9** et **5.7**. □

**6.** Rappelons que l'on a défini en **1.5** pour tout $\alpha \in \Delta_{\sigma,\mu}$ un élément $h_\alpha$ dans $M \cap M_\alpha^1$ et un nombre réel $t_\alpha$. Posons $\widetilde{\alpha} := H_M({h_\alpha}^{t_\alpha})$. Désignons par $\Lambda_{\mathcal{O}}$ le $\mathbb{Z}$-module libre dans $a_M$ égal à l'image de $M^\sigma/M^1$ par $H_M$, par $\Lambda_{\mathcal{O}}^\vee$ le $\mathbb{Z}$-module libre inclus dans $a_M^*$ qui est par $\langle \cdot, \cdot \rangle$ en dualité parfaite avec $\Lambda_{\mathcal{O}}$, par $\Sigma_{\mathcal{O}}$ l'ensemble des $\widetilde{\alpha}$, $\alpha \in \Sigma_{\mathcal{O},\mu}$, et, finalement, par $\Sigma_{\mathcal{O}}^\vee$ l'ensemble des multiples $\alpha^*$ de $\alpha$, $\alpha \in \Sigma_{\mathcal{O},\mu}$, tels que $\langle \alpha^*, \widetilde{\alpha} \rangle = 2$.

**6.1 Proposition:** *Le quadruplet $(\Lambda_{\mathcal{O}}, \Sigma_{\mathcal{O}}, \Lambda_{\mathcal{O}}^\vee, \Sigma_{\mathcal{O}}^\vee)$ est une donnée radicielle. Le système de racines sous-jacent est réduit. Le groupe de Weyl de $\Sigma_{\mathcal{O}}$ est canoniquement isomorphe à $W_{\mathcal{O}}$, et l'ensemble $\Delta_{\mathcal{O}} = \{H_M(h_\alpha^{t_\alpha}), \alpha \in \Delta_{\mathcal{O},\mu}\}$ forme une base de $\Sigma_{\mathcal{O}}$.*

*Par ailleurs, si $\Sigma_\mu'$ est une composante irréductible de $\Sigma_{\mathcal{O},\mu}$, et si $\Sigma_\mu'$ n'est pas de type $C_n$, alors l'ensemble des $\widetilde{\alpha}$, $\alpha \in \Sigma_\mu'$, est une composante irréductible de*

*$\Sigma_{\mathcal{O}}$ qui est du même type que $\Sigma'$. Dans le cas contraire, la composante irréductible correspondante de $\Sigma_{\mathcal{O}}$ est de type $B_n$.*

D'abord un lemme:

**6.2 Lemme:** *Soit $\alpha \in \Sigma_{\mathcal{O},\mu}$. Notons $m_\alpha$ le plus grand nombre $> 0$, tel que $\chi_{\frac{2\pi i}{m_\alpha \log q}\alpha} \in \mathfrak{X}^{\mathrm{nr}}(M_\alpha)$. Alors $H_M(h_\alpha) = \frac{m_\alpha}{2}\alpha^\vee$.*

*Preuve:* Écrivons $H_M(h_\alpha) = \frac{m}{2}\alpha^\vee$ avec $m > 0$. Alors

$$1 = \chi_{\frac{2\pi i}{m_\alpha \log q}\alpha}(h_\alpha) = q^{-\langle \frac{2\pi i}{m_\alpha \log q}\alpha, H_M(h_\alpha)\rangle} = q^{-\frac{2\pi i}{\log q}\frac{m}{m_\alpha}}.$$

Il s'ensuit que $m \in m_\alpha\mathbb{Z}$. Mais, alors, grâce à la maximalité de $m_\alpha$, il faut que $m_\alpha = m$. □

*Preuve:* (de la proposition) Fixons une composante $\Sigma_{\mathcal{O},\mu,i}$ de $\Sigma_{\mathcal{O},\mu}$, et notons $\{\alpha_1, \ldots, \alpha_d\}$ sa base. Rappelons que $d = d_i$, sauf si $\Sigma_{\mathcal{O},\mu,i}$ est de type $A_d$. Alors $d = d_i - 1$. Un élément de $\mathfrak{X}^{\mathrm{nr}}(M)$ est de la forme $|\det_{m_1}|^{s_{1,1}} \otimes \cdots \otimes |\det_{m_1}|^{s_{1,d_1}} \otimes |\det_{m_2}|^{s_{2,1}} \otimes \cdots \otimes |\det_{m_2}|^{s_{2,d_2}} \otimes \cdots \otimes |\det_{m_r}|^{s_{r,1}} \otimes \cdots \otimes |\det_{m_r}|^{s_{r,d_r}} \otimes 1$, où les $s_{i,j}$ sont des nombres complexes. Pour $j = 1, \ldots, d_i - 1$, $\chi_{\alpha_j}$ est donnée par $s_{i,j} = 1/m_i$, $s_{i,j+1} = -1/m_i$, les autres exposants $s_{i',j'}$ étant nuls. Si $d = d_i$, alors $\chi_{\alpha_d}$ est donné par $s_{i,d} = 1/m_i$ et $s_{i,d-1} = 0$ (resp. $s_{i,d-1} = 1/m_i$), si $\Sigma_{\mathcal{O},\mu,i}$ est de type $B_d$ (resp. $D_d$), et par $s_{i,d} = 2/m_i$ si $\Sigma_{\mathcal{O},\mu,i}$ est de type $C_d$, les autres exposants étant nuls.

Fixons une uniformisante $\widetilde{\omega}$ de $F$. On peut choisir $h_{\alpha_j} = \mathrm{diag}(1, 1, \ldots, 1, \widetilde{\omega}, \widetilde{\omega}^{-1}, 1, \ldots, 1)$ pour $j = 1, \ldots, d_i - 1$, et, si $d = d_i$, $h_{\alpha_d} = \mathrm{diag}(1, \ldots, 1, \widetilde{\omega}, 1, \ldots, 1)$ si $\Sigma_{\mathcal{O},\mu,i}$ est de type $B_d$ ou $C_d$ , et $h_{\alpha_d} = \mathrm{diag}(1, 1, \ldots, 1, \widetilde{\omega}, \widetilde{\omega}, 1, \ldots, 1)$ si $\Sigma_{\mathcal{O},\mu,i}$ est de type $D_d$, $\widetilde{\omega}$ se trouvant chaque fois à la dernière place sur la diagonale de la $j^{\text{ème}}$ (ou $d^{\text{ème}}$) copie de $\mathrm{GL}_{m_i}$.

Posons $\epsilon = 2$, si $\Sigma_{\mathcal{O},\mu,i}$ est de type $B_d$, $\epsilon = 1$ sinon. Alors, on en déduit $m_{\alpha_j} = 2m_i^{-1}$ pour $j = 1, \ldots, d_i - 1$, et, si $d = d_i$, $m_{\alpha_d} = \frac{2}{\epsilon}m_i^{-1}$. Si $t_i$ désigne l'ordre du stabilisateur de $\sigma_i$, on a par ailleurs $t_{\alpha_j} = t_i$ pour $j = 1, \ldots, d$. Par suite, grâce au lemme **6.2**,

$$\Delta_{\mathcal{O},i} = \{\frac{t_i}{m_i}\alpha_1^\vee, \ldots, \frac{t_i}{m_i}\alpha_{d_i-1}^\vee, \frac{t_i}{\epsilon m_i}\alpha_d^\vee\}.$$

Ce sont des éléments de $\Lambda_{\mathcal{O}} = H_M(M^\sigma/M^1)$, puisque par construction ces éléments sont de la forme $H_M(h_\alpha^{t_\alpha})$, $\alpha \in \Delta_{\mathcal{O},\mu}$, alors que $h_\alpha^{t_\alpha} \in M^\sigma$ d'après **1.13**, en remarquant que le groupe $\mathrm{Stab}(\mathcal{O})$ opère d'après la définition **1.5** trivialement sur les $h_\alpha^{t_\alpha}$.

On pose $\Delta_{\mathcal{O},i}^\vee = \{\frac{m_i}{t_i}\alpha_1, \ldots, \frac{m_i}{t_i}\alpha_{d_i-1}, \frac{\epsilon m_i}{t_i}\alpha_{d_i}\}$. On a $\Delta_{\mathcal{O},i}^\vee \subseteq \Lambda_{\mathcal{O}}^\vee$: observons d'abord que les groupes $\mathrm{GL}_{m_i}/(GL_{m_i})^1$ sont cycliques. Le groupe $(\mathrm{GL}_{m_i})^{\sigma_i}/(\mathrm{GL}_{m_i})^1$ en est un sous-groupe d'indice $t_i$ (car $(\mathrm{GL}_{m_i}/(\mathrm{GL}_{m_i})^1)^{t_i} = (\mathrm{GL}_{m_i})^\sigma/$ $(\mathrm{GL}_{m_i})^1)$), et le groupe $M^\sigma/M^1$ est un produit de tels groupes. Il est clair que les caractères rationnels $m_i\alpha_j$, considérés comme éléments de $a_M^*$, envoient l'image de $M$ par $H_M$ dans $\mathbb{Z}$. Comme les éléments de $\Delta_{\mathcal{O},i}^\vee$ sont triviaux sur les facteurs $\mathrm{GL}_{m_j}$, $j \neq i$, il en résulte suit que les $\frac{m_i}{t_i}\alpha_j$ envoient $M^\sigma/M^1$ dans $\mathbb{Z}$, i.e. ce sont des éléments de $\Lambda_{\mathcal{O}}^\vee$.

On voit que le système des racines obtenu est du même type que $\Sigma_{\mathcal{O},\mu,i}$ sauf dans le cas où celui-ci est de type $C_d$. Alors, on trouve un système équivalent au système dual de $\Sigma_{\mathcal{O},\mu,i}$. Il est donc de type $B_d$. □

**7.** On va maintenant faire le lien avec les algèbres de Hecke avec paramètres.

**7.1** Rappelons d'abord la notion d'algèbre de Hecke avec paramètres définie dans [L1].

Soit $(\Lambda, \Lambda^\vee, \Sigma, \Sigma^\vee, \Delta)$ un quintuplet, où $\Lambda$ et $\Lambda^\vee$ sont des groupes abéliens libres de type fini en dualité par une application $\mathbb{Z}$-bilinéaire $\langle\cdot,\cdot\rangle : \Lambda \times \Lambda^\vee \to \mathbb{Z}$, $\Sigma \subseteq \Lambda$ un système de racines, $\Delta$ une base de $\Sigma$ et $\Sigma^\vee \subseteq \Lambda^\vee$ le système de racines dual de $\Sigma$, la dualité étant donné par $\langle\cdot,\cdot\rangle$.

Désignons par $W(\Sigma)$ le groupe de Weyl de $\Sigma$. Notons $\Sigma_1, \dots, \Sigma_r$ les composantes irréductibles de $\Sigma$ et $\Delta_i := \Delta \cap \Sigma_i$. Pour tout $i \in \{1, \dots, r\}$, supposons donné un ensemble $\{q_\alpha, \alpha \in \Delta_i\}$ de nombres réels $> 1$ tels que $q_\alpha = q_\beta$, si $\alpha$ et $\beta$ sont conjugués par un élément du groupe de Weyl de $\Sigma_i$. Si $\Sigma_i$ est de type $B_n$, on se donne en outre un nombre réel $q_i > 1$.

Lorsque $\alpha$ est dans $\Sigma$, désignons par $s_\alpha$ la symétrie élémentaire de $\Lambda$ associée à $\alpha$, et, pour $\alpha, \beta \in \Sigma$, par $m(\alpha,\beta)$ l'ordre de $s_\alpha s_\beta$.

Considérons le groupe $B_0(\Sigma)$ de générateurs $U_{s_\alpha}$, $\alpha \in \Delta$, avec les relations de tresse $U_{s_\alpha}U_{s_\beta}U_{s_\alpha}\cdots = U_{s_\beta}U_{s_\alpha}U_{s_\beta}\cdots$ ($m(\alpha,\beta)$ facteurs) pour tous $\alpha, \beta \in \Delta$. Remarquons qu'alors $U_{s_{\alpha_1}}\cdots U_{s_{\alpha_r}} = U_{s_{\beta_1}}\cdots U_{s_{\beta_r}}$, lorsque $s_{\alpha_1}\cdots s_{\alpha_r}$ et $s_{\beta_1}\cdots s_{\beta_r}$ sont des décompositions réduites d'un même élément $w$ de $W(\Sigma)$ [Sp, 8.3.3]. On peut donc poser $U_w := U_{s_{\alpha_1}}\cdots U_{s_{\alpha_r}}$. Notons $\mathcal{H}_0(\Sigma)$ le quotient de l'algèbre de groupe de $B_0(\Sigma)$ par l'idéal bilatère engendré par les éléments $(U_{s_\alpha}+1)(U_{s_\alpha}-q_\alpha)$, $\alpha \in \Delta$.

Notons $C = \mathbb{C}[\Lambda]$ l'algèbre de groupe de $\Lambda$ et $Z_\lambda$ l'élément de $C$ associé à $\lambda \in \Lambda$.

On appellera *algèbre de Hecke avec paramètres* $\{q_\alpha\}$ *et* $\{q_i\}$, et on notera $\mathcal{H} := \mathcal{H}(\Sigma,\ \{q_\alpha\}, \{q_i\})$ la $\mathbb{C}$-algèbre qui, en tant que $\mathbb{C}$-espace vectoriel est engendrée par $U_wZ_\lambda$, $w$ et $\lambda$ parcourant respectivement $W(\Sigma)$ et $\Lambda$. La multiplication dans $\mathcal{H}(\Sigma, \{q_\alpha\}, \{q_i\})$ est déduite de celle dans $\mathcal{H}_0(\Sigma)$ et de celle dans $C$ avec la règle de commutation

$$Z_\lambda U_{s_\alpha} - U_{s_\alpha} Z_{s_\alpha\lambda} = \begin{cases} (q_\alpha - 1)\frac{Z_\lambda - Z_{s_\alpha(\lambda)}}{1-Z_{-\alpha}}, \text{ si } \alpha^\vee \notin 2\Lambda^\vee, \\ (q_\alpha - 1 + Z_{-\alpha}((q_\alpha q_i)^{1/2} - (q_\alpha q_i^{-1})^{1/2}))\frac{Z_\lambda - Z_{s_\alpha(\lambda)}}{1-Z_{-2\alpha}}, \text{sinon;} \end{cases}$$

pour $\alpha \in \Delta$ et $\lambda \in \Lambda$.

Ceci est bien défini parce que $\alpha^\vee \in 2\Lambda^\vee$ est équivalent à dire que la composante irréductible $\Sigma_i$ de $\Sigma$ à laquelle $\alpha$ appartient est de type $B_n$ et que $\alpha$ est la racine courte dans $\Delta \cap \Sigma_i$.

Notons $\mathcal{Z}$ le centre de $\mathcal{H}$, $K(\mathcal{Z})$ son corps des fractions, et posons $\mathcal{H}_{K(\mathcal{Z})} = \mathcal{H} \otimes_{\mathcal{Z}} K(\mathcal{Z})$. Suivant [L1, 3.12], on a $\mathcal{H}_{K(\mathcal{Z})} = \bigoplus_{w\in W(\Sigma)} K(C)U_w$, où $K(C)$ désigne le corps des fractions de $C$.

**7.2** On va maintenant revenir à notre algèbre $\mathrm{End}_G(i_P^G E_{B_\mathcal{O}})$ et établir des propriétés des opérateurs $T_w$ pour montrer qu'il s'agit bien d'une algèbre de Hecke avec paramètres.

**Lemme:** *Soient $\alpha, \alpha' \in \Delta_\mathcal{O}$, $s = s_\alpha$, $s' = s_{\alpha'}$, et $m \in \mathbb{Z}$ tels que $(ss')^m = 1$. Supposons ou que $b_s = b_{s'} = 0$, ou que $s$ et $s'$ commutent, ou bien que $\alpha$ et $\alpha'$ engendrent un système de racines de type $B_2$, que $\alpha'$ soit la racine courte et que $b_s = 0$. Alors $T_sT_{s'}T_s \cdots = T_{s'}T_sT_{s'} \cdots$, le nombre de facteurs de chaque côté étant $m$.*

*Remarque:* Il semblerait que les conclusions du lemme ci-dessus deviennent fausses, si on omet une des hypothèses.

*Preuve:* Le système de racines $\Sigma_1$ engendré par $\alpha$ et $\alpha'$ est réduit de rang 2. S'il est réductible, il est de type $A_1 \times A_1$, et $J_s$ et $J_{s'}$ commutent. Mais, alors $T_s$ et $T_{s'}$ commutent également, et l'assertion est triviale.

Sinon, $\Sigma_1$ est de type $A_2$, $B_2$ ou $G_2$. Sans perte de généralité, on peut supposer que la racine $\alpha'$ est au plus aussi longue que $\alpha$. Notons $m$ l'ordre de $s_\alpha s_{\alpha'}$. On pose $q_\alpha = q^{a_s}$, $q_{\alpha'} = q^{a_{s'}+b_{s'}}$ et $q_0 = q^{a_{s'}-b_{s'}}$. Considérons le quintuplet $(\Lambda_\mathcal{O}, \Lambda_\mathcal{O}^\vee, \Sigma_1, \Sigma_1^\vee, \{\alpha, \alpha'\})$ et l'algèbre de Hecke à paramètres $\mathcal{H} = \mathcal{H}(\Sigma_1, \{q_\alpha, q_{\alpha'}\}, \{q_0\})$ associée. Remarquons que $\mathbb{C}[\Lambda_\mathcal{O}]$ s'identifie à $B_\mathcal{O}$. Notons $\mathcal{Z}$ le centre de $\mathcal{H}$. L'algèbre $\mathcal{H}_{K(\mathcal{Z})}$ est donc un $K(B_\mathcal{O})$-module.

Avec les notations de **7.1**, posons

$$j_\alpha = \frac{q_\alpha - Z_\alpha}{1 - Z_\alpha} \qquad \text{et} \qquad j_{\alpha'} = \frac{(q_{\alpha'}^{1/2} q_0^{1/2} - Z_{\alpha'})(Z_{\alpha'} + q_{\alpha'}^{1/2} q_0^{-1/2})}{1 - Z_{\alpha'}^2}.$$

Pour $P \in C$, $\beta \in \{\alpha, \alpha'\}$ et $t = s_\beta$, les règles de commutation impliquent que

$$P(U_\beta - q_\beta) - (U_\beta - q_\beta)\, {}^tP = ({}^tP - P) j_\beta. \qquad (*)$$

Si $\beta \in \Sigma_1$ est conjugué à $\alpha$ (resp. $\alpha'$), posons $j_\beta = j_\alpha(Y_\beta)$ (resp. $j_\beta = j_{\alpha'}(Y_\beta)$).

Dans $\mathcal{H}_{K(\mathcal{Z})}$, posons $S_s = U_s - q_\alpha + j_\alpha$, $S_{s'} = U_{s'} - q_{\alpha'} + j_{\alpha'}$, et, lorsque $w \in W(\Sigma_1)$ et que $w = s_{\alpha_1} \cdots s_{\alpha_r}$ est une décomposition réduite de $w$, $S_w = S_{\alpha_1} \cdots S_{\alpha_r}$. Ceci ne dépend pas de la décomposition réduite choisie [R, 4.3]. Alors, pour tout

$w, w' \in W(\Sigma_1)$,

$$S_w S_{w'} = (\prod_\beta j_\beta j_{-\beta}) S_{ww'},$$

le produit portant sur l'ensemble des $\beta \in \Sigma_1^+$ tels que $w^{-1}\beta < 0$ et ${w'}^{-1}w^{-1}\beta > 0$ [R, 4.3].

On va montrer qu'il existe un unique homomorphisme d'algèbres $\mathcal{H}_{K(\mathcal{Z})} \to \mathrm{End}_G(i_P^G E_{B_\mathcal{O}})_{K(\mathcal{Z}_\mathcal{O})}$ qui envoie $S_s$ sur $R_s$, $S_{s'}$ sur $R_{s'}$ et qui induit l'identité sur $K(B_\mathcal{O})$. Il en suivra la relation $T_s T_{s'} T_s \cdots = T_{s'} T_s T_{s'} \cdots$, le nombre de facteurs de chaque côté étant égal à l'ordre $ss'$, puisque $U_s \mapsto T_s$ et que $U_{s'} \mapsto T_{s'}$.

Comme les opérateurs $S_s$ et $S_{s'}$ engendrent ensemble avec $K(B_\mathcal{O})$ l'algèbre $\mathcal{H}_{K(\mathcal{Z})}$, l'unicité est évidente. Il reste à prouver l'existence, i.e. que $R_s, R_{s'}$ vérifient les mêmes relations que $S_s, S_{s'}$.

Les relations de commutation (*) pour les opérateurs $R_s$ et $R_{s'}$ résultent directement de leur définition **5.2** avec le lemme **4.8**. Il en est de même des identités $R_s^2 = j_\alpha j_{-\alpha}$ et $R_{s'}^2 = j_{\alpha'} j_{-\alpha'}$ avec le lemme **4.7**. Il reste à vérifier que, lorsque $w \in W(\Sigma_1)$ et que $w = s_{\alpha_1} \cdots s_{\alpha_r}$ est une décomposition en symétries simples de $w$, alors $R_{s_{\alpha_1}} \cdots R_{s_{\alpha_r}}$ est indépendant du choix de cette décomposition. Pour cela, on va étudier les différents cas. En fait, si $m(s, s')$ désigne l'ordre de $ss'$, il suffit de montrer que $R_s R_{s'} R_s \cdots = R_{s'} R_s R_{s'} \cdots$, le nombre de facteurs de chaque côté étant $m(s, s')$ [Sp, 8.3.3].

Si $\Sigma_1$ est de type $A_2$ ou $G_2$, alors $b_{s'} = 0$ par hypothèse et $c_{s'} = c_s$, puisque $s$ et $s'$ sont conjugués. Les relations ci-dessus sont alors une conséquence directe du corollaire **3.4**.

Supposons maintenant $\Sigma_1$ de type $B_2$. Il faut montrer que $R_{s_1} R_{s_2} R_{s_1} R_{s_2} = R_{s_2} R_{s_1} R_{s_2} R_{s_1}$. Comme chaque facteur $R_{s_i}$ apparaît avec la même multiplicité, les facteurs scalaires sont égaux. On est donc ramené à l'égalité $X_\alpha J_{s_\alpha} X_{\alpha'} J_{s_{\alpha'}} X_\alpha J_{s_\alpha} X_{\alpha'} J_{s_{\alpha'}} = X_{\alpha'} J_{s_{\alpha'}} X_\alpha J_{s_\alpha} X_{\alpha'} J_{s_{\alpha'}} X_\alpha J_{s_\alpha}$ qui résulte du lemme **2.7**. □

**7.3 Proposition:** *Lorsque $\Sigma_i$ est une composante irréductible de $\Sigma_\mathcal{O}$ et que $\alpha \in \Delta \cap \Sigma_i$, alors $b_{s_\alpha} \neq 0$ implique que $\Sigma_i$ est de type $B_n$ et que $\alpha$ correspond à la racine courte de $B_n$.*

*Preuve:* Il résulte des travaux de Bernstein-Zelevinsky [BZ2, 4.2] et, dans le cas d'une algèbre simple, de [T] que $b_{s_\alpha} = 0$ si la fonction $\mu^{M_\alpha}$ est égale à celle pour un groupe linéaire général ou le groupe multiplicatif d'une algèbre simple. Ceci est également vrai pour $\mathrm{SL}_2(F)$ et $\mathrm{PGL}_2(F)$. Or, dans les autres cas, $\Sigma_i$ est, d'après **1.13** et **6.1**, de type $B_n$ et $\alpha$ est la racine courte. □

**7.4** Lorsque $\Sigma_i$ est une composante irrédictible de $\Sigma_\mathcal{O}$, désignons par $\alpha_i$ l'unique racine dans $\Delta \cap \Sigma_i$ telle que $\alpha_i^\vee \in 2\Lambda_\mathcal{O}^\vee$, si une telle racine existe, et posons $a_i = a_{s_{\alpha_i}}$ et $b_i = b_{s_{\alpha_i}}$.

**Proposition:** *L'algèbre* $\bigoplus_{w\in W_{\mathcal{O}}} B_{\mathcal{O}}\, T_w$ *est une algèbre de Hecke avec paramètres* $\{q^{a_{s_\alpha}+b_{s_\alpha}}\}$ *et* $\{q^{a_i-b_i}\}$ *et donnée radicielle basique* $(\Lambda_{\mathcal{O}}, \Sigma_{\mathcal{O}}, \Lambda_{\mathcal{O}}^\vee, \Sigma_{\mathcal{O}}^\vee, \Delta_{\mathcal{O}})$*, i correspondant aux différentes composantes irréductibles de* $\Sigma_{\mathcal{O}}$*.*

*Preuve:* La seule propriété qui reste à vérifier est la règle de commutation. C'est un calcul élémentaire. □

**7.5** Remarquons que les valeurs possibles pour les paramètres $a_s$ et $b_s$ sont connues pour les groupes considérés ici par l'intermédiaire de la correspondance de Langlands: l'ordre du groupe des caractères non ramifiés qui stabilisent la classe d'isomorphie d'une représentation irréductible cuspidale de $\mathrm{GL}_n(F)$ est un diviseur de $n$. (Ceci résulte de la correspondance locale de Langlands [HT].) On a un résultat similaire dans le cas des algèbres simples par la correspondance de Jacquet-Langlands [DKV].

Dans le cas où $G$ est égal à $\mathrm{GL}_n(F)$ ou le groupe multiplicatif d'une algèbre simple, on en déduit les valeurs de $a_s$ et $b_s$ directement des travaux cités dans la preuve de **7.3**. Pour $H_k$ un groupe symplectique ou orthogonal et $\sigma\otimes\tau$ une représentation irréductible cuspidale de $\mathrm{GL}_n(F)\times H_k$ telle que $\mu(\sigma\otimes\tau)=0$, C. Moeglin [M] a récemment déduit des travaux de J. Arthur que le nombre $a>0$ tel que $\sigma|\det_n|_F^a\otimes\tau$ soit un pôle de $\mu$ est un demi-entier. Par ailleurs, si $\chi$ est un caractère unitaire non ramifié tel que la fonction $\mu$ s'annulle également en $\sigma\chi\otimes\tau$, alors le nombre réel $b>0$ tel que $(\sigma\chi)|\det|_F^b\otimes\tau$ soit un pôle de $\mu$ multiplié par l'ordre du stabilisateur $t$ de $\sigma$ est un entier si et seulement si $ta$ l'est. A part cela, à moins que $\sigma$ et $\sigma\chi$ soient isomorphes, la valeur de $b$ n'est en général pas conditionnée par celle de $a$. Le lecteur pourra consulter [H3] où ceci est rendu complètement explicite et appliqué aux algèbres de Hecke avec paramètres obtenus ici.

Déterminer les valeurs des paramètres $a_s$ et $b_s$ directement sur $\sigma$, par exemple par la théorie des types, semble être un problème très dur et difficile (voir par exemple [B] pour un résultat partiel).

**7.6 Proposition:** *Soit* $w\in W_{\mathcal{O}}$ *et* $r\in R(\mathcal{O})$*. Alors* $r^{-1}wr\in W_{\mathcal{O}}$ *et* $T_wJ_r = J_rT_{r^{-1}wr}$*.*

*Preuve:* Il suffit de considérer le cas où $w$ est une symétrie simple $s_\alpha$ et de montrer que $r^{-1}s_\alpha r$ est également une symétrie simple. Or, comme par définition $r$ laisse $\Sigma_{\mathcal{O}}\cap\Sigma(P)$ invariant, il en est de même pour $\Delta_{\mathcal{O}}$, et, par suite, $r^{-1}s_\alpha r = s_{r\alpha}$ est une symétrie simple. □

**7.7 Théorème:** *L'algèbre* $\mathrm{End}_G(i_P^G E_{B_{\mathcal{O}}})$ *est isomorphe au produit semi-direct* $\mathbb{C}[R\ (\mathcal{O})]\ltimes\mathcal{H}(\Sigma_{\mathcal{O}}, \{q^{a_{s_\alpha}+b_{s_\alpha}}\}, \{q^{a_i-b_i}\})$ *de l'algèbre de groupe de* $R(\mathcal{O})$ *avec l'algèbre de Hecke avec paramètres* $\{q^{a_{s_\alpha}+b_{s_\alpha}}\}$ *et* $\{q^{a_i-b_i}\}$*.*

*Preuve:* C'est une conséquence immédiate de **5.10**, **7.4** et de **7.6**. □

**7.8 Corollaire:** *La catégorie* $Rep_G({}^W\mathcal{O})$ *est isomorphe à la catégorie des modules à droite sur l'algèbre* $\mathbb{C}[R(\mathcal{O})] \ltimes \mathcal{H}(\Sigma_{\mathcal{O}}, \{q^{a_{s_\alpha}+b_{s_\alpha}}\}, \{q^{a_i-b_i}\})$.

*Preuve:* C'est une conséquence immédiate de **7.7** combiné avec le théorème de Bernstein [Ro, 1.6] mentionné dans l'introduction. □

**7.9** *Remarque:* L'isomorphisme de catégories de **7.8** est compatible avec l'induction parabolique dans le sens suivant: remarquons d'abord que le théorème **7.7** et le corollaire **7.8** se généralisent facilement aux sous-groupes de Levi des groupes considérés ici. En effet, un tel sous-groupe de Levi est un produit direct de groupes de ce type, et les algèbres considérées sont donc des produits des algèbres associées aux facteurs.

Soit alors $P' = M'U'$ un sous-groupe parabolique de $G$ tel que $P' \supseteq P$ et $M' \supseteq M$. Pour $\alpha \in \Delta_{\mathcal{O}}^{M'}$ et $r \in R^{M'}(\mathcal{O})$, notons $T_{s_\alpha}^{M'}$ et $J_r^{M'}$ les opérateurs $T_{s_\alpha}$ et $J_r$ définis relatifs à $M'$. On a une injection canonique de $\mathrm{End}_{M'}(i_{P\cap M'}^{M'} E_{B_{\mathcal{O}}})$ dans $\mathrm{End}_G(i_P^G E_{B_{\mathcal{O}}})$ défini par le foncteur d'induction parabolique $i_{P'}^G$. On observe alors que cette inclusion correspond à l'inclusion naturelle des algèbres de Hecke avec paramètres respectives (ainsi que des algèbres de groupe fini). En effet, cette inclusion de $\mathrm{End}_{M'}(i_{P\cap M'}^{M'} E_{B_{\mathcal{O}}})$ dans $\mathrm{End}_G(i_P^G E_{B_{\mathcal{O}}})$ envoie un opérateur $T_{s_\alpha}^{M'}$, $\alpha \in \Delta_{\mathcal{O}}^{M'}$, sur $T_{s_\alpha}$ et un opérateur $J_r^{M'}$, $r \in R^{M'}(\mathcal{O})$, sur $J_r$, la construction de ces opérateurs étant compatible avec l'induction parabolique (pour les opérateurs d'entrelacement voir [W, IV.1 (14)]) et se faisant par rapport à un sous-groupe de Levi contenu dans $M'$.

Rappelons, comme déjà remarqué dans l'introduction, que l'isomorphisme de catégorie entre $Rep_G({}^W\mathcal{O})$ et la catégorie des $\mathrm{End}_G(i_P^G E_{B_{\mathcal{O}}})$-modules à droite est lui-aussi compatible avec l'induction parabolique [Ro, 5.3]. L'isomorphisme de catégorie **7.8** est donc compatible avec l'induction parabolique.

On voit de même qu'il est compatible avec le foncteur de Jacquet: en effet, le foncteur de Jacquet $r_{P'}^G$ définit un foncteur de la catégorie $Rep_G({}^W\mathcal{O})$ à la catégorie $Rep_{M'}({}^{W^{M'}}\mathcal{O})$. D'autre part, par l'inclusion de $\mathrm{End}_{M'}(i_{P\cap M'}^{M'} E_{B_{\mathcal{O}}})$ dans $\mathrm{End}_G(i_P^G E_{B_{\mathcal{O}}})$, tout $\mathrm{End}_G(i_P^G E_{B_{\mathcal{O}}})$-module à droite a une structure de $\mathrm{End}_{M'}$ $(i_{P\cap M'}^{M'} E_{B_{\mathcal{O}}})$-module à droite, ce qui induit un foncteur entre ces deux catégories de modules. Par [Ro, 5.3], celui-ci commute avec le foncteur de Jacquet.

Clermont Université, Université Blaise-Pascal, Laboratoire de mathématiques, BP 10448, 63000 Clermont-Ferrand, France, and CNRS, UMR 6620, Laboratoire de Mathématiques, BP 80026, 63177 Aubière, France

*E-mail address*: heiermann@math.univ-bpclermont.fr

# INTERTWINING OPERATORS AND HECKE ALGEBRAS WITH PARAMETERS OF A P-ADIC REDUCTIVE GROUP: THE CASE OF CLASSICAL GROUPS

VOLKER HEIERMANN
(TRANSLATED BY YUJIE XU)

ABSTRACT. For $G$, a symplectic or orthogonal $p$-adic group (not necessarily split) or an inner form of a general linear $p$-adic group, we compute the endomorphism algebras of some induced projective generators á la Bernstein of the category of smooth representations of $G$ and show that these algebras are isomorphic to the semi-direct product of a Hecke algebra with parameters by a finite group algebra. Our strategy and parts of our intermediate results apply to a general reductive connected $p$-adic group.

*In memory of Jacqueline*

Fix a nonarchimedean local field $F$. Let $G$ be (the group of points of) a connected reductive group defined over $F$. It was shown by J. Bernstein [Ber84] that the category $\mathrm{Rep}(G)$ of complex smooth representations of $G$ is the direct product of certain full subcategories, denoted below as $\mathrm{Rep}_G({}^W\mathcal{O})$. These are generally conjectured to be isomorphic to categories of right modules over the Hecke algebras with parameters (or more generally over the *extended* Hecke algebra, i.e. the product of such an algebra with an algebra of finite group, possibly twisted by a 2-cocycle). The interest of this comes from the fact that the representation theory of Hecke algebras with parameters has been extensively studied by G. Lusztig and others. (Let's cite for example [Kat09], [Lus02] and [OS10].)

One approach to prove this equivalence of category is the "theory of types" of C.Bushnell and Ph. Kuszko [BK98]. On another side, J. Bernstein made explicit in 1992 in a course at Harvard university certain representations $\Pi_{\mathcal{O}}$ of $G$ which possesses the property that $\mathrm{Rep}_G({}^W\mathcal{O})$ is isomorphic to the category of right modules over the intertwining algebra $\mathrm{End}_G(\Pi_{\mathcal{O}})$.

The goal of this work is to calculate the algebra $\mathrm{End}_G(\Pi_{\mathcal{O}})$ explicitly to show that it is indeed a Hecke algebra with extended parameters. For that, I restrict myself to the case where $G$ is a symplectic or orthogonal group (split or not) or again where $G$ is an inner form of $\mathrm{GL}_n(F)$. We shall remark that, except in the case of an inner form of $\mathrm{GL}_n(F)$, this equivalence of category is to the knowledge of the author not yet accessible by the theory of types.

In fact, the hypotheses that I need should also be satisfied for other forms of these classical groups. Moreover, the approach and a good part of these intermediary results are general. Certain proofs were moreover motivated by obtaining a general result and could well be simplified under the hypotheses imposed in this paper.

Note that the parameters of these Hecke algebras with parameters obtained here have been explained in [Hei10] (see also section 7.5).

Let's be more precise for now: let $P = MU$ be a parabolic subgroup of $G$ and $(\sigma, E)$ be a cuspidal irreducible representation of $M$. Denote by $\mathfrak{X}^{\mathrm{nr}}(M)$ the group of unramified characters of $M$, $W$ the Weyl group of $G$ defined relative to a maximal split torus contained in $M$ and $\mathcal{O}$ the set of equivalence classes of representations of the form $\sigma \otimes \chi$, $\chi \in \mathfrak{X}^{\mathrm{nr}}(M)$. Thus, ${}^W\mathcal{O}$ is the orbit of $\mathcal{O}$ for the action by $W$, ${}^W\mathcal{O} = \{{}^w\mathcal{O} | w \in W\}$.

Denote by $i_P^G$ the functor of normalized parabolic induction and by $r_P^G$ its left adjoint. The category $\mathrm{Rep}_G({}^W\mathcal{O})$ is thus the full subcategory of the category $\mathrm{Rep}(G)$ of complex smooth representations of $G$ whose objects are the representations $\pi$ which satisfy the following property: the set of equivalence classes of irreducible subquotients of $r_{P'}^G\pi$ is contained in ${}^W\mathcal{O}$ if $P'$ is associated to $P$, and it does not contain any irreducible cuspidal subquotient otherwise.

Denote by $M^1$ the intersection of kernels of unramified characters of $M$ and $(\sigma_1, E_1)$ an irreducible component of the restriction of $\sigma$ to $M^1$. Denote by $\mathrm{ind}_{M^1}^M$ the functor of compact induction. Therefore, it was shown by Bernstein [Roc02, 1.6] that the category $\mathrm{Rep}_G({}^W\mathcal{O})$ is isomorphic to the category of right $\mathrm{End}_G(i_P^G \mathrm{ind}_{M^1}^M E_1)$-modules. Note that $\mathrm{ind}_{M^1}^M E_1$ does not depend on a choice of $(\sigma_1, E_1)$. This isomorphism of categories is moreover compatible with the parabolic induction and the Jacquet functor [Roc02, 5.3].

The object of study of this work is thus the algebra $\mathrm{End}_G(i_P^G \mathrm{ind}_{M^1}^M E_1)$. For this I restrict myself to the case where $G$ is a symplectic or orthogonal group (split or not) or again where $G$ is an inner form of $\mathrm{GL}_n(F)$. In this situation, I give a basis for this algebra as a module over a ring of regular functions on $\mathcal{O}$ (cf. proposition 7.4) to deduce that this algebra is isomorphic to semi-direct product of a Hecke algebra with parameters (in the sense of Lusztig [Lus89]) with the algebra of a finite group (cf. theorem 7.7). (The finite group is a kind of "$R$-group".) We finish our article with a few remarks saying that the isomorphism of categories which is deduced from our isomorphism of algebras is compatible with the parabolic induction and the Jacquet functor (cf. 7.9).

The prior development of the article is as follows: paragraph 1 is used to fix some definitions and notations and to establish or recall basic results on standard intertwining operators and on smooth representations of a $p$-adic group. Paragraph 2 is a preliminary to paragraph 3 where we deduce standard intertwining operators from Harish-Chandra generators of the intertwining algebra of the "large" representation $i_P^G \operatorname{ind}_{M^1}^M E_{|M^1}$. However these have poles. In particular, they only generate algebra on the field of rational functions of $\mathcal{O}$. In paragraph 4, we specialize in the intertwining algebra of $i_P^G \operatorname{ind}_{M^1}^M E_1$, and we explain the generators which still have poles. In paragraph 5, we define the operators $T_w$, without poles and motivated by certain relations that the generators of a Hecke algebra with parameters should satisfy, to then show that these operators generate the algebra $\operatorname{End}_G(i_P^G \operatorname{ind}_{M^1}^M E_1)$ (possibly together with some of the intertwining operators above related to "$R$-group"). We remark that the proof that these operators $T_w$ generate the algebra $\operatorname{End}_G(i_P^G \operatorname{ind}_{M^1}^M E_1)$ is crucial in this work. On this occasion, a property of intertwining operators proven in [Hei01] will be useful. Paragraph 6 is used to define the root data that will be attached to the description of $\operatorname{End}_G(i_P^G \operatorname{ind}_{M^1}^M E_1)$ as Hecke algebra with parameters. In paragraph 7, the definition of Lusztig's Hecke algebra with parameters will be recalled before establishing the last indispensable properties of these operators $T_w$.

The author thanks J.-L. Waldspurger for having indicated to him in the example of $\mathrm{SL}_2$ the possibility of finding Hecke algebras with parameters from intertwining operators. He also benefited from discussions with P. Schneider and E.-W. Zink on different aspects of this article, and he thanks G. Henniart and J.-L. Waldspurger for remarks on a preliminary version of this paper. The critiques of different reviewers ultimately helped improve the final version.

## 1.

We keep the hypotheses and notations of the introduction. In addition, we fix a minimal parabolic subgroup $P_0$ contained in $P$, a Levi subgroup $M_0$ of $P_0$ contained in $M$, $P_0 = M_0U_0$, and a maximal split (over $F$) torus $A_0$ of $M_0$. The parabolic subgroup $P$ is thus a standard parabolic subgroup (i.e. $P \supseteq P_0$), and its Levi subgorup $M$ is a standard Levi (i.e. $M \supseteq M_0$). The maximal split torus contained in the center of $M$ will be denoted by $A_M$. We denote $W = W^G (= N_G(A_0)/Z_G(A_0))$ the Weyl group of $G$ defined relative to $A_0$, and we equip it with the length fixed by the choice of $P_0$. We will denote this order $\ell^G$ or simply $\ell$. We denote by $K$ a maximal compact subgroup of $G$ which is in good position relative to $A_0$ by Bruhat-Tits theory. Moreover, when $H$ is a group, we will write $X(H)$ for the group $\operatorname{Hom}(H, \mathbb{C}^\times)$.

The set of nontrivial roots of $A_M$ in the Lie algebra of $G$ will be designated by $\Sigma(A_M)$, the subset of roots which act on the Lie algebra of $U$ by $\Sigma(P)$, and the set of roots reduced by $\Sigma_{\mathrm{red}}(P)$. We have a bijection $\alpha \mapsto M_\alpha$ between $\Sigma_{\mathrm{red}}(P)$ and the set of Levi subgroups of $G$ which contain $M$ and which are minimal for this property.

Moreover, we denote by $a_M$ the real Lie algebra of $A_M$, by $a_M^*$ its dual, by $a_{M,\mathbb{C}}^*$ its complexification, by $|\cdot|_F$ the modulus of $F$, by $q$ the cardinality of residue field of $F$ and by $H_M$ te map $M \to a_M$ which satisfies $q^{-\langle H_M(m),\alpha\rangle} = |\alpha(m)|_F$ for every rational character $\alpha$ of $M$ and every $m \in M$. The kernel of $H_M$ is equal to $M^1$. If $s$ is a complex number, we denote $\chi_{\alpha\otimes s}$ the unramified character $m \mapsto |\alpha(m)|^s$ of $M$. This map extends to a homomorphism of groups $a_{M,\mathbb{C}}^* \to \mathfrak{X}^{\mathrm{nr}}(M)$ which is surjective.

1.1

We denote by $W(M)$ the set of representatives in $W$ of elements of quotient group $\{w \in W | w^{-1}Mw = M\}/W^M$ which are of minimum length in their class modulo $W^M$ on the right. Observe that $W(M)$ is a subgroup of $W$. (This is also the subgroup of $W$, formed of elements which fix $M$ and $P_0 \cap M$.) We denote by $W(M, \mathcal{O})$ the subgroup of $W(M)$ formed of elements that stabilize $\mathcal{O}$.

For $w \in W(M)$, put $\ell_M(w) := |\Sigma_{\mathrm{red}}(P) \cap \Sigma_{\mathrm{red}}(w\overline{P}w^{-1})|$.

**Proposition** *For two elements $w$ and $w'$ of $W(M)$ satisfying $\ell_M(ww') = \ell_M(w') + \ell_M(w)$, it is necessary and sufficient that $\ell(ww') = \ell(w) + \ell(w')$.*

*Proof.* First note that, for every $w \in W(M)$, $w\Sigma_{\mathrm{red}}(M \cap P_0) = \Sigma_{\mathrm{red}}(M \cap P_0)$ and that $\ell(w) = |\Sigma_{\mathrm{red}}(P_0) \cap \Sigma_{\mathrm{red}}(w\overline{P_0})|$. First it follows that $\Sigma^M(A_0) \cap \Sigma_{\mathrm{red}}(P_0) \cap \Sigma_{\mathrm{red}}(w\overline{P_0}) = \varnothing$. Therefore, an element $\alpha$ of $\Sigma(A_0)$ is in $\Sigma_{\mathrm{red}}(P_0) \cap \Sigma_{\mathrm{red}}(w\overline{P_0})$, if and only if $\alpha|_{A_M} \in \Sigma(P) \cap \Sigma(w\overline{P})$. Consequently, for $w, w' \in W(M)$, $\ell_M(ww') = \ell_M(w) + \ell_M(w')$ equals

$$|\Sigma_{\mathrm{red}}(P_0) \cap \Sigma_{\mathrm{red}}(ww'\overline{P_0})| = |\Sigma_{\mathrm{red}}(P_0) \cap \Sigma_{\mathrm{red}}(w\overline{P_0})| + |\Sigma_{\mathrm{red}}(P_0) \cap \Sigma_{\mathrm{red}}(w'\overline{P_0})|,$$

thence the proposition by the expression for $\ell(w)$ recalled at the beginning. □

For every $w \in W$, denote $P(w) = M(w)U(w)$ the minimal standard parabolic subgroup such that $P, wP \subseteq P(w)$. (If $M_\alpha$ is a standard Levi subgroup and $w$ the nontrivial element of $W^{M_\alpha}(M)$, thus $M(w) = M_\alpha$.) We can choose representatives $\overline{w}$ of elements of $W$ in $G$ such that $\overline{w} \in K \cap M(w)$, which is what we will do from now on. In the following we will identify the elements of $W$ with their representatives in $K$. It will always be clear from the context, whether $w$ denotes an element of $W$ or of $G$, by making the following conventions: for $w, w' \in W$, the symbol $ww'$, considered as an element of $K$, denotes the product of the representative of $w$ in $K$ with that of $w'$ and **not** the representative of $ww'$. Moreover, the symbol $w^{-1}$ corresponds to the inverse (of the representative) of $w$ in the group $G$ which is indeed an element of $K \cap M_w$.

Note that our construction of operators of $\mathrm{End}_G(i_P^G E_B)$ will essentially be independent of the choice of the set of representatives of $W$. We will however in **1.15** impose some additional restrictions on these representatives.

When $P_1 = M_1U_1$ is a semi-standard parabolic subgroup of $G$ (i.e. $M_1$ contains $A_0$), $(\pi, V)$ a smooth representation of $M_1$ and $w \in W$, we denote by $\lambda(w)$ the isomorphism $i_{P_1}^G V \rightarrow i_{wP_1}^G wV, v \mapsto v(w^{-1}\cdot)$, which depends on the choice of the representative of $w$ in $G$. We can also define $\lambda(g)$ for every element $g$ of $G$. If $g \in M_1$, then $\lambda(g)$ is an isomorphism $i_{P_1}^G V \rightarrow i_{P_1}^G gV$ which is induced by functoriality of the isomorphism $\pi(g^{-1}) : V \rightarrow gV$. We will also write $i_{P_1}^G(\pi(g^{-1}))$ (not to confuse with $i_P^G\pi(g^{-1})$).

1.2

Recall that $\mathcal{O}$ is the quotient of a complex torus by a finite group. The function $\mu$ of Harish-Chandra is a rational function defined over $\mathcal{O}$. We write $\mu^H$, if it is defined with respect to a reductive subgroup $H$ of $G$ which contains $M$. We omit this exponent, if $G = H$.

**Theorem** [Sil79, 5.4.2.2 and 5.4.2.3] (Harish-Chandra) *Let $\alpha \in \Sigma(P)$ and let $\sigma_1$ be a cuspidal irreducible representation of $M$.*
*(a) If $\mu^{M_\alpha}(\sigma_1) = 0$, thus there exists a unique nontrivial element $s_\alpha$ in $W^{M_\alpha}(M)$ such that $s_\alpha(P \cap M_\alpha) = \overline{P} \cap M_\alpha$ and $s_\alpha\sigma_1 \simeq \sigma_1$.*
*(b) Suppose that there exists a unique nontrivial element $s_\alpha$ in $W^{M_\alpha}(M)$ such that $s_\alpha(P \cap M_\alpha) = \overline{P} \cap M_\alpha$ and $s_\alpha\sigma_1 \simeq \sigma_1$. Thus, in order that $\mu^{M_\alpha}(\sigma_1) \neq 0$, it is necessary and sufficient that the representation $i_{P\cap M_\alpha}^{M_\alpha}\sigma_1$ is reducible. The representation $i_{P\cap M_\alpha}^{M_\alpha}\sigma_1$ is thus direct sum of two non-isomorphic irreducible representations.*

1.3.

Denote by $a_M^{M_\alpha}$ the subspace of $a_M$ which is orthogonal to $a_{M_\alpha}^*$.

**Proposition** *The set $\Sigma_{\mathcal{O},\mu} := \{\alpha \in \Sigma_{\mathrm{red}}(A_M) | \mu^{M_\alpha}$ has a zero on $\mathcal{O}\}$ is a root system. For $\alpha \in \Sigma_{\mathcal{O},\mu}$, denote by $s_\alpha$ the unique element of $W^{M_\alpha}(M, \mathcal{O})$ which conjugates $P \cap M_\alpha$ and $\overline{P} \cap M_\alpha$. The Weyl group $W_{\mathcal{O}}$ of $\Sigma_{\mathcal{O},\mu}$ identifies with the subgroup of $W(M, \mathcal{O})$ generated by the reflections $s_\alpha$. For every $\alpha \in \Sigma_{\mathcal{O},\mu}$, denote by $\alpha^\vee$ the unique element of $a_M^{M_\alpha}$ which satisfies $\langle \alpha, \alpha^\vee \rangle = 2$. Thus $\Sigma_{\mathcal{O},\mu}^\vee := \{\alpha^\vee | \alpha \in \Sigma_{\mathcal{O},\mu}\}$ is the set of coroots of $\Sigma_{\mathcal{O},\mu}$, the duality being that between $a_M$ and $a_M^*$.*

*The set $\Sigma(P) \cap \Sigma_{\mathcal{O},\mu}$ is the set of positive roots for a certain order on $\Sigma_{\mathcal{O},\mu}$.*

*Proof.* The proof of the first part of the proposition is analogous to that of the proposition **4.2** in [Hei04] (see also [Sil81]), after verifying that $\Sigma_{\mathcal{O}}$ is stable under the action of $W(M,\mathcal{O})$. This follows from the $W(M,\mathcal{O})$-invariant of the Harish-Chandra function $\mu$ [Wal03, V.2.1].

Let $\alpha \in \Sigma_{\mathcal{O},\mu}$. The sub-vector-space of $a_M^*$ generated by $\alpha$ is $a_M^{M_\alpha *}$. This one is the orthogonal of $a_{M_\alpha}$. Consequently, $\alpha^\vee$ must belong to $a_M^{M_\alpha}$. As $\langle \alpha, \alpha^\vee \rangle = 2$, this is thus indeed the coroot.

The last assertion comes from the fact that $\Sigma(P)$ defines an order on $\Sigma_{\text{red}}(A_M)$, thus on $\Sigma_{\mathcal{O},\mu}$. □

1.4

**Definition** We set $\Sigma_{\mathcal{O},\mu}(P) = \Sigma(P) \cap \Sigma_{\mathcal{O},\mu}$, and we denote by $\Delta_{\mathcal{O},\mu}$ the basis of $\Sigma_{\mathcal{O},\mu}$ determined by the order for which the set of positive roots is $\Sigma_{\mathcal{O},\mu}(P)$. The length over $W_{\mathcal{O}}$ will be denoted by $\ell_{\mathcal{O}}$.

We also fix a unitary representation $\sigma$ whose equivalence class belongs to $\mathcal{O}$ such that $\mu^{M_\alpha}(\sigma) = 0$ for every $\alpha \in \Delta_{\mathcal{O},\mu}$. (This is possible, since $\Delta_{\mathcal{O},\mu}$ is a basis and that $\mu^{M_\alpha}$ is invariant for the action by $\mathfrak{X}^{\text{nr}}(M_\alpha)$.)

1.5

For $\alpha \in \Sigma_{\mathcal{O},\mu}$, fix an element $h_\alpha$ of $M \cap M_\alpha^1$ such that $H_M(h_\alpha)$ is a multiple of $\alpha^\vee$ by a real number $> 0$ and that, for every $\chi \in \mathfrak{X}^{\text{nr}}(M)$,$\chi(h_\alpha) = 1$ equivalent to $\chi \in \mathfrak{X}^{\text{nr}}(M_\alpha)$. This is indeed defined, since the kernel of $H_M$ is equal to $M^1$ and that $M \cap M_\alpha^1/M^1$ is a free $\mathbb{Z}$-module of rank 1 (cf. [Hei04, 1.2]). We also see that $H_M(h_\alpha)$ is unique.

Denote by $\text{Stab}(\mathcal{O})$ the group of unramified characters $\chi$ of $M$, satisfying $\sigma \otimes \chi \simeq \sigma$ (it only depends on $\mathcal{O}$), and $t_\alpha$ the smallest integer $\geq 1$ such that $\chi(h_\alpha^{t_\alpha}) = 1$ equivalent to $\chi \in \mathfrak{X}^{\text{nr}}(M_\alpha)\,\text{Stab}(\mathcal{O})$. Denote by $b_{h_\alpha}$ the map $\mathfrak{X}^{\text{nr}}(M) \to \mathbb{C}$, $\chi \mapsto \chi(h_\alpha)$, and put $Y_\alpha = b_{h_\alpha}$ and $X_\alpha = Y_\alpha^{t_\alpha}$. (The reason for this double notation will become clear in the section 2.)

**Lemma** *For $w \in W(M)$, we have ${}^wY_\alpha = Y_{w\alpha}$. In particular, ${}^{s_\alpha}Y_\alpha = Y_\alpha^{-1}$.*

*Proof.* Put $H_M(h_\alpha) = m_\alpha \alpha^\vee$ with $m_\alpha > 0$. Remark that $Y_\alpha = b_{h_\alpha}$ and thus ${}^wY_\alpha = {}^wb_{h_\alpha} = b_{wh_\alpha w^{-1}}$. As $H_M(wh_\alpha w^{-1}) = m_\alpha(w\alpha^\vee)$, it remains to prove that $m_\alpha = m_{w\alpha}$. As, for $\lambda \in a_M^*$,

$$\chi_{w\lambda}(wh_\alpha w^{-1}) = {}^w\chi_\lambda(wh_\alpha w^{-1}) = \chi_\lambda(h_\alpha),$$

we have $\chi_{w\lambda}(wh_\alpha w^{-1}) = 1$ if and only if $\chi_\lambda \in \mathfrak{X}^{\text{nr}}(M_\alpha)$. But this is equivalent to $\chi_{w\lambda} = {}^w\chi_\lambda \in \mathfrak{X}^{\text{nr}}(M_{w\alpha})$. Thus, $H_M(h_{w\alpha}) = H_M(wh_\alpha w^{-1})$, and, consequently, $m_\alpha = m_{w\alpha}$. □

1.6

The following result was shown by A. Silberger [Sil80, 1.6] (see also the remark in the proof of the proposition **4.1** in [Hei04]). (Note that the statement of theorem **1.6** in [Sil80] is relative to a function denoted by $U(\sigma, \cdot)$, but this one is effectively the function $\mu$ (cf. last paragraph on the page 575 of [Sil80].)

**Proposition** *Let $\alpha \in \Sigma_{\mathrm{red}}(P)$ and $s = s_\alpha$. If $\mu^{M_\alpha}$ is not constant, thus $\alpha \in \Sigma_{\mathcal{O},\mu}$, and we can find a constant $c'_s > 0$ and real numbers $a_s \geq 0$ and $b_s \geq 0$ such that we have the identity of rational functions*

$$\mu^{M_\alpha}(\sigma \otimes \cdot) = c'_s \frac{(1 - X_\alpha(\cdot))(1 - X_\alpha^{-1}(\cdot))}{(1 - X_\alpha(\cdot)q^{-a_s})(1 - X_\alpha^{-1}(\cdot)q^{-a_s})} \frac{(1 + X_\alpha(\cdot))(1 + X_\alpha^{-1}(\cdot))}{(1 + X_\alpha(\cdot)q^{-b_s})(1 + X_\alpha^{-1}(\cdot)q^{-b_s})}.$$

1.7

*Remark* As $\Delta_{\mathcal{O},\mu}$ is linearly independent, we can choose $\sigma$ such that, for every $\alpha \in \Delta_{\mathcal{O},\mu}$, the real numbers $a_{s_\alpha}$ and $b_{s_\alpha}$ of the statement of the proposition above satisfies $a_{s_\alpha} \geq b_{s_\alpha}$, which we will now assume. In particular, $a_{s_\alpha} > 0$.

1.8

When $P' = MU'$ is another parabolic subgroup, denote by $J_{P|P'}$ the intertwining operator defined in [Wal03, IV]. This is a rational operator

$$\mathfrak{X}^{\mathrm{nr}}(M) \to \mathrm{Hom}_{\mathbb{C}}(i_{P'\cap K}^K E, i_{P\cap K}^K E), \chi \mapsto J_{P|P'}(\sigma \otimes \chi)$$

which induces at every regular point $\chi$ a homomorphism between the representations $i_{P'}^G(\sigma \otimes \chi)$ and $i_P^G(\sigma \otimes \chi)$.

**Lemma** *Let $\alpha \in \Delta_{\mathcal{O},\mu}$, $s = s_\alpha$, and suppose that $M_\alpha$ is a standard Levi subgroup of $G$. The coordinates of the operator $J_{P|sP}$ are rational functions in $Y_\alpha$. The pols of $J_{P|sP}$ are precisely the zeros of $\mu^{M_\alpha}$. Every pole is of order 1 and its residue is bijective. Furthermore, $J_{P|sP}J_{sP|P}$ is equal to $(\mu^{M_\alpha})^{-1}$ up to multiplication by a nonzero constant.*

*Proof.* The rationality property results from [Wal03, IV.1.1], by remarking that $J_{P|sP}$ is invariant by $\mathfrak{X}^{\mathrm{nr}}(M_\alpha)$. Thanks to the functoriality property that the intertwining operators satisfy relative to parabolic induction, we return to the case where $P$ is a maximal parabolic subgroup. Therefore, $sP = \overline{P}$, and $J_{P|sP}J_{sP|P}$ is equal to $\mu^{-1}$ multiplied by a nonzero constant.

Let $\sigma_1$ be a pole of $J_{P|\overline{P}}$. Therefore, by [Sil79, 5.4.2.1], $\sigma_1$ is a unitary representation. As, according to [Wal03, V.2.3], the operator $\mu J_{P|\overline{P}}$ is regular at $\sigma_1$, we have indeed $\mu(\sigma_1) = 0$.

Reciprocally, suppose $\mu(\sigma_1) = 0$. Therefore, as $\mu^{-1} = J_{P|\overline{P}} J_{\overline{P}|P}$, at least one of the operators $J_{P|\overline{P}}$ or $J_{\overline{P}|P}$ must have a pole at $\sigma_1$. Now, such a pole is at most of order 1 (cf. [Wal03, IV.1.2]). As the poles of $\mu^{-1}$ are of order 2, $J_{P|\overline{P}}$ must have a pole of order 1.

Finally, $\mu(\sigma) = 0$ implies by the Harish-Chandra theorem **1.2** that $s\sigma \simeq \sigma$ and that $i_P^G\sigma$ and $i_{\overline{P}}^G\sigma$ are irreducibles. As the residue of $J_{P|\overline{P}}$ at $\sigma$ is an intertwining operator $i_{\overline{P}}^G\sigma \to i_P^G\sigma$, this must be bijective. □

1.9

Let $w, w' \in W$ and denote $P_{ww'}, P_{w'}$ be the standard parabolic subgroups of $G$ of Levi subgroups $ww'Mw'^{-1}w^{-1}$ and $w'Mw'^{-1}$ respectively. To simplify the notations, denote, for $w' \in W$ and $P'$ a semi-standard parabolic subgroup, by $w'^{-1}P'$ the semi-standard parabolic subgroup equal to $w'^{-1}\overline{P_{w'}}w'$. We write in the following $\mu_{P,w,w'}$ or more simply $\mu_{w,w'}$ for the rational function on $\mathfrak{X}^{\mathrm{nr}}(M)$ such that $\mu_{w,w'}(\chi) = \prod_\alpha \mu^{M_\alpha}(\sigma \otimes \chi)$, the product being on $\Sigma(P) \cap \Sigma(w'^{-1}\overline{P_{w'}}) \cap \Sigma(w'^{-1}w^{-1}P_{ww'})$. (Remark that $P_{w'} = P_{ww'} = P$, if $w, w' \in W(M, \mathcal{O})$.)

**Proposition** *Let $w, w' \in W$. Therefore, as rational function of $\chi \in \mathfrak{X}^{\mathrm{nr}}(M)$,*

$$\lambda(ww')J_{w'^{-1}w^{-1}P_{ww'}|P}(\sigma\otimes\chi) = \mu_{w,w'}(\chi)\lambda(w)J_{w^{-1}P_{ww'}|P_{w'}}(w'\sigma\otimes w'\chi)\lambda(w')J_{w'^{-1}P_{w'}|P}(\sigma\otimes\chi).$$

*Proof.* This is an immediate consequence of the composition law for the intertwining operators [Wal03, IV.3]. □

1.10

**Proposition** *Let $P' = MU'$ be another parabolic subgroup of Levi $M$. If $\Sigma(P) \cap \Sigma(\overline{P'})$ has an empty intersection with $\Sigma_{\mathcal{O},\mu}(P)$, the operator $J_{P'|P}$ is indeed defined and bijective at every point of $\mathcal{O}$.*

*Proof.* As $J_{P'|P}$ decomposes into elementary intertwining operators which come from intertwining operators relative to the $M_\alpha$, $\alpha \notin \Sigma_{\mathcal{O},\mu}(P)$, we come back to the case where $P$ is a maximal parabolic subgroup of $G$ and $P' = \overline{P}$ with $\mu$ constant. According to [Sil79, 5.4.2.1], the poles of $J_{\overline{P}|P}$ are representations which become unitary up to torsion by an unramified character of $G$. As $\mu$ is constant and that the operator $\mu J_{\overline{P}|P}$ is regular at every unitary representation of $\mathcal{O}$ [Wal03, V.2.3], $J_{\overline{P}|P}$ must itself be regular over every $\mathcal{O}$, the operator $J_{\overline{P}|P}$ being polynomial over $\mathfrak{X}^{\mathrm{nr}}(G)$. (This can be seen directly from the definition of $J_{\overline{P}|P}$.) It is bijective at every point of $\mathcal{O}$, since $\mu$ does not cancel. □

1.11

**Definition** We put $R(\mathcal{O}) = \{w \in W(M,\mathcal{O}) | w\alpha \in \Sigma(P)$ for every $\alpha \in \Sigma(P) \cap \Sigma_{\mathcal{O},\mu}\}$.

1.12

**Proposition** *The group $R(\mathcal{O})$ is a subgroup of $W(M,\mathcal{O})$. We have $W(M,\mathcal{O}) = R(\mathcal{O}) \ltimes W_{\mathcal{O}}$.*

*Proof.* Set $R = R(\mathcal{O})$. As $\Sigma_{\mathcal{O},\mu}$ is $W(M,\mathcal{O})$-invariant (cf. proof of **1.3**), $R$ is by definition a subgroup of $W(M,\mathcal{O})$. The group $W_{\mathcal{O}}$ is generated by the symmetries $s_\alpha$, $\alpha \in \Delta_{\mathcal{O}}$. As $ws_\alpha w^{-1} = s_{w\alpha}$ and that $w\alpha \in \Sigma_{\mathcal{O},\mu}$, $W_{\mathcal{O}}$ is a distinguished subgroup of $W(M,\mathcal{O})$. As $W_{\mathcal{O}}$ permutes the Weyl chambers in $\Sigma_{\mathcal{O},\mu}$, we have $W_{\mathcal{O}} \cap R = \{1\}$. It remains to show that $W(M,\mathcal{O}) = RW_{\mathcal{O}}$.

Now, let $w \in W(M,\mathcal{O})$. Therefore, $w(\Sigma(P) \cap \Sigma_{\mathcal{O},\mu}) = \Sigma(wPw^{-1}) \cap \Sigma_{\mathcal{O},\mu}$. This is the set of positive roots in $\Sigma_{\mathcal{O},\mu}$ for a certain order on $\Sigma_{\mathcal{O},\mu}$. As $W_{\mathcal{O}}$ permutes transitively the different orders on $\Sigma_{\mathcal{O}}$, there exists $w' \in W_{\mathcal{O}}$ such that $w'w(\Sigma(P) \cap \Sigma_{\mathcal{O},\mu}) = \Sigma(P) \cap \Sigma_{\mathcal{O},\mu}$, i.e. $w'w \in R$. □

1.13

**Proposition** *Suppose that $G$ is a symplectic or orthogonal group, and denote $d = \mathrm{rank}_F(G)$.*
(a) *Replacing $\sigma$ by another element of $\mathcal{O}$ and conjugating $(M,\sigma)$ by an element of $G$, we can suppose $M = \underline{M}(F)$ with*

$$\underline{M} = \mathrm{GL}_{k_1} \times \cdots \times \mathrm{GL}_{k_1} \times \mathrm{GL}_{k_2} \times \cdots \times \mathrm{GL}_{k_2} \times \cdots \times \mathrm{GL}_{k_r} \times \cdots \times \mathrm{GL}_{k_r} \times \underline{H}_k,$$

*where $\underline{H}_k$ denotes a semisimple group of absolute rank $k$ of the same type as $G$, and*

$$\sigma = \sigma_1 \otimes \cdots \otimes \sigma_1 \otimes \sigma_2 \otimes \cdots \otimes \sigma_2 \otimes \cdots \otimes \sigma_r \otimes \cdots \otimes \sigma_r \otimes \tau,$$

*the inertial classes of $\sigma_i$ being pairwise distinct, as well as $\sigma_i \simeq \sigma_i^\vee$ if $\sigma_i$ and $\sigma_i^\vee$ are in the same inertial orbit.*
(b) *Denote by $d_i$ the number of factors equal to $\sigma_i$ and identify $A_{\underline{M}}$ with $\mathbb{T} = \mathbb{G}_m^{d_1} \times \mathbb{G}_m^{d_2} \times \cdots \times \mathbb{G}_m^{d_r}$. Denote by $\alpha_{i,j}$ the rational characters of $A_{\underline{M}}$ (identified with $\mathbb{T}$) which sends an element $x = (x_{1,1}, \cdots, x_{1,d_1}, x_{2,1}, \cdots, x_{2,d_2}, \cdots, x_{r,1}, \cdots, x_{r,d_r})$ to $x_{i,j}x_{i,j+1}^{-1}$, if $j < d_i$, and to $x_{i,d_i}$, if $j = d_i$. The root system $\Sigma_{\mathcal{O},\mu}$ is the direct sum of $r$ components irreducible or empty $\Sigma_{\mathcal{O},\mu,i}$, $i = 1, \cdots, r$, defined in the following way:*
*Suppose firstly either $k \neq 0$ or $G$ has root system of type $B_d$.*

(i)*If the function $s \mapsto \mu(\sigma_i|\det_{k_i}|^s \otimes \tau)$ (defined relative to $\mathrm{GL}_{k_i} \times \underline{H}_k$ and $\underline{H}_{k_i+k}$) has a pole over $\mathbb{C}$, thus a basis of $\Sigma_{\mathcal{O},\mu,i}$ is given by $\{\alpha_{i,1}, \cdots, \alpha_{i,d_i}\}$, and this system is of type $B_{d_i}$.*

(ii)*If the function* $s \mapsto \mu(\sigma_i|\det_{k_i}|^s \otimes \tau)$ *(defined relative to* $\mathrm{GL}_{k_i} \times \underline{H}_k$ *and* $\underline{H}_{k_i+k}$*) is regular on* $\mathbb{C}$ *and* $\sigma_i \simeq \sigma_i^\vee$*, thus a basis of* $\Sigma_{\mathcal{O},\mu,i}$ *is given by* $\{\alpha_{i,1}, \cdots, \alpha_{i,d_i-1}, \alpha_{i,d_i-1} + 2\alpha_{i,d_i}\}$*, and this system is of type* $D_{d_i}$.

(iii)*If not, a basis of* $\Sigma_{\mathcal{O},\mu,i}$ *is given by* $\{\alpha_{i,1}, \cdots, \alpha_{i,d_i-1}\}$*, and this system is of type* $A_{d_i-1}$.

*Suppose for now* $k = 0$ *and either* $G$ *has root system of type* $C_d$ *or* $G$ *has root system of type* $D_d$ *and* $k_i \geq 2$*. Therefore, in the case (i) above, a basis of* $\Sigma_{\mathcal{O},\mu,i}$ *is given by* $\{\alpha_{i,1}, \cdots, 2\alpha_{i,d_i}\}$*, and the system is of type* $C_{d_i}$*. In the other case, the situation remains unchanged.*

*Suppose finally* $k = 0$*,* $G$ *has root system of type* $D_d$ *and* $k_i = 1$*. Therefore, if* $\sigma_i \simeq \sigma_i^\vee$*, a basis of* $\Sigma_{\mathcal{O},\mu,i}$ *is given by* $\{\alpha_{i,1}, \cdots, \alpha_{i,d_i-1}, \alpha_{i,d_i-1} + 2\alpha_{i,d_i}\}$ *and this system is of type* $D_{d_i}$*. If not, the system is of type* $A_{d_i-1}$ *and a basis is given by* $\{\alpha_{i,1}, \cdots, \alpha_{i,d_i-1}\}$.

*Proof.* The Levi subgroup of classical groups are described in [Gau66, Kap. 5]. By permutation of factors, we put $\sigma$ in the desired form. Suppose that there exists an unramified character $\chi_i$ of $\mathrm{GL}_{k_i}$ such that $\sigma_i^\vee \simeq \sigma_i \otimes \chi_i$. We can write $\chi_i = |\det_{k_i}|_F^{s_i}$, where $s_i$ is a complex number. Put $\chi_i^{1/2} = |\det_{k_i}|^{s_i/2}$. Therefore $(\sigma_i \otimes \chi_i^{1/2})^\vee \simeq \sigma_i^\vee \otimes \chi_i^{-1/2}\chi_i = \sigma_i \otimes \chi_i^{1/2}$.

Identify $A_{\underline{M}}$ with torus $\mathbb{T}$, write $x = (x_{1,1}, \cdots, x_{1,d_1}, x_{2,1}, \cdots, x_{2,d_2}, \cdots, x_{r,1}, \cdots, x_{r,d_r})$ for the general element of $\mathbb{T}$ and let's examine the different scenarios.

Suppose firstly that either $k \neq 0$ or $G$ has a root system of type $B_d$. The reduced roots in $\Sigma(A_M)$ are thus $x \mapsto x_{i,j}^{\pm 1} x_{i',j'}^{\pm 1}$, $(i,j) \neq (i',j')$, as well as $x \mapsto x_{i,j}^{\pm 1}$. (As $k \neq 0$, there indeed exists, in the case of a root system of type $C_d$ or $D_d$, roots of restriction to $A_M$ equal to $x \mapsto x_{i,j}^{\pm 1}$.) If $\alpha \in \Sigma(A_M)$ corresponds to a root $x \mapsto x_{i,j} x_{i',j'}^{-1}$, $(i,j) \neq (i',j')$, thus the function $\mu^{M_\alpha}$ is equal to that defined from the representation $\sigma_i \otimes \sigma_{i'}$ of Levi subgroup $\mathrm{GL}_{k_i} \times \mathrm{GL}_{k_{i'}}$ of $\mathrm{GL}_{k_i+k_{i'}}$. If $\alpha \in \Sigma(A_M)$ correspond to a root $x \mapsto x_{i,j} x_{i',j'}$, $(i,j) \neq (i',j')$, $i \leq i'$, thus the function $\mu^{M_\alpha}$ is equal to that defined from the representation $\sigma_i \otimes \sigma_{i'}^\vee$ of Levi subgroup $\mathrm{GL}_{k_i} \times \mathrm{GL}_{k_{i'}}$ of $\mathrm{GL}_{k_i+k_{i'}}$. In the other case, the function $\mu^{M_\alpha}$ is that defined by the representation $\sigma_i \otimes \tau$ of Levi subgroup $\mathrm{GL}_{k_i} \times H_m$ of $H_{m+k_i}$. According to the different cases presented in (i), (ii) and (iii), we thus deduce from Harish-Chandra theorem **1.2** and from results of Bernstein-Zelevinsky [BZ77, 4.2] the indicated result.

For $k = 0$ and $G$ with root system of type $C_d$, the reduced roots in $\Sigma(A_M)$ are $x \mapsto x_{i,j}^{\pm 1} x_{i',j'}^{\pm 1}$, $(i,j) \neq (i',j')$, as well as $x \mapsto x_{i,j}^{\pm 2}$. We thus conclude as above.

For $k = 0$ and $G$ with root system finally of type $D_d$, the reduced root in $\Sigma(A_M)$ are $x \mapsto x_{i,j}^{\pm 1} x_{i',j'}^{\pm 1}$, $(i,j) \neq (i',j')$, as well as, if $k_i \geq 2$, $x \mapsto x_{i,j}^{\pm 2}$. We thus conclude as above, taking into account that the function $\mu$ is indeed known for the split groups of rank 1. □

1.14

**Proposition** *Suppose that $G$ is the group of rational points of an inner form of $\mathrm{GL}_n$, i.e. $G = \mathrm{GL}_n(D)$, where $D$ is a division algebra with center $F$.*

*Replacing $\sigma$ by another element of $\mathcal{O}$ and conjugating $(M, \sigma)$ by an element of $G$, we can suppose $M$ equal to*

$$\mathrm{GL}_{k_1}(D)\times\cdots\times\mathrm{GL}_{k_1}(D)\times\mathrm{GL}_{k_2}(D)\times\cdots\times\mathrm{GL}_{k_2}(D)\times\cdots\times\mathrm{GL}_{k_r}(D)\times\cdots\times\mathrm{GL}_{k_r}(D)$$

*and*

$$\sigma = \sigma_1 \otimes \cdots \otimes \sigma_1 \otimes \sigma_2 \otimes \cdots \otimes \sigma_2 \otimes \cdots \otimes \sigma_r \otimes \cdots \otimes \sigma_r,$$

*the inertial classes of $\sigma_i$ being pairwise distinct.*

*Denote by $d_i$ the number of factors equal to $\sigma_i$ and identify $A_{\underline{M}}$ with $\mathbb{T} = \mathbb{G}_m^{d_1} \times \mathbb{G}_m^{d_2} \times \cdots \times \mathbb{G}_m^{d_r}$. Denote by $\alpha_{i,j}$ the rational characters of $A_{\underline{M}}$ (identified with $\mathbb{T}$) which sends an element $x = (x_{1,1}, \cdots, x_{1,d_1}, x_{2,1}, \cdots, x_{2,d_2}, \cdots, x_{r,1}, \cdots, x_{r,d_r})$ to $x_{i,j}x_{i,j+1}^{-1}$. The root system $\Sigma_{\mathcal{O},\mu}$ is the direct sum of $r$ components irreducible or empty $\Sigma_{\mathcal{O},\mu,i}$, $i = 1, \cdots, r$, of type $A_{d_i-1}$ and of basis $\{\alpha_{i,1}, \cdots, \alpha_{i,d_i-1}\}$ respectively.*

*Proof.* For split $G$, this results from the work of Bernstein-Zelevinsky [BZ77, 4.2] with the Harish-Chandra theorem **1.2**. In the general case, we deduce it from the work of Bernstein-Zelevinsky with the trace formula [DKV84, Tad90]. □

1.15

**Proposition** *(i) If $G$ is a linear group or the multiplicative group of a division algebra, then $R(\mathcal{O}) = 1$.*
*(ii) If not, $R(\mathcal{O}) \neq 1$, if and only if the conclusions of (ii) of proposition **1.13** are verified for at least an index $i$ with, if $G$ is of type $D_n$, moreover $k_i$ even or else $H \neq 1$ and $\tau$ invariant by the outer automorphism. The group $R(\mathcal{O})$ is thus the direct product of groups $R(\mathcal{O})_i$ indexed by the $i$ for which the conclusions above are verified with $R(\mathcal{O})_i = \{1, s_{\alpha_{i,d_i}}\}$ ($s_{\alpha_{i,d_i}}$ exchanging $\alpha_{i,d_{i-1}}$ and $\alpha_{i,d_{i-1}} + 2\alpha_{i,d_i}$ and satisfying $s_{\alpha_{i,d_i}}\sigma_i \simeq \sigma_i^\vee$).*

*Proof.* It is easy to see that the conditions of theorem **1.2** are assembled exactly in the case described in the statement, taking into account [BJ03, 3.4] for the groups of type $D_n$. □

The proposition **1.15** as well as the proposition **1.12** lead us to add the following two hypotheses on the choice of the representatives of the elements of $W(M, \mathcal{O})$: concerning the representatives of the elements in $R(\mathcal{O})$, we suppose that, if $r$ is a product of elements $s_{\alpha_i}$ in $R(\mathcal{O})_i$ pairwise distinct in $R(\mathcal{O})$, then the same is true for the representatives. Concerning an arbitrary element $w$ of $W(M,\mathcal{O})$, we suppose that, if $w = w_{\mathcal{O}} r$ with $w_{\mathcal{O}} \in W_{\mathcal{O}}$ and $r \in R(\mathcal{O})$, then the representative of $w$ is the product of those of $w_{\mathcal{O}}$ with that of $r$.

1.16

The following proposition is a particular case of a well-known result for the finite groups [Hen88, Lemma 7.5]. Note that it is probably not valid for an arbitrary reductive group over $F$, although the proof remains valid in cases other than that of a symplectic or orthogonal group or an inner form of a general linear group considered here.

**Proposition** *Let $(\sigma_1, E_1)$ be an irreducible component of $\sigma|_{M^1}$. Thus $\sigma_1$ extends to a representation $\sigma_2$ of $M^\sigma = \{m \in M | {}^m\sigma_1 \simeq \sigma_1\}$ such that $\sigma|_{M^\sigma} = \bigoplus_{m \in M/M^\sigma} {}^m\sigma_2$. In particular, $\sigma = \mathrm{Ind}_{M^\sigma}^M \sigma_2$ and ${}^m\sigma_2 \not\simeq \sigma_2$ if $m \notin M^\sigma$. Moreover, $\sigma|_{M^1}$ is direct sum of pairwise non-isomorphic irreducible representations.*

*Proof.* We first assume that $A_M$ is of rank 1. First note that, if $H$ is an open closed subgroup of finite index of a locally profinite group $G$, then $\mathrm{ind}_H^G = \mathrm{Ind}_H^G$ and, as in the case of finite groups, a representation $(\sigma, E)$ is induced by a representation $(\sigma_2, E_1)$, if and only if $V = \bigoplus_{g \in G/H} \sigma(g) E_1$.

Observe for now that $A_M \subseteq M^\sigma$. The quotient $M/M^\sigma$ is thus a finite commutative group. As $\sigma$ is cuspidal, we can write $E = E_{\sigma_1} \oplus E'$ such that $E_{\sigma_1}$ is direct sum of representations isomorphic to $\sigma_1$ and that no subquotient of $E'$ is not isomorphic to $\sigma_1$ [BZ76, 2.44]. The group $M^\sigma$ thus acts on $E_{\sigma_1}$, and we have $E = \bigoplus_{m \in M/M^\sigma} \sigma(m) E_{\sigma_1}$. Hence, by the remark above, $\mathrm{Ind}_{M^\sigma}^M E_{\sigma_1} = E$. AS $E$ is irreducible, the representation of $M^\sigma$ in $E_{\sigma_1}$ must be irreducible. We will now show that in fact $E_{\sigma_1} = E_1$.

First let's see that $\sigma_1$ extends to an irreducible representation $\sigma_2$ of $M^\sigma$. First, $\sigma_1$ extends to an irreducible representation $\widetilde{\sigma}_1$ of $A_M M$, given by $am \mapsto \chi_\sigma(a)\sigma_1(m)$, $\chi_\sigma$ denotes the central character of $\sigma$. As the quotient $M^\sigma / A_M M^1$ is by hypothesis a finite cyclic group, it is thus easy to define an extension to $M^\sigma$. (Given an element $m'$ of $M^\sigma$ whose image generates the cyclic quotient group $M^\sigma/A_M M^1$, we can always find an $A_M M^1$-equivariant isomorphism $\varphi_{m'}$ : $E_1 \to m'^{-1} E_1$ such that $m'^i m \mapsto \varphi_{m'}^i \widetilde{\sigma}_1(m)$ for $m \in A_M M^1$ gives the desired extension.)

Let us now show that $\mathrm{Ind}_{A_M M^1}^{M^\sigma} \widetilde{\sigma}_1 = \bigoplus_{\chi \in X(M^\sigma/A_M M^1)} \sigma_2 \otimes \chi$. The representations $\sigma_2 \otimes \chi$, $\chi \in X(M^\sigma/A_M M^1)$ are pairwise non-isomorphic, since every isomorphism $\sigma_2 \otimes \chi \mapsto \sigma_2 \otimes \chi'$ induces an automorphism of $\widetilde{\sigma}_1$ and, by the Schur lemma, it would therefore be scalar, which is impossible. By Frobenius reciprocity, the irreducible sub-representations of $\mathrm{Ind}_{A_M M^1}^{\sigma}\widetilde{\sigma}_1$ are irreducible representations of $M^\sigma$ with restriction to $A_M M^1$ equal to $\widetilde{\sigma}_1$ and their multiplicity is 1. Now, the arguments of the proof of [BZ76, 3.29] show that every irreducible representation of $M^\sigma$ with restriction to $A_M M^1$ equal to $\widetilde{\sigma}_1$ is isomorphic to $\sigma_2\otimes\chi$ for a $\chi\in X(M^\sigma/A_M M^1)$. As $\mathrm{Ind}_{A_M M^1}^{M^\sigma}\widetilde{\sigma}_1$ is unitary, since $\widetilde{\sigma}_1$ is so, we have indeed the indicated decomposition.

It follows that $E_{\sigma_1}=E_1$: the representation $\rho$ of $M^\sigma$ in $E_{\sigma_1}$ is irreducible, with restriction equal to a multiple of $\sigma_1$. The Frobenius reciprocity thus gives $\mathrm{Hom}_M(\rho,\mathrm{Ind}_{A_M M^1}^{M_\sigma}\widetilde{\sigma}_1)\neq 0$, this implies that $\rho$ is isomorphic to $\sigma_2\otimes\chi$ for a $\chi\in X(M^\sigma/A_M M^1)$. Consequently, $E_{\sigma_1}=E_1$ and we can suppose in the following that $\rho=\sigma_2$.

As a result, noting that $M^\sigma$ is distinguished in $M$, that $\sigma=\mathrm{Ind}_{M^\sigma}^{M}\sigma_2$ and $\sigma|_{M^\sigma}=\bigoplus_{m\in M/M^\sigma}{}^m\sigma_2$. As $\sigma$ is irreducible, the ${}^m\sigma_2$ are pairwise non-isomorphic. Since moreover $({}^m\sigma_2)|_{M^1}={}^m\sigma_1$ with $m$ running through a system of representatives of $M/M^\sigma$, $\sigma|_{M^1}$ is direct sum of pairwise non-isomorphic irreducible representations.

Now let's remove the hypothesis that $A_M$ is of rank 1. As $G$ is a symplectic or orthogonal group, or the multiplicative of a simple algebra, $M$ is of the form $\mathrm{GL}_{m_1}(D)\times\cdots\times\mathrm{GL}_{m_r}(D)\times H$, where $D$ is a division algebra with center $F$ and where $H$ is a semisimple group of the same type as $G$ or the trivial group. We deduce that $M^1$ is of the form $\mathrm{GL}_{m_1}(D)^1\times\cdots\times\mathrm{GL}_{m_r}(D)^1\times H$, that $\sigma$ is isomorphic to a product $\sigma_1\otimes\cdots\otimes\sigma_r\otimes\tau$ with $\sigma_i,\tau$ cuspidal representations, and that $M^\sigma$ is of the form $\mathrm{GL}_{m_1}(D)^{\sigma_1}\times\cdots\times\mathrm{GL}_{m_r}(D)^{\sigma_r}\times H$. We are thus reduced to the case where $M=\mathrm{GL}_m(D)$. But, here $A_M$ is of rank 1, and the assertion is true by what preceded. □

1.17

When $\chi$ is an unramified character of $M$, denote by $E_\chi$ the space $E$ equipped with the representation $\sigma\otimes\chi$.

**Corollary** (with the notations of proposition 1.16) *We have* $\mathrm{Stab}(\mathcal{O})=X(M/M^\sigma)$ *and, for every* $\chi\in X(M/M^\sigma)$*, the map* $\phi_{\sigma,\chi}:E\to E_\chi$ *which sends an element* $e\in\sigma(m)E_1$ *to* $\chi(m)e$ *is an isomorphism.*

*Proof.* Let $\chi\in\mathfrak{X}^{\mathrm{nr}}(M)$ such that $\sigma\otimes\chi\simeq\sigma$. Therefore, there exists $m\in M$ such that $\sigma_2\otimes\chi\simeq{}^m\sigma_2$. Consequently, $\sigma_1\simeq(\sigma_2\otimes\chi)|_{M^1}\simeq({}^m\sigma_2)|_{M^1}={}^m\sigma_1$, i.e. we have $m\in M^\sigma$ and thus $\sigma_2\otimes\chi\simeq\sigma_2$. Let $\varphi_\chi$ be an automorphism of $E_1$ which intertwines $\sigma_2\otimes\chi$ and $\sigma_2$. Therefore it intertwines $\sigma_1$ with itself, i.e. $\varphi_\chi$ is the multiplication by a scalar. It then follows that $\chi|_{M^\sigma}\equiv 1$. Thus $\mathrm{Stab}(\mathcal{O})\subseteq X(M/M^\sigma)$. Conversely, it is an immediate verification that for every $\chi \in X(M/M^\sigma)$ the map $\phi_{\sigma,\chi}$ defines isomorphisms between $\sigma$ and $\sigma \otimes \chi$. $\square$

## 2

Denote $B = B_M$ the ring of polynomials on the complex affine variety $\mathfrak{X}^{\mathrm{nr}}(M)$. This is an integral ring, and we denote by $K(B)$ its fraction field.

Write $E_B$ (resp. $E_{K(B)}$) for the complex vector space $E \otimes_{\mathbb{C}} B$ (resp. $E_B \otimes_B K(B) \simeq E \otimes_{\mathbb{C}} K(B)$). Let's equip it with the action of $M$ given by $\sigma_B : M \to \mathrm{Aut}_{\mathbb{C}}(E_{K(B)})$, $\sigma_B(m)(e \otimes b) := \sigma(m)e \otimes bb_m$. Here, $b_m : \mathfrak{X}^{\mathrm{nr}}(M) \to \mathbb{C}$ is defined by $b_m(\chi) = \chi(m)$.

By restriction to $K$, the space $i_P^G E_B$ (resp. $i_P^G E_{K(B)}$) is isomorphic to $i_{K\cap P}^K E \otimes_{\mathbb{C}} B$ (resp. $i_{K\cap P}^K E \otimes_{\mathbb{C}} K(B)$). We also identify $i_{K\cap P}^K E$ to its image in $i_P^G E_B$ (resp. $i_P^G E_{K(B)}$) by the canonical homomorphism $v \mapsto v \otimes 1$.

### 2.1

Note first that $E_B$ is canonically isomorphic to $\mathrm{ind}_{M^1}^M E$. Indeed, denote by $\mathcal{R}(M/M^1)$ a system of representatives of $M/M^1$. We obtain a canonical isomorphism between $B$ and $\mathrm{ind}_{M^1}^M \mathbb{C} = \mathbb{C}[M/M^1]$, by sending an element of the form $b_m, m \in \mathcal{R}(M/M^1)$, to the map $M \to \mathbb{C}$ which is 1 on $m^{-1}M^1$ and 0 outside of this set. We identify these two spaces using this isomorphism. We have a canonical isomorphism $E \otimes \mathrm{ind}_{M^1}^M \mathbb{C} \to \mathrm{ind}_{M^1}^M E|_{M^1}$, by sending $e \otimes b$ to the map $v_{e\otimes b}$ in $\mathrm{ind}_{M^1}^M E|_{M^1}$ which sends $m$ to $b(m)\sigma(m)e$. The inverse isomorphism is given by $\mathrm{ind}_{M^1}^M E|_{M^1} \to E \otimes \mathrm{ind}_{M^1}^M \mathbb{C}$, $v \mapsto \sum_{m \in \mathcal{R}(M/M^1)} \sigma(m^{-1})v(m) \otimes b_{m^{-1}}$. (Here, $b_{m^{-1}}$ is considered as element of $\mathrm{ind}_{M^1}^M \mathbb{C}$ by means of the isomorphism above.)

### 2.2

For $w \in W^G$, denote by $wE_B$ (resp. $wE_{K(B)}$) the space $E_B$ (resp. $E_{K(B)}$) equipped with the representation $w\sigma_B$ of $wMw^{-1}$. (We have thus $w\sigma_B(wmw^{-1}) = \sigma_B(m)$.) Denote by $\tau_w$ the $\mathbb{C}$-vector-space-automorphism of $E \otimes_{\mathbb{C}} K(B)$ which sends $e \otimes b$ to $e \otimes {}^w b$, ${}^w b(\chi) := b(w^{-1}\chi)$. (We thus have ${}^w b_m = b_{wmw^{-1}}$ for $m \in M$.) It defines an $M$-equivariant isomorphism between $wE_{K(B)}$ and $(wE)_{K(B)}$ and we will again denote by $\tau_w$ the $G$-equivariant isomorphism $i_P^G w(E_{K(B)}) \to i_P^G (wE)_{K(B)}$ which is deduced by functoriality.

### 2.3

For $b \in K(B)$, denote by $b_\chi$ the element of $K(B)$ given by $b_\chi(\chi') = b(\chi\chi')$. We have $(b_m)_\chi = \chi(m)b_m$. The following result is a direct verification:

**Lemma** *Let $\chi \in \mathfrak{X}^{\mathrm{nr}}(M)$ and denote by $E_\chi$ the space $E$ equipped with the representation $\sigma \otimes \chi$ of $M$. The map $\rho_\chi : E \otimes_{\mathbb{C}} B \to E_\chi \otimes_{\mathbb{C}} B$, $e \otimes b \mapsto e \otimes b_\chi$, defines an isomorphism between the smooth representations of $M$. We have $\rho_\chi^{-1} = \rho_{\chi^{-1}}$ and $\tau_w \circ \rho_\chi = \rho_{w\chi}\tau_w$.*

2.4

We will define for every $w \in W_{\mathcal{O}}$ an isomorphism $\rho_w : i_P^G wE \to i_P^G E$ which sends $i_P^G(wE)_B$ to $i_P^G E_B$. For this, let's set

$$\rho_w = \left[\left(\left(\prod_{\alpha \in \Sigma_{\mathcal{O},\mu}(P) \cap w^{-1}\Sigma_{\mathcal{O},\mu}(\overline{P})} (Y_\alpha(\chi) - 1)\right)\lambda(w)J_{w^{-1}P|P}(\sigma \otimes \chi)\right)_{|\chi=1}\right]^{-1},$$

if this expression is well-defined. We write $\rho_{\sigma,w}$, if we want to emphasize the dependence on $\sigma$. It is moreover clear that the operators $\rho_{\sigma,w}$ and $\tau_w$ commute.

**Proposition** *The map $\rho_w$ comes from an isomorphism $wE \simeq E$, and it defines for every $w \in W_{\mathcal{O}}$ an isomorphism $i_P^G wE \to i_P^G E$. It satisfies the following compatibility properties: for every $w' \in W$ with $w'Mw'^{-1}$ equal to standard Levi subgroup of a standard parabolic subgroup $P_{w'}$ of $G$, $\rho_w$ induces an isomorphism $i_{P_{w'}}^G(w'wE)_{K(B)} \to i_{P_{w'}}^G(w'E)_{K(B)}$, and we have the commutation formula*

$$\lambda(w')J_{w'^{-1}P_{w'}|P}(\sigma \otimes \cdot)\rho_w = \rho_w \lambda(w')J_{w'^{-1}P_{w'}|P}(w\sigma \otimes \cdot).$$

*Moreover, for $w' \in W_{\mathcal{O}}$, we have, with $w''$ equal to representative of $ww'$ in $K$,*

$$\rho_w\rho_{w'} = \left(\mu_{w,w'}\prod_\alpha (Y_\alpha - 1)^{-1}(Y_\alpha^{-1} - 1)^{-1}\right)_{|\chi=1}\lambda(w''(ww')^{-1})\rho_{ww'},$$

*the product taken over the set $\Sigma_{\mathcal{O},\mu}(P) \cap w'^{-1}\Sigma_{\mathcal{O},\mu}(\overline{P}) \cap w'^{-1}w^{-1}\Sigma_{\mathcal{O},\mu}(P)$, as well as that, for $r \in R(\mathcal{O})$ and $s_\alpha$ a simple rflection in $W_{\mathcal{O}}$,*

$$\rho_{r\sigma,s_{r\alpha}} = \lambda(rs_\alpha r^{-1}s_{r\alpha}^{-1})\rho_{\sigma,s_\alpha}.$$

*Proof.* Suppose first $w = s_\alpha$ with $\alpha \in \Delta_{\mathcal{O},\mu}$. Therefore, $\mu^{M_\alpha}(\sigma) = 0$, and it results from the Harish-Chandra theorem **1.2** that $i_{P\cap M_\alpha}^{M_\alpha}\sigma$ is irreducible and that $s_\alpha\sigma \simeq \sigma$.

If $M_\alpha$ is standard, we deduce from **1.8** that the operator $J_{s_\alpha P|P}$ considered as rational function in $Y_\alpha$ has a pole of order 1 at $Y_\alpha = 1$ and that the operator $\mathrm{Res}_{Y_\alpha=1}(\lambda(s_\alpha)J_{s_\alpha P|P}(\sigma \otimes (\cdot))) : i_P^G\sigma \to i_P^G s_\alpha\sigma$ is bijective. Consequently, $\rho_{s_\alpha}$ is well-defined and bijective. It is induced by functoriality by an isomorphism between the irreducible representations $i_{P\cap M_\alpha}^{M_\alpha} s_\alpha\sigma$ and $i_{P\cap M_\alpha}^{M_\alpha}\sigma$. Now, such an isomorphism is, by Schur lemma, uniquely determined up to a constant. In consequence, $\rho_{s_\alpha}$ is itself induced by functoriality by an isomorphism $s_\alpha\sigma \to \sigma$. The two compatibility properties result from the functoriality property for the intertwining operators.

If $M_\alpha$ is not standard, so $\alpha$ is contained in an irreducible component of type $D_{d_i}$ of $\Sigma_{\mathcal{O},\mu}$, and there exists by **1.15** and by **1.13** an element $r \in R(\mathcal{O})$ such that $\alpha = r\beta$ with $\beta \in \Delta_{\mathcal{O},\mu}$, $M_\beta$ standard and $r\sigma \simeq \sigma$, or this irreducible component is of type $B_{d_i}$ or $C_{d_i}$ and $\alpha$ is the short (resp. long) root.

In the first case, following **1.9**, we decompose

$$(\#)\quad \begin{aligned}&\lambda(s_\alpha)J_{s_\alpha P|P}(\sigma\otimes\chi)\\ &\quad=\mu_{r,s_\beta r^{-1}}(\chi)\mu_{s_\beta,r^{-1}}(\chi)\lambda(s_\alpha r s_\beta^{-1}r^{-1})\lambda(r)J_{r^{-1}P|P}(s_\beta r^{-1}\sigma\otimes s_\beta r^{-1}\chi)\\ &\quad\quad\times\lambda(s_\beta)J_{s_\beta P|P}(r^{-1}\sigma\otimes r^{-1}\chi)\lambda(r^{-1})J_{rP|P}(\sigma\otimes\chi).\end{aligned}$$

We have $\Sigma(P)\cap\Sigma(r\overline{P})\cap\Sigma(rs_\beta P)=\varnothing$, since $r\beta\in\Sigma(P)$ and $\Sigma(P)\cap\Sigma(s_\beta\overline{P})=\{\beta\}$. The function $\mu_{s_\beta,r}$ is thus constant. Likewise, the set $\Sigma(P)\cap\Sigma(rs_\beta\overline{P})\cap\Sigma(s_\alpha P)$ has an empty intersection with $\Sigma_{\mathcal{O},\mu}$, which shows that $\mu_{r,s_\beta r^{-1}}$ is constant. As $r\sigma\simeq\sigma$ and that $J_{rP|P}(\sigma\otimes\chi)$ is according to **1.10** well-defined and bijective for every $\chi$, we see, applying the above to $\rho_{r^{-1}\sigma,s_\beta}$ and using the equality ${}^rX_\beta = X_\alpha$ (cf. **1.5**), that the operator $\rho_{s_\alpha}$ is well-defined and bijective. As $\rho_{s_\beta}$ is induced by funcoriality from an isomorphism $s_\beta E\to E$, we deduce the same for $\rho_{s_\alpha}$, after applying the compatibility properties of $\rho_{s_\beta}$ to the decomposition formula for $\rho_{s_\alpha}^{-1}$ above. The compatibility property for $\rho_{s_\alpha}$ follows.

If $M_\alpha$ is not standard and $\alpha$ belongs to an irreducible component of type $B_{d_i}$ or $C_{d_i}$, then denote $w$ the element of minimal length in $W$ which sends an element $x$ in $A_M$ to the element of $A_0$ deduced from $x$, by exchanging the components $(x_{i,1},\cdots,x_{i,d_i})$ and $(x_{r,1},\cdots,x_{r,d_r})$. The group $wMw^{-1}$ is thus a standard Levi subgroup of $G$ and the root $\beta = w\alpha$ of $A_{wMw^{-1}}$ is in $\Delta_{w\mathcal{O},\mu}$ with $(wMw^{-1})_\beta$ standard. Denote by $P_w$ the standard parabolic subgroup of $G$ of Levi subgroup $wMw^{-1}$. The roots in $\Sigma_{\mathrm{red}}(P)\cap\Sigma_{\mathrm{red}}(w^{-1}\overline{P_w})$ are of the form $x\mapsto x_{i,j_1}x_{i_2,j_2}^{-1}$ with $i<i_2\le r$ and $x\mapsto x_{i_1,j_1}x_{r,j_2}^{-1}$, with $i\le i_1<r$. The intersection with $\Sigma_{\mathcal{O},\mu}$ is thus empty. It is the same for that of $w(\Sigma_{\mathrm{red}}(P)\cap\Sigma_{\mathrm{red}}(w^{-1}\overline{P_w}))$ with $\Sigma_{s_\beta w\mathcal{O},\mu}$. In consequence, the operators $\lambda(w)J_{w^{-1}P_w|P}(\sigma)$ and $\lambda(w^{-1})J_{wP|P_w}(s_\beta w\sigma)$ are according to **1.10** well-defined and bijective. The functions $\mu_{w^{-1},s_\beta w}$ and $\mu_{s_\beta,w}$ are constant. By the decomposition formula **1.9**, we thus find, as a rational function in $\chi$,

$$\begin{aligned}
&\lambda(s_\alpha)J_{s_\alpha P|P}(\sigma\otimes\chi)\\
&\quad=\mu_{w^{-1},s_\beta w}\mu_{s_\beta,w}\lambda(s_\alpha w^{-1}s_\beta^{-1}w)\\
&\qquad\times\lambda(w^{-1})J_{wP|P_w}(s_\beta w(\sigma\otimes\chi))\lambda(s_\beta)J_{s_\beta P_w|P_w}(w(\sigma\otimes\chi))\lambda(w)J_{w^{-1}P_w|P}(\sigma\otimes\chi).
\end{aligned}$$

As $Y_\beta({}^w\chi)={}^{w^{-1}}Y_\beta(\chi)=Y_\alpha(\chi)$, it follows that

$$\rho_{\sigma,s_\alpha}=\mu_{w^{-1},s_\beta w}\mu_{s_\beta,w}\lambda(w^{-1})J_{wP|P_w}(s_\beta w\sigma)\lambda(s_\beta)\rho_{w\sigma,s_\beta}\lambda(w)J_{w^{-1}P_w|P}(\sigma),$$

i.e. following what has been previously established for $\rho_{w\sigma,s_\beta}$, $\rho_{\sigma,s_\alpha}^{-1}$ is well-defined and bijective. Indeed, we deduce from compatibility properties for $\rho_{w\sigma,s_\beta}$ that $\rho_{\sigma,s_\alpha}$ is a product of $\rho_{w\sigma,s_\beta}$ by a nonzero scalar. As $s_\alpha E=w^{-1}s_\beta(wE)$, the compatibility property for $\rho_{\sigma,s_\alpha}$ are immediate.

The last assertion of the proposition results directly from the decomposition (#) for $r\in R(\mathcal{O})$, after having multiplied it with $X_\alpha-1$, evaluating at 1 and after using the equality ${}^rX_\beta=X_\alpha$.

Now suppose we have shown the equality

$$\rho_{w'}^{-1}\rho_w^{-1}=\Big(\mu_{w,w'}^{-1}\prod_\alpha(Y_\alpha-1)(Y_\alpha^{-1}-1)\Big)_{|\chi=1}\lambda(ww'w''^{-1})\rho_{ww'}^{-1}. \qquad (*)$$

with $w''=ww'$, when $w,w'$ are elements in $W_\mathcal{O}$ such that $\rho_w$ and $\rho_{w'}$ are well-defined, bijective and such that $\rho_{w'}$ satisfies the two compatibility properties, the product taken over the set $\Sigma_{\mathcal{O},\mu}(P)\cap w'^{-1}\Sigma_{\mathcal{O},\mu}(\overline{P})\cap w'^{-1}w^{-1}\Sigma_{\mathcal{O},\mu}(P)$.

First let us deduce that $\rho_w$ is well-defined and bijective, when $w\in W_\mathcal{O}$, and that the compatibility properties are indeed satisfied. For this, let's do a recurrence on the length of $w$. The case $\ell_\mathcal{O}(w)=1$ has already been established. Let $w''=ws_\alpha\in W_\mathcal{O}$ with $\alpha\in\Delta_{\mathcal{O},\mu}$, $\ell_\mathcal{O}(w'')=\ell_\mathcal{O}(w)-1$. By recurrence hypothesis, the operators $\rho_w$ and $\rho_{s_\alpha}$ are well-defined, bijective, and satisfy the compatibility properties. We can thus apply (*), and we find

$$\rho_{w''}^{-1}=\lambda(w''s_\alpha^{-1}w^{-1})\rho_{s_\alpha}^{-1}\rho_w^{-1},$$

the set $\Sigma_{\mathcal{O},\mu}(P)\cap s_\alpha\Sigma_{\mathcal{O},\mu}(\overline{P})\cap w^{-1}s_\alpha^{-1}\Sigma_{\mathcal{O},\mu}(P)$ being empty thanks to the hypothesis $\ell(ws_\alpha)=\ell(w)+\ell(s_\alpha)$. As $w''s_\alpha^{-1}w^{-1}\in M$, it follows that $\rho_{w''}^{-1}$ is bijective, that $\rho_{w''}$ is well-defined and that it satisfies the compatibility properties.

It remains to show the equality (*). Let $w, w' \in W_{\mathcal{O}}$ satisfy the indicated hypotheses. We have

$$
\begin{aligned}
&\rho_{w'}^{-1}\rho_w^{-1}\\
&= \rho_{w'}^{-1}\left(\left(\prod_{\alpha}(Y_\alpha - 1)\right)\lambda(w)J_{w^{-1}P|P}(\sigma\otimes\chi)\right)_{|\chi=1}\\
&= \left(\left(\prod_{\alpha}(Y_\alpha - 1)\right)\lambda(w)J_{w^{-1}P|P}(w'\sigma\otimes\chi)\right)_{|\chi=1}\rho_{w'}^{-1}\\
&= \left(\left(\prod_{\alpha}({}^{w'^{-1}}Y_\alpha - 1)\right)\left(\prod_{\beta}(Y_\beta - 1)\right)\lambda(w)J_{w^{-1}P|P}(w'\sigma\otimes w'\chi)\right.\\
&\left.\times\lambda(w')J_{w'^{-1}P|P}(\sigma\otimes\chi)\right)_{|\chi=1}\\
&= \left(\left(\prod_{\alpha}({}^{w'^{-1}}Y_\alpha - 1)\right)\left(\prod_{\beta}(Y_\beta - 1)\right)\mu_{w,w'}^{-1}\lambda(ww')J_{w'^{-1}w^{-1}P|P}(\sigma\otimes\chi)\right)_{|\chi=1},
\end{aligned}
$$

the product over $\alpha$ taken in the set $\Sigma_{\mathcal{O},\mu}(P)\cap w^{-1}\Sigma_{\mathcal{O},\mu}(\overline{P})$ and that over $\beta$ taken in the set $\Sigma_{\mathcal{O},\mu}(P)\cap w'^{-1}\Sigma_{\mathcal{O},\mu}(\overline{P})$. Note the set equality

$$
\begin{aligned}
&(w'^{-1}\Sigma_{\mathcal{O},\mu}(P)\cap w'^{-1}w^{-1}\Sigma_{\mathcal{O},\mu}(\overline{P}))\cup(\Sigma_{\mathcal{O},\mu}(P)\cap w'^{-1}\Sigma_{\mathcal{O},\mu}(\overline{P}))\\
&=(\Sigma_{\mathcal{O},\mu}(P)\cap w'^{-1}w^{-1}\Sigma_{\mathcal{O},\mu}(\overline{P}))\\
&\cup\pm(\Sigma_{\mathcal{O},\mu}(P)\cap w'^{-1}\Sigma_{\mathcal{O},\mu}(\overline{P})\cap w'^{-1}w^{-1}\Sigma_{\mathcal{O},\mu}(P)),
\end{aligned}
$$

the union being disjoint. As ${}^{w'^{-1}}Y_\alpha = Y_{w'^{-1}\alpha}$ by the lemma **1.5**, we can indeed deduce the equality $(*)$. □

2.5

We will now fix an isomorphism $\rho_w : wE \to E$ first for $w \in R(\mathcal{O})$ and then for $w \in W(M,\mathcal{O})$.

Recall that $R(\mathcal{O})$ is a product of subgroups $R(\mathcal{O})_i$ of order 2 (cf. **1.15**(ii)). Denote by $s_{\alpha_i}$ the nontrivial element of $R(\mathcal{O})_i$. Choose for every $i$ an isomorphism $\rho_{s_{\alpha_i}} : s_{\alpha_i}E \to E$ such that $\rho_{s_{\alpha_i}}^2\sigma(s_{\alpha_i}^{-2}) = \mathrm{Id}_E$. (This is always possible, even if it means multiplying $\rho_{s_{\alpha_i}}$ by a nonzero complex number. The isomorphism $\rho_{s_{\alpha_i}}$ is uniquely determined up to product by $\pm 1$.)

Change the representative of $s_{\alpha_i}$ on the right by an element of $m \in M$ then come back to compose $\rho_{s_{\alpha_i}}$ on the right with $\sigma(m)$, as we can easily verify.

It is clear that the $\rho_{s_{\alpha_i}}$ switch between them, since the $s_{\alpha_i}$ acts on the different components of $\sigma$. For $r \in R(\mathcal{O})$, $r = s_{\alpha_{i_1}} \cdots s_{\alpha_{i_\ell}}$, we put $\rho_r = \rho_{s_{\alpha_{i_1}}} \cdots \rho_{s_{\alpha_{i_\ell}}}$. By the above, it is clear that this is independent of the decomposition chosen for $r$. Multiply the representative of $r$ on the right by an element $m \in M$ amounts to composing $\rho_r$ on the right with $\sigma(m)$. Indeed, this amounts to multiplying each element $s_{\alpha_i}$ in the decomposition of $r$ on the right by a certain element $m_i$ of $M_{s_{\alpha_i}}$ by our convention in **1.5** on the choice of representatives of elements in $R(\mathcal{O})$.

Finally, if $w \in W(M, \mathcal{O})$, we can by **1.12** write $w = rw_{\mathcal{O}}$ with $r \in R(\mathcal{O})$ and $w_{\mathcal{O}} \in W_{\mathcal{O}}$ uniquely determined by $w$, and we thus set $\rho_w = \rho_r \rho_{w_{\mathcal{O}}}$.

The isomorphism $i_P^G wE \to i_P^G E$ deduced from $\rho_w$ by functoriality will (by abuse of notation) be denoted again by $\rho_w$. Multiply the representative of $w$ on the right by an element $m$ of $M$ amounts to composing $\rho_w$ on the right by $i_P^G(\sigma(m))$.

2.6

**Definition** Let $\chi \in \mathrm{Stab}(\mathcal{O})$ and recall the isomorphism $\phi_{\sigma,\chi} : E \to E_\chi$ defined in **1.17**. We denote by $\phi_\chi$ the automorphism of $E_B$ which sends $e \otimes b$ to $(\phi_{\sigma,\chi} e) \otimes b_{\chi^{-1}}$. We write again $\phi_\chi$ for the automorphisms of $i_P^G E_B$ and of $i_P^G E_{K(B)}$ deduced from $\phi_\chi$ by functoriality.

2.7

The following lemma will be useful in the following:

**Lemma** *Let $\alpha, \alpha \in \Delta_{\mathcal{O}}$ be simple roots, $s = s_\alpha$ and $s' = s_{\alpha'}$. Let $m_s \in M_\alpha^1 \cap M$ and $m_{s'} \in M_{\alpha'}^1 \cap M$. Thus*

$$b_{m_s} {}^{s}b_{m_{s'}} {}^{ss'}b_{m_s} {}^{ss's}b_{m_{s'}} \cdots = b_{m_{s'}} {}^{s'}b_{m_s} {}^{s's}b_{m_{s'}} {}^{s'ss'}b_{m_s} \cdots ,$$

*the number of factors of each side being equal to the order of $ss'$.*

*Proof.* This is a simple calculation, by considering the different root systems of rank 2; indeed, by **1.5**, there exists integers $k$ and $k'$ such that $b_{m_s} = Y_\alpha^k$ and $b_{m_{s'}} = Y_{\alpha'}^{k'}$.

Moreover, we have ${}^{s}b_{m_s} = Y_\alpha^{-k}$, ${}^{s'}b_{m_{s'}} = Y_{\alpha'}^{-k'}$, ${}^{s'}b_{m_s} = Y_{s'(\alpha)}^k$ and ${}^{s}b_{m_{s'}} = Y_{s(\alpha')}^{k'}$. As an example, let's deal with the case where $\alpha$ and $\alpha'$ generate a root system of type $A_2$. The order of $ss'$ is thus 3, and we have $s(\alpha') = \alpha + \alpha' = s'(\alpha)$. We can deduce ${}^{s'}b_{m_s} = Y_\alpha^k Y_{\alpha'}^k$, ${}^{s}b_{m_{s'}} = Y_\alpha^{k'} Y_{\alpha'}^{k'}$, ${}^{ss'}b_{m_s} = Y_{\alpha'}^k$ and ${}^{s's}b_{m_{s'}} = Y_\alpha^{k'}$, hence

$$b_{m_s} {}^{s}b_{m_{s'}} {}^{ss'}b_{m_s} = Y_\alpha^k Y_\alpha^{k'} Y_{\alpha'}^{k'} Y_{\alpha'}^k = Y_\alpha^{k'} Y_\alpha^k Y_{\alpha'}^k Y_\alpha^{k'} = b_{m_{s'}} {}^{s'}b_{m_s} {}^{s's}b_{m_{s'}}.$$

□

**3**

Consider $\mathrm{End}_G(i_P^G E_B)$ as a left $B$-module. We will build for every $w \in W(M,\mathcal{O})$ an intertwining operator $A_w \in \mathrm{Hom}_G(i_P^G E_B, i_P^G E_{K(B)})$ which extends in a canonical way to an element of $\mathrm{End}_G(i_P^G E_{K(B)})$.

3.1

For $\chi \in \mathfrak{X}^{\mathrm{nr}}(M)$, denote by $B_\chi$ the ideal of polynomials $b \in B$ with $b(\chi) = 0$ and $sp_\chi$ the specialization map $E \otimes B \to E$ and $E \otimes K(B) \to E$, $e \otimes b \mapsto b(\chi)e$, as well as the induced maps $i_P^G E_B \to i_P^G E$ and $i_P^G E_{K(B)} \to i_P^G E$. (In the case of $E_{K(B)}$, the map $sp_\chi$ is only defined on the subspace of elements which are regular at $\chi$.) Let $w \in W(M)$. According to [Wal03, IV.1.1], there exists $b_w \in B$ and a homomorphism of $G-B$-modules

$$J_{B,w} : i_P^G E_B \to i_P^G(wE_B)$$

such that, for every $\chi \in \mathfrak{X}^{\mathrm{nr}}(M)$, $v \in i_P^G E_B$,

$$b_w(\chi)(\lambda(w)J_{w^{-1}P|P}(\sigma \otimes \chi)sp_\chi v) = sp_\chi J_{B,w} v.$$

The homomorphism $b_w^{-1} J_{B,w}$ defines an element of $\mathrm{Hom}_G(i_P^G E_B, i_P^G(wE_{K(B)}))$ that we will denote by $J_{K(B),w}$. Suppose for now $w \in W(M,\mathcal{O})$. Thus we denote by $A_w$ the element of $\mathrm{Hom}_G(i_P^G E_B, i_P^G E_{K(B)})$ obtained, we compose $J_{K(B),w}$ on the left with $\rho_w \circ \tau_w$. It extends in a canonical way to an element of the algebra $\mathrm{End}_G(i_P^G E_{K(B)})$.

**Proposition** *Let $w \in W(M,\mathcal{O})$. Multiply the representative of $w$ on the right by an element $m \in M$ amounts to multiplying the operator $A_w$ on the left by ${}^w b_m^{-1}$. In particular, $A_w$ does not depend on the choice of a representative in $K$ of $w$.*

*Proof.* By the definition above, that in **2.5** and the remark in **1.1**, we have $\rho_{wm} = \rho_w i_P^G(\sigma(m))$ and $J_{K(B),wm} = i_P^G(\sigma_B(m^{-1}))J_{K(B),w}$, hence the result. □

3.2

**Lemma** *Let $w \in W(M,\mathcal{O})$, $v \in i_P^G E_B$ and $\chi \in \mathfrak{X}^{\mathrm{nr}}(M)$ such that $\sigma \otimes \chi$ is a regular point for $J_{w^{-1}P|P}$.*

*Therefore*

$$sp_\chi A_w v = \rho_w \lambda(w) J_{w^{-1}P|P}(\sigma \otimes w^{-1}\chi) sp_{w^{-1}\chi} v.$$

*Proof.* We have

$$sp_\chi A_w v = sp_\chi \rho_w \tau_w J_{K(B),w} v = \rho_w sp_\chi \tau_w J_{K(B),w} v.$$

But, the canonical map $B \to B/B_\chi$ composed on the right with $\tau_w$ has kernel $B_{w^{-1}\chi}$. This composite map is thus equal to $sp_{w^{-1}\chi}$. The expression above is thus equal to $\rho_w sp_{w^{-1}\chi} J_{K(B),w} v$. This is indeed equal to the expression of the statement. □

3.3

**Proposition** *Let $w, w' \in W_{\mathcal{O}}$. Thus,*

$$A_w A_{w'} = {}^{ww'}\mu_{w,w'}^{-1} \left( \mu_{w,w'} \prod_\alpha (Y_\alpha - 1)^{-1} (Y_\alpha^{-1} - 1)^{-1} \right)_{|\chi=1} A_{ww'},$$

*the product taken over the set $\Sigma_{\mathcal{O},\mu}(P) \cap w'^{-1}\Sigma_{\mathcal{O},\mu}(\overline{P}) \cap w'^{-1}w^{-1}\Sigma_{\mathcal{O},\mu}(P)$.*

*Proof.* It suffices to show that, for every $\chi \in \mathfrak{X}^{\mathrm{nr}}(M)$, the two operators have the same specialisation at $\chi$. Set $w'' = ww'$. According to lemma **3.2** and the commutation property of $\rho_w$ remarked in **2.4**, the specialization at $\chi$ of the operator on the left is equal to

$$\begin{aligned}
&\rho_w \lambda(w) J_{w^{-1}P|P}(\sigma \otimes w^{-1}\chi) \rho_{w'} \lambda(w') J_{w'^{-1}P|P}(\sigma \otimes w'^{-1}w^{-1}\chi) sp_{w'^{-1}w^{-1}\chi} \\
&\quad = \rho_w \rho_{w'} \lambda(w) J_{w^{-1}P|P}(w'\sigma \otimes w^{-1}\chi) \lambda(w') J_{w'^{-1}P|P}(\sigma \otimes w'^{-1}w^{-1}\chi) sp_{w'^{-1}w^{-1}\chi} \\
&\quad = \rho_w \rho_{w'} \lambda(w) \lambda(w')^{ww'}\mu_{w,w'}^{-1}(\chi) J_{w'^{-1}w^{-1}P|P}(\sigma \otimes w'^{-1}w^{-1}\chi) sp_{w'^{-1}w^{-1}\chi} \\
&\quad = \rho_w \rho_{w'} sp_{w^{-1}\chi} \tau_{w'} \lambda(ww'w''^{-1})^{ww'}\mu_{w,w'}^{-1}(\chi) J_{K(B),ww'} \\
&\quad = \rho_w sp_{w^{-1}\chi} \rho_{w'} \tau_{w'} \lambda(ww'w''^{-1})^{ww'}\mu_{w,w'}^{-1}(\chi) J_{K(B),ww'} \\
&\quad = sp_\chi \rho_w \rho_{w'} \tau_{ww'} \lambda(ww'w''^{-1})^{ww'}\mu_{w,w'}^{-1}(\chi) J_{K(B),ww'}.
\end{aligned}$$

It remains thus to replace $\rho_w \rho_{w'}$ by the expression given in the proposition **2.4**, hence the proposition. □

3.4

**Corollary** *Let $w, w' \in W_{\mathcal{O}}$ such that $\ell_{\mathcal{O}}(ww') = \ell_{\mathcal{O}}(w) + \ell_{\mathcal{O}}(w')$. Thus $A_w A_{w'} = A_{ww'}$. Moreover, if $s = s_\alpha$ is a simple reflection in $W_{\mathcal{O}}$, then there exists a complex number $c_s'' \neq 0$ such that $A_s^2 = c_s''(\mu^{M_\alpha})^{-1}$.*

*Proof.* If $\ell_{\mathcal{O}}(ww') = \ell_{\mathcal{O}}(w) + \ell_{\mathcal{O}}(w')$, then $\Sigma_{\mathcal{O},\mu}(P) \cap w'^{-1}\Sigma_{\mathcal{O},\mu}(\overline{P}) \cap w'^{-1}w^{-1}\Sigma_{\mathcal{O},\mu}(P)$ is the empty set, and the rational function $\mu_{w,w'}$ is constant. The first assertion is thus an immediate consequence of **3.3**. The second assertion results from **3.3** and from the fact that $\mu_{s,s}$ is equal to $\mu^{M_\alpha}$ multiplied by a constant and that $A_1 = \mathrm{Id}$. □

3.5

**Proposition** (i) *For $w \in W_{\mathcal{O}}$ and $r \in R(\mathcal{O})$, we have*

$$A_r A_w = A_{rw} = A_{rwr^{-1}} A_r.$$

(ii) *For every $r, r' \in R(\mathcal{O})$, we have $A_r A_{r'} = A_{rr'}$.*
(iii) *For every $b \in B$ and every $w \in W(M, \mathcal{O})$, $A_w b = {}^w b A_w$.*

*Proof.* (i) As $r\Delta_{\mathcal{O},\mu} = \Delta_{\mathcal{O},\mu}$ for $r \in R(\mathcal{O})$, it suffices, following the corollary **3.4**, to consider the case where $w$ is a simple element $s_\alpha$ of $W_{\mathcal{O}}$. By applying the same reasoning as in the proof of the proposition **3.3**, we find for $\chi$ a generic element of $\mathfrak{X}^{\mathrm{nr}}(M)$ that

$$sp_\chi(A_r A_{s_\alpha}) = \rho_r \rho_{s_\alpha} \lambda(r)\lambda(s_\alpha) J_{s_\alpha^{-1} r^{-1} P|P}(\sigma \otimes s_\alpha^{-1} r^{-1}\chi) sp_{s_\alpha^{-1} r^{-1}\chi} = sp_\chi(A_{rs_\alpha})$$

and

$$sp_\chi(A_{rs_\alpha r^{-1}} A_r) = \rho_{s_{r\alpha}} \rho_r \lambda(s_{r\alpha})\lambda(r) J_{r^{-1} s_{r\alpha}^{-1} P|P}(\sigma \otimes r^{-1} s_{r\alpha}^{-1}\chi) sp_{r^{-1} s_{r\alpha}^{-1}\chi}.$$

These two expressions are equal, if and only if

$$\rho_r \rho_{s_\alpha} \lambda(r)\lambda(s_\alpha) = \rho_{s_{r\alpha}} \rho_r \lambda(s_{r\alpha})\lambda(r),$$

which is equivalent to

$$\rho_r^{-1} \rho_{s_{r\alpha}}^{-1} \rho_r \rho_{s_\alpha} = \lambda(s_{r\alpha} r s_\alpha^{-1} r^{-1}).$$

Now, as the left-hand side of this last equality is equal to $\rho_{r\sigma, s_{r\alpha}}^{-1} \rho_{\sigma, s_\alpha}$, this results from the last assertion of the proposition **2.4**.

(ii) It suffices to consider the case where $r = s_{\alpha_i}$ is the nonzero element of an $R(\mathcal{O})_i$. We will first consider the case where the projection of $r'$ onto $R(\mathcal{O})_i$ is zero. By the same reasoning as in the proof of the proposition **3.3**, we find for $\chi$ generic

$$sp_\chi(A_r A_{r'}) = \rho_r \rho_{r'} \lambda(r)\lambda(r') J_{r'^{-1} r^{-1} P|P}(\sigma \otimes r'^{-1} r^{-1}\chi) sp_{r'^{-1} r^{-1}\chi}.$$

As $\rho_r \rho_{r'} = \rho_{rr'}$ if the projection of $r'$ onto $R(\mathcal{O})_i$ is zero, this expression is equal to $sp_\chi(A_{rr'})$. We see also that the operators $A_{s_{\alpha_i}}$ and $A_{r'}$ commute.

Suppose for now the projection of $A_{r'}$ onto $R(\mathcal{O})_i$ is non-zero. It suffices thus, according to commutativity, to consider the case $r = r'$. For $\chi$ generic in $\mathfrak{X}^{\mathrm{nr}}(M)$, we find, by an analogous reasoning as that above,

$$sp_\chi(A_r^2) = \rho_r \rho_r \lambda(r^2) J_{P|P}(\sigma \otimes r^2\chi) sp_{r^2\chi} = sp_\chi,$$

since $r^2 \in M \cap K$ and that $\rho_r$ was chosen such that $\rho_r^2 \lambda(r^2) = \mathrm{Id}$.

(iii) This is immediate. □

3.6

**Proposition** *The space* $\mathrm{Hom}_M(E_B, E_{K(B)})$ *is isomorphic to* $\bigoplus_{\chi\in\mathrm{Stab}(\mathcal{O})} K(B)\phi_\chi$ *as* $K(B)$*-vector-spaces.*

*Proof.* We have $E_B = \mathrm{ind}_{M^1}^M E_{|M^1}$. It is therefore necessary to determine $\mathrm{Hom}_M(\mathrm{ind}_{M^1}^M E, E_{K(B)})$. Observe that $(E_{|M^1})^\vee = (E^\vee)_{|M^1}$, since $M^1$ is an open subgroup of $M$. As $M^1$ and $M$ are unimodular, Frobenius reciprocity relative to compact induction [BZ76, 2.29] thus gives $\mathrm{Hom}_M(\mathrm{ind}_{M^1}^M E_{|M^1}, E_{K(B)}) \simeq \mathrm{Hom}_{M^1}(E_{|M^1}, (E_{K(B)})_{|M^1})$. We define an isomorphism $\beta : \mathrm{Hom}_M(\mathrm{ind}_{M^1}^M E_{|M^1}, E_{K(B)}) \to \mathrm{Hom}_{M^1}(E_{|M^1}, (E_{K(B)})_{|M^1})$, $\varphi \mapsto \beta(\varphi)$, we put $\beta(\varphi)(e) = \varphi(v_e)$, where $v_e$ denotes the element of $\mathrm{ind}_{M^1}^M E$ with value 0 for $m \in M - M^1$ and with value $e$ at 1. Denote $E_{\overline{m}}, \overline{m} \in M/M^\sigma$, the different irreducible components of $E_{|M^1}$. They are pairwise non-isomorphic according to the proposition **1.16**. It follows that $\mathrm{Hom}_{M^1}(E_{|M^1}, (E_{K(B)})_{|M^1}) = \bigoplus_{\overline{m}\in M/M^\sigma} \mathrm{Hom}_{M^1}(E_{\overline{m}}, E_{\overline{m}} \otimes K(B))$. Choose a $\mathbb{C}$-basis $f_i$, $i \in I$, of $K(B)$. Therefore, as $M^1$-module, $E_{\overline{m}} \otimes K(B) = \bigoplus_{i\in I} E_{\overline{m}} \otimes f_i$. We thus deduce from Schur lemma that the elements of $\mathrm{Hom}_{M^1}(E_{\overline{m}}, E_{\overline{m}} \otimes K(B))$ are of the form $e \mapsto e\otimes f$ with $f \in K(B)$. As a result, $\mathrm{Hom}_{M^1}(E_{|M^1}, (E_{K(B)})_{|M^1})$ is a $K(B)$-module isomorphic to $K(B)^{M/M^\sigma}$.

To prove the proposition, it thus suffices, thanks to the $K(B)$-linearity of $\beta$, to show that the $\beta(\phi_\chi)$ form a basis of $K(B)$-module $\mathrm{Hom}_{M^1}(E_{|M^1}, (E_{K(B)})_{|M^1})$. We identify $\mathrm{ind}_{M^1}^M E_{|M^1}$ and $E\otimes B$. We have $\beta(\phi_\chi)(e) = \phi_\chi(e\otimes 1) = \phi_{\sigma,\chi}(e)\otimes 1$. The projection of $\beta(\phi_\chi)$ onto $\mathrm{Hom}_{M^1}(E_{\overline{m}}, E_{\overline{m}}\otimes K(B))$ is thus the multiplication by $\chi(m)$. We thus deduce from the duality of groups $\mathrm{Stab}(\mathcal{O})$ and $M/M^\sigma$ and from the linear independence of characters that the $\beta(\phi_\chi)$ form indeed a basis of $K(B)$-module $\mathrm{Hom}_{M^1}(E_{|M^1}, (E_{K(B)})_{|M^1})$. □

3.7

**Proposition** *The intertwining operators* $\phi_\chi A_w$*, where* $\chi \in \mathrm{Stab}(\mathcal{O})$ *and* $w \in W(M, \mathcal{O})$*, are* $K(B)$*-linearly independent in* $\mathrm{Hom}_G(i_P^G E_B, i_P^G E_{K(B)})$ *equipped with the structure of* $K(B)$*-vector space induced from that of* $i_P^G E_{K(B)}$*.*

*Proof.* Suppose for the sake of contradiction that the operators $\phi_\chi A_w$ are $K(B)$-linearly dependent. There thus exists elements $b_{w,\chi}$ of $B$ such that $\sum_{w\in W(M,\mathcal{O}),\chi\in\mathrm{Stab}(\mathcal{O})} b_{w,\chi}\phi_\chi A_w = 0$, and we can suppose the number $\ell$ of elements $b_{w,\chi} \neq 0$ to be minimal. As none of the $\phi_\chi A_w$ is zero, there exist at least two indices $(w_1, \chi_1)$ and $(w_2, \chi_2)$ with $b_{w_1,\chi_1} \neq 0 \neq b_{w_2,\chi_2}$. We can find $b \in B$ with $({}^{w_1}b)_{\chi_1^{-1}} = b$ and $({}^{w_2}b)_{\chi_2^{-1}} \neq b$. We have $\sum_{w\in W(M,\mathcal{O})} bb_{w,\chi}\phi_\chi A_w = 0$ and $0 = \sum_{w\in W(M,\mathcal{O})} b_{w,\chi}\phi_\chi A_w b = \sum_{w\in W(M,\mathcal{O})} b_{w,\chi}({}^w b)_{\chi^{-1}}\phi_\chi A_w$, where we applied the easy commutation rule $A_w b = {}^w b A_w$. Subtracting these two equations, we obtain a linear combination of $\phi_\chi A_w$ which is zero and whose fewer than $\ell$ coefficients are non-zero. This gives a contradiction. □

3.8

**Theorem** *As $K(B)$-vector-spaces, we have*

$$\mathrm{Hom}_G(i_P^G E_B, i_P^G E_{K(B)}) = \bigoplus_{w\in W(M,\mathcal{O}),\chi\in \mathrm{Stab}(\mathcal{O})} K(B)\phi_\chi A_w.$$

*Proof.* By Frobenius reciprocity, we have

$$\mathrm{Hom}_G(i_P^G E_B, i_P^G E_{K(B)}) = \mathrm{Hom}_M(r_P^G i_P^G E_B, E_{K(B)}).$$

By the geometric lemma [BZ77, 2.12], $r_P^G i_P^G E_B$ admits a filtration by the subspaces $\mathcal{F}_w$, $w \in W^M\backslash W^G/W^M$, whose subquotients are isomorphic to $wE_B$. We deduce a filtration of $K(B)$-module $\mathrm{Hom}_M(r_P^G i_P^G E_B, E_{K(B)})$ whose subquotients are, considering the fact **3.6**, isomorphic to $\mathrm{Hom}_M(wE_B, E_{K(B)})$.

As $\mathrm{Hom}_M(wE_B, E_{K(B)}) \neq 0$, if and only if $w \in W(M,\mathcal{O})$ and thus, by the proposition **3.6**, $\mathrm{Hom}_M(wE_B, E_{K(B)})$ is of dimension $|\mathrm{Stab}(\mathcal{O})|$, we deduce that the dimension of $\mathrm{Hom}_G(i_P^G E_B, i_P^G E_{K(B)})$ is equal to $|W(M,\mathcal{O})||\mathrm{Stab}(\mathcal{O})|$.

As the intertwining operators $\phi_\chi A_w$, $w \in W(M,\mathcal{O})$, $\chi \in \mathrm{Stab}(\mathcal{O})$ are linearly independent by the proposition **3.7**, the theorem follows. □

## 4

For now, we fix an irreducible component $E_1$ of $E_{|M^1}$, and we denote by $\sigma_1$ the representation of $M^1$ in this space. We will be interested in the algebra $\mathrm{End}_G(i_P^G(\mathrm{ind}_{M^1}^M E_1))$. Note that neither this algebra nor the representation $\mathrm{ind}_{M^1}^M E_1$ does not depend on the choice of $(\sigma_1, E_1)$. Denote by $B_{\mathcal{O}}$ the subalgebra of $B$ formed of polynomials which are invariant under translation by elements of $\mathrm{Stab}(\mathcal{O})$. According to corollary **1.17**, this is the ring of regular functions of the quotient affine variety $X(M/M^1)/X(M/M^\sigma)$. It is thus identified with $\mathbb{C}[M^\sigma/M^1]$. In particular, this is a factorial ring, since $M^\sigma/M^1$ is a free $\mathbb{Z}$-module of finite type and of the same rank as $M/M^1$. Recall that we denoted $\mathcal{R}(M/M^1)$ a system of representatives of $M/M^1$.

4.1

**Lemma** *The canonical isomorphism* $\operatorname{ind}_{M^1}^{M} E_{|M^1} \to E \otimes B$ *of* ***2.1*** *sends* $\operatorname{ind}_{M^1}^{M} E_1$ *to the set* $\sum_{m \in \mathcal{R}(M/M^1)} \sigma(m) E_1 \otimes b_m$. *In particular, this space is a* $B_{\mathcal{O}}$*-submodule of* $E \otimes B$.

*Proof.* This is immediate. □

We identify in the following $\operatorname{ind}_{M^1}^{M} E_1$ and $\sum_{m \in \mathcal{R}(M/M^1)} \sigma(m) E_1 \otimes b_m$ by means of the isomorphism in **2.1**. We write $E_{B_{\mathcal{O}}}$ for $\operatorname{ind}_{M^1}^{M} E_1$, although this space is not isomorphic to $E \otimes B_{\mathcal{O}}$.

4.2

Denote by $K(B_{\mathcal{O}})$ the fraction field of $B_{\mathcal{O}}$ and set $(E_{B_{\mathcal{O}}})_{K(B_{\mathcal{O}})} = E_{B_{\mathcal{O}}} \otimes_{B_{\mathcal{O}}} K(B_{\mathcal{O}})$.

**Lemma** *We have*

$$\operatorname{Hom}_M(E_{B_{\mathcal{O}}}, E_{K(B_{\mathcal{O}})}) \simeq K(B_{\mathcal{O}}).$$

*Proof.* As in the proof of proposition **3.6**, we deduce from Frobenius reciprocity that

$$\begin{aligned}
&\operatorname{Hom}_M(\operatorname{ind}_{M^1}^{M} E_1, (\operatorname{ind}_{M^1}^{M} E_1)_{K(B_{\mathcal{O}})}) \\
&= \operatorname{Hom}_{M^1}(E_1, (\operatorname{ind}_{M^1}^{M} E_1 \otimes_{B_{\mathcal{O}}} K(B_{\mathcal{O}}))_{|M^1}) \\
&= \operatorname{Hom}_{M^1}(E_1, \left( \bigoplus_{m \in \mathcal{R}(M/M^{\sigma})} \sigma(m) E_1 \otimes b_m \right) \otimes_{B_{\mathcal{O}}} K(B_{\mathcal{O}})) \\
&= \operatorname{Hom}_{M^1}(E_1, E_1 \otimes K(B_{\mathcal{O}})) \\
&= K(B_{\mathcal{O}}).
\end{aligned}$$

□

4.3

**Proposition** *We have*

$$\bigcap_{\eta \in \operatorname{Stab}(\mathcal{O})} \ker(\phi_\eta - \operatorname{Id}) = (E_{B_{\mathcal{O}}})_{K(B_{\mathcal{O}})}.$$

*In particular, the automorphism* $\phi_\eta$ *is trivial on* $E_{B_{\mathcal{O}}}$.

*Proof.* Fix the system of representatives $\mathcal{R} = \mathcal{R}(M/M^1)$ of $M/M^1$ and $\mathcal{R}' = \mathcal{R}(M/M^{\sigma})$ of $M/M^{\sigma}$ respectively. Every element of $(E_{B_{\mathcal{O}}})_{K(B_{\mathcal{O}})}$ is written in the form $b^{-1}\sum_{m\in\mathcal{R}}\sigma(m)e_m\otimes b_m$ with $e_m\in E_1$ and $b\in B_{\mathcal{O}}$. We verify right away that such an element is in the kernel of $\phi_\eta-\mathrm{Id}$ for every $\eta$.

Conversely, let $v\in E_{K(B)}$. It can be written in the form $b^{-1}\sum_{(m,m')\in\mathcal{R}\times\mathcal{R}'}\sigma(m')e_{m,m'}\otimes b_m$ with $e_{m,m'}\in E_1$ and $b\in B$. Let $\eta\in\mathrm{Stab}(\mathcal{O})$. Therefore,

$$\phi_\eta(v)=b_{\eta^{-1}}^{-1}\sum_{(m,m')\in\mathcal{R}\times\mathcal{R}'}\sigma(m')\eta(m')\eta(m)^{-1}e_{m,m'}\otimes b_m.$$

Thus, $\phi_\eta(v)=v$ implies that

$$b_{\eta^{-1}}\sum_{(m,m')\in\mathcal{R}\times\mathcal{R}'}\sigma(m')e_{m,m'}\otimes b_m=b\sum_{(m,m')\in\mathcal{R}\times\mathcal{R}'}\sigma(m')\eta(m'm^{-1})e_{m,m'}\otimes b_m.$$

Consequently, for every $m'\in\mathcal{R}'$,

$$(*)\qquad b_{\eta^{-1}}\sum_{m\in\mathcal{R}}e_{m,m'}\otimes b_m=b\sum_{m\in\mathcal{R}}\eta(m'm^{-1})e_{m,m'}\otimes b_m.$$

We will first show that every primary component of $b$ divides $b_{\eta^{-1}}$. As $b$ and $b_{\eta^{-1}}$ are Laurent monomials of the same degree, it will follow that they do not differ more than by a constant, and it will thus be clear that $b^{-1}b_{\eta^{-1}}=\eta(m'm^{-1})$ if $e_{m,m'}\neq 0$.

Let $b_0$ be an irreducible factor of $b$. We can suppose that $b_0$ does not divide $\sum_{(m,m')\in\mathcal{R}\times\mathcal{R}'}\sigma(m')e_{m,m'}\otimes b_m$ in $E_B$. There thus exists $m'\in\mathcal{R}'$ such that $b_0$ does not divide $\sum_{m\in\mathcal{R}}e_{m,m'}\otimes b_m$. Fix a basis $(e_i)_{i\in I}$ of $E_1$ and write $e_{m,m'}=\sum_i c_i(m,m')e_i$ with $c_i(m,m')\in\mathbb{C}$. Thus the equality $(*)$ implies that

$$b_{\eta^{-1}}\sum_{m\in\mathcal{R}}\sum_{i\in I}c_i(m,m')e_i\otimes b_m=b\sum_{m\in\mathcal{R}}\sum_{i\in I}c_i(m,m')\eta(m'm^{-1})e_i\otimes b_m.$$

We can deduce $b_{\eta^{-1}}\sum_{m\in\mathcal{R}}c_i(m,m')b_m=b\sum_{m\in\mathcal{R}}c_i(m,m')\eta(m'm^{-1})b_m$ for every $i$. Moreover, there exists at least an $i$ such that $b_0$ does not divide $\sum_{m\in\mathcal{R}}c_i(m,m')b_m$ in $B$. As a result, the $b_0$-primary component of $b$ divides $b_{\eta^{-1}}$.

We have thus $b^{-1}b_{\eta^{-1}}=\eta(m'm^{-1})$, if $e_{m,m'}\neq 0$. Suppose that there exists $(m_0,m_0')\in\mathcal{R}\times\mathcal{R}'$ such that $e_{m_0,m_0'}\neq 0$. Set $h=m_0'm_0^{-1}$. Therefore, $\frac{b(\chi\eta^{-1})}{b(\chi)}$ is constant in $\chi$ with value $\eta(h)$.

Write $b=\sum_{m\in M/M^1}c_mb_m$ with $c_m\in\mathbb{C}$. Therefore $b_{\eta^{-1}}=\sum_{m\in M/M^1}c_m\eta(m)^{-1}b_m$. As $b_{\eta^{-1}}=\eta(h)b$, we find

$$\sum_{m\in M/M^1}c_m\eta(m)^{-1}b_m=\sum_{m\in M/M^1}\eta(h)c_mb_m.$$

Considering $M/M^1$ as group of characters of $\mathfrak{X}^{\mathrm{nr}}(M)$, we deduce from the linear independence of characters that $\eta(h)c_m = \eta^{-1}(m)c_m$ for every $m$ and every $\eta \in \mathrm{Stab}(\mathcal{O})$. Thus, $c_m \neq 0$ implies that $m \in h^{-1}M^\sigma$, hence $b = \sum_{m\in M^\sigma/M^1} c_m b_{h^{-1}m} = b_{h^{-1}}b'$ with $b' := \sum_{m\in M^\sigma/M^1} c_m b_m \in B_{\mathcal{O}}$.

As $e_{m,m'} \neq 0$ implied by what preceded that $\eta(m'm^{-1}) = \eta(h)$ for every $\eta \in \mathrm{Stab}(\mathcal{O})$, we have thus $m'm^{-1} \in hM^\sigma$. We can deduce $v = b'^{-1}b_h \sum_{m\in\mathcal{R}} \sigma(mh)e_m \otimes b_m$ with $e_m \in E_1$. However, this is $b'^{-1} \sum_{m\in\mathcal{R}} \sigma(mh)e_m \otimes b_{mh}$, and, consequently, $v$ has the indicated form. □

4.4

**Lemma** *Suppose $M$ is maximal Levi of $G$ and $W_{\mathcal{O}} \neq 1$. Let $s$ be the unique element $\neq 1$ of $W_{\mathcal{O}}$ and $(\sigma^1, E^1)$ an irreducible component of $\sigma|_{M\cap G^1}$. Thus $s\sigma^1 \simeq \sigma^1$.*

*Proof.* As $K$ and $U$ are contained in $G^1$, the space ${}^K_{P\cap K}E^1$ is well-defined as $G^1$-invariant subspace of $i^K_{P\cap K}E$. Set $\Sigma_{\mathcal{O}} = \{\pm\alpha\}$. As $\overline{U} \subseteq G^1$, it follows directly from the definition of intertwining operators that the operators $J_{\overline{P}|P}(\sigma\otimes\chi_{\lambda\alpha})$, $\lambda \in \mathbb{C}$, sends the space $i^K_{P\cap K}E^1$ to thespace $i^K_{\overline{P}\cap K}E^1$, if $\lambda$ is regular. The operator $J_{P|\overline{P}}$ possesses the analogous property. As the composite $J_{P|\overline{P}}(\sigma\otimes\chi_{\lambda\alpha})J_{\overline{P}|P}(\sigma\otimes\chi_{\lambda\alpha})$ is scalar, it has a pole of order 2 at $\lambda = 0$, and the operators $J_{P|\overline{P}}(\sigma\otimes\chi_{\lambda\alpha})$ and $J_{\overline{P}|P}(\sigma\otimes\chi_{\lambda\alpha})$ possesses at most poles of order 1, it follows that the restriction of $J_{\overline{P}|P}(\sigma\otimes\chi_{\lambda\alpha})$ to $i^K_{P\cap K}E^1$ admits a pole of order 1 at $\lambda = 0$. We will deduce that $s\sigma^1 \simeq \sigma^1$. More precisely, we will show that $J_{\overline{P}|P}(\sigma\otimes\chi_{\lambda\alpha})$ is regular at $\lambda = 0$, if $s\sigma^1 \not\simeq \sigma^1$. If $G^1 = G$, this is proven in [Wal03, IV.2.2], since $\sigma$ is unitary cuspidal and that $s$ normalizes $M$. As in general $G^1 \neq G$ and that it is not necessarily a reductive group, but only a locally profinite group, we must generalize this proof, which requires some preliminary preparations.

To simplify, set $P^1 = (G^1\cap M)U$. Note first that in the theory of smooth representations of locally profinite groups [BZ76, 2.], the representation of $G^1$ in $i^K_{P\cap K}E^1$ is the representation induced by the representation $(\delta_P)^{1/2}_{|G^1\cap M}\sigma^1$ trivially extended to $(G^1\cap M)U$. We will denote $i^{G^1}_{P^1}\sigma^1$. Set $V^1 = i^K_{P\cap K}E^1$. Since $\overline{U} \subseteq G^1$, the vector space $V^1(\overline{U})$ generated by the elements of the form $v^1 - (i^{G^1}_{(G^1\cap M)U}\sigma^1)(\overline{u})v^1$, $v^1 \in V^1$, $\overline{u} \in \overline{U}$, is contained in $V^1$. Moreover, it is $G^1\cap M$-invariant. It follows that the Jacquet functor $r^G_{\overline{P}}$ sends $V^1$ to a $G^1\cap M$-invariant-subspace of $r^G_{\overline{P}}i^K_{P\cap K}E$ that it provides the action of $G^1\cap M$ given by the representation $r^G_{\overline{P}}i^G_P\sigma$. We write $r^{G^1}_{\overline{P}^1}V^1$ and $r^{G^1}_{\overline{P}^1}i^{G^1}_{P^1}\sigma^1$. We verify the Frobenius reciprocity

$$\mathrm{Hom}_{G^1}(i_{P^1}^{G^1}\sigma^1, i_{\overline{P}^1}^{G^1}\sigma^1) = \mathrm{Hom}_{G^1\cap M}(r_{\overline{P}^1}^{G^1} i_{P^1}^{G^1}\sigma^1, \sigma^1), \qquad \Phi \mapsto r_{\overline{P}^1}^{G^1}\Phi,$$

hence $\left(r_{\overline{P}^1}^{G^1}\Phi\right)\left(r_{\overline{P}^1}^{G^1}v^1\right) := (\Phi(v^1))(1)$. The ring $B^1 = \mathbb{C}[M\cap G^1/M^1]$ is identified with the ring of regular functions on $\mathfrak{X}^{\mathrm{nr}}(M\cap G^1) := \mathfrak{X}^{\mathrm{nr}}(M)_{|M\cap G^1}$. Set $E_{B^1}^1 = E^1\otimes B^1$ and denote by $\sigma_{B^1}^1$ the representation of $M\cap G^1$ in $E_{B^1}^1$ given by $\sigma_{B^1}^1(m^1)(e^1\otimes b^1) = \sigma^1(m^1)e^1\otimes b^1 b_{m_1}$. The above results apply equally to the representation $V_{B^1}^1 := i_{P^1}^{G^1}E_{B^1}^1$, considered as $G^1$-invariant subspace of $i_P^G E_B$. Denote by $V^2$ and $V_{B^1}^2$ the vector sub-spaces of $V^1$ (resp. $V_{B^1}^1$) formed of elements with support in $P^1\overline{P}^1$. The map $p_1 : V^2 \to E^1$, $v^2 \mapsto \int_{\overline{U}} v^2(\overline{u})d\overline{u}$, induces an isomorphism $r_{\overline{P}^1}^{G^1}V^2 \to E^1$ of representations of $M\cap G^1$. The map $p_s : V^1 \to sE^1$, $v^1 \mapsto v^1(s)$, defines an isomorphism $r_{\overline{P}^1}^{G^1}V^1/r_{\overline{P}^1}^{G^1}V^2 \to sE^1$. We have the analogous situation for $V_{B^1}^1$ and $V_{B^1}^2$.

We will now proceed to the proof that $J_{\overline{P}|P}(\sigma\otimes\chi_{\lambda\alpha})$ is regular at $\lambda = 0$, if $s\sigma^1 \not\simeq \sigma^1$: let $n$ be the smallest integer such that $\lambda^n J_{\overline{P}|P}(\sigma\otimes\chi_{\lambda\alpha})$ is regular at $\lambda = 0$. Denote this operator $J(\lambda)$ and consider the operator $r_{\overline{P}^1}^{G^1}J(0)$ in $\mathrm{Hom}_{M\cap G^1}(r_{\overline{P}^1}^{G^1} i_{P^1}^{G^1}\sigma^1, \sigma^1)$, obtained by Frobenius reciprocity. As the quotient $r_{\overline{P}^1}^{G^1}V^1/r_{\overline{P}^1}^{G^1}V^2$ is isomorphic to $s\sigma^1$ and that this representation is by assumption not isomorphic to $\sigma^1$, the restriction of $r_{\overline{P}^1}^{G^1}J(0)$ to $r_{\overline{P}^1}^{G^1}V^2$ is nontrivial. There thus exists $v^2 \in V_{B^1}^2$ such that

$$0 \neq (r_{\overline{P}^1}^{G^1}J(0))(r_{\overline{P}^1}^{G^1}sp_1v^2) = (J(0)sp_1v^2)(1).$$

However, the map $(J_{\overline{P}|P}(\sigma\otimes\chi_{\lambda\alpha})v^2)(1)$ is regular at every $\lambda\in\mathbb{C}$ [Wal03, p. 283], since $v^2$ is supported in $P\overline{P}$. Consequently, $n = 0$. □

4.5

**Lemma** *For every $w \in W(M,\mathcal{O})$, there exists $m_w$ in $M$ such that $b_{m_w}A_w$ leaves the space $(i_P^G E_{B_{\mathcal{O}}})_{K(B_{\mathcal{O}})}$ invariant. If $w = s_\alpha$ with $\alpha\in\Delta_{\mathcal{O},\mu}$ where $w\in R(\mathcal{O})_i$, we can choose $m_w$ in $M\cap M_\alpha^1$.*

*Proof.* By definition, $A_w = \rho_w\tau_w b_w^{-1}J_{B,w}$ with $b_w\in B^\times$. Indeed, we can choose $b_w \in B_{\mathcal{O}}^\times$, since the value of $J_{w^{-1}P|P}$ at a representation $\sigma'$ only depends in a certain sense on the isomorphism class of $\sigma'$ [Wal03, IV.1]. The operator $\lambda(w)^{-1}J_{B,w}$ depending only on the isomorphism class of $\sigma$ [Wal03, IV.1], it commutes with the operators $\phi_\eta$, $\eta\in\mathrm{Stab}(\mathcal{O})$. We thus deduce from proposition **4.3** that $J_{B,w}$ sends the space $i_P^G E_{B_{\mathcal{O}}}$ to the space $(i_P^G wE_{B_{\mathcal{O}}})_{K(B_{\mathcal{O}})}$. Consider the restriction of the isomorphism $\rho_w\tau_w : i_P^G wE_B \to i_P^G E_B$ to $i_P^G wE_{B_{\mathcal{O}}}$.

Let $v := \sum_{m\in\mathcal{R}(M/M^1)} \sigma(m)e_m \otimes b_m$ in $E_{B_\mathcal{O}}(e_m \in E_1)$. We have

$$\tau_w(v) = \sum_{m\in\mathcal{R}(M/M^1)} \sigma(m)e_m \otimes b_{wmw^{-1}} = \sum_{m\in\mathcal{R}(M/M^1)} (w\sigma)(m)e_{w^{-1}mw} \otimes b_m.$$

The isomrphism $\rho_w$ sends $(w\sigma)(m)e_{w^{-1}mw}$ to an element in $\sigma(m)\rho_{\sigma,w}E_1$, hence

$$\rho_w\tau_w v \in \sum_{m\in\mathcal{R}(M/M^1)} \sigma(m)\rho_{\sigma,w}E_1 \otimes b_m.$$

The space $\rho_{\sigma,w}E_1$ is the subspace of $E_{|M^1}$ equipped with the representation $w^{-1}\sigma_1$. By **1.16**, there exists $m_w \in M$ such that $\sigma(m_w)\rho_{\sigma,w}E_1$ is equal to $E_1$. It follows that the operator $A_w$ leaves $i_P^G E_{K(B_\mathcal{O})}$ invariant, if we multiply the representative of $w$ in $G$ on the right by $w^{-1}m_w w$. Now, this amounts to **3.1** by multiplying $A_w$ by $b_{m_w}^{-1}$.

Suppose for now $\alpha \in \Delta_{\mathcal{O},\mu}$ and $s = s_\alpha \in W_\mathcal{O}$. Denote by $\sigma^1$ the irreducible component of $\sigma_{|M\cap M_\alpha^1}$ whose restriction to $M^1$ contains $\sigma_1$. As $s\sigma^1 \simeq \sigma^1$ by lemma **4.4**, there exists $m_s \in M \cap M_\alpha^1$ such that $s\sigma_1 \simeq m_s\sigma_1$. The operator $b_{m_s}A_s$ leaves the space $i_P^G E_{K(B_\mathcal{O})}$ invariant, by the above.

If $w \in R(\mathcal{O})_i$, thus the effect of $w$ on $\sigma$ concerns the only component $\sigma_i$ of $\sigma$ with $w\sigma_i \simeq \sigma_i$. We can thus choose $m_w$ in the group $\mathrm{GL}_{k_i}$ on which $\sigma_i$ is defined. It therefore belongs to a component of $M_\alpha$ which is in the notations of the proposition **1.13** equal to $\underline{H}_{k_i+k}$. As this component is semi-simple, $m_w$ belongs to $M_\alpha^1$. □

4.6

**Definition** If $s$ is a simple symmetry in $W_\mathcal{O}$, we fix an element $m_s$ of $M\cap M_\alpha^1$ such that the conclusions of lemma **4.5** are verified and set $J_s = b_{m_s}A_s$. This is an element of $\mathrm{Hom}_G(i_P^G E_{B_\mathcal{O}}, i_P^G E_{K(B_\mathcal{O})})$. For $r \in R(\mathcal{O})_i$, fix $m_r \in R(\mathcal{O})$ such that the conclusions of Lemma **4.5** relative to $A_r$ are verified and set $J_r = b_{m_r}A_r$. For $w \in W_\mathcal{O}$, $w = s_1\cdots s_\ell$, define $J_w = J_{s_1}\cdots J_{s_\ell}$. (It will follow from the lemma below that this definition is correct.) For $r \in R(\mathcal{O})$, we do the same from $J_{s_{\alpha_i}}$, $\alpha_i \in R(\mathcal{O})_i$. If finally $w \in W(M,\mathcal{O})$, $w = rw_\mathcal{O}$, with $r \in R(\mathcal{O})$ and $w_\mathcal{O} \in W_\mathcal{O}$, we write $J_w = J_rJ_{w_\mathcal{O}}$.

4.7

**Lemma** (i) *Let $w \in W_\mathcal{O}$ and let $w = s_1\cdots s_\ell$ be a decomposition of $w$ into simple symmetries. The operator $J_{s_1}\cdots J_{s_\ell}$ only depends on $w$ and not on the decomposition of $w$.*
(ii) *If $s = s_\alpha$ is a simple reflection in $W_\mathcal{O}$, thus $J_s^2 = c_s''(\mu^{M_\alpha})^{-1}$, where $c_s''$ is*

*the complex number defined in* ***3.4***.
(iii) *For every* $r, r' \in R(\mathcal{O})$*, we have* $J_r J_{r'} = J_{rr'}$.

*Proof.* For (i), according to [Spr98, 8.3.3], it suffices to show the following assertion: if $s$ and $s'$ are simple symmetries in $W_{\mathcal{O}}$ and if $m(s,s')$ denotes the orderof $ss'$, thus

$$J_s J_{s'} J_s \cdots = J_{s'} J_s J_{s'} \cdots ,$$

the number of factors of each side being equal to $m(s,s')$. Considering the composition result established for the operators $A_s$ in **3.5** and the definitions of $J_s$ and of $J_{s'}$, this reduces to proving that

$$b_{m_s} {}^{s}b_{m_{s'}} {}^{ss'}b_{m_s} {}^{ss's}b_{m_{s'}} \cdots = b_{m_{s'}} {}^{s'}b_{m_s} {}^{s's}b_{m_{s'}} {}^{s'ss'}b_{m_s} \cdots$$

Now, this results from Lemma **2.7**.

For (ii), the proof of **3.4** remains valid, since $J_s^2$ is the restriction of $A_s^2$.

By iteration, it suffices to verify the equality (iii) when $r = s_{\alpha_i}$ is the non-trivial element of an $R(\mathcal{O})_i$. If the projection of $r'$ onto $R(\mathcal{O})_i$ is trivial, thus the assertion results directly from the definition **4.6**. If not, it suffices to show that $J_r^2 = 1$. By **3.5**, we have

$$J_r^2 = b_{m_r} A_r b_{m_r} A_r = b_{m_r} {}^{r}b_{m_r}.$$

As ${}^{r}b_{m_r} = b_{m_r}^{-1}$, we indeed have $b_{m_r} {}^{r}b_{m_r} = 1$. □

4.8

**Lemma** *The operators* $J_w$*,* $w \in W(M, \mathcal{O})$*, are* $K(B_{\mathcal{O}})$*-linearly independent in* $\mathrm{Hom}_G(i_P^G E_{B_{\mathcal{O}}}, i_P^G E_{K(B_{\mathcal{O}})})$*, and, we have, for* $b$ *in* $K(B_{\mathcal{O}})$*, the commutation relation* $J_w b = {}^{w}b J_w$.

*Proof.* The proof of **3.7** generalizes, after noticing that the operators $J_w$ are nonzero. Now, this results from the fact that $sp_\chi$ sends the space $E_{B_{\mathcal{O}}}$, considered as subspace of $E_B$, to the space $E$ equipped with the representation $\sigma \otimes \chi$. As $sp_\chi \circ J_{w,K(B)} = \lambda(w) J_{w^{-1}P|P}(\sigma \otimes \chi) sp_\chi$ and that $J_{w^{-1}P|P}(\sigma \otimes \chi) \neq 0$ for $\chi$ regular, the operator $J_w$ is indeed nonzero.

The commutation relation is an immediate consequence of the definition of $J_w$ and from the commutation relation for the operators $A_w$ (cf. **3.5**). □

4.9

**Theorem** *As* $K(B_{\mathcal{O}})$*-modules,*

$$\mathrm{Hom}_G(i_P^G E_{B_{\mathcal{O}}}, i_P^G E_{K(B_{\mathcal{O}})}) = \bigoplus_{w \in W(M,\mathcal{O})} K(B_{\mathcal{O}}) J_w.$$

*Proof.* By Frobenius reciprocity, we have

$$\mathrm{Hom}_G(i_P^G E_{B_\mathcal{O}}, i_P^G E_{K(B_\mathcal{O})}) = \mathrm{Hom}_M(r_P^G i_P^G E_{B_\mathcal{O}}, E_{K(B_\mathcal{O})})$$

By the geometric lemma [BZ77, 2.12], $r_P^G i_P^G E_{B_\mathcal{O}}$ admits a filtration by the subspaces $\mathcal{F}_w$, $w \in W^M\backslash W^G/W^M$, whose subquotients are isomorphic to $wE_{B_\mathcal{O}}$. We deduce a filtration of $K(B_\mathcal{O})$-module $\mathrm{Hom}_M(r_P^G i_P^G E_{B_\mathcal{O}}, E_{K(B_\mathcal{O})})$ whose subquotients are isomorphic to $\mathrm{Hom}_M(wE_{B_\mathcal{O}}, E_{K(B_\mathcal{O})})$.

Given the $K(B_\mathcal{O})$-linear independence of operators $J_w$, $w \in W(M,\mathcal{O})$, it remains to prove that

$$\mathrm{Hom}_G(wE_{B_\mathcal{O}}, (E_{B_\mathcal{O}})_{K(B_\mathcal{O})}) \simeq \begin{cases} 0, & \text{if } w \notin W(M,\mathcal{O}); \\ K(B_\mathcal{O}), & \text{if not.} \end{cases}$$

As $wE_{B_\mathcal{O}}$ is a quotient of $wE_B$ and that $(E_{B_\mathcal{O}})_{K(B_\mathcal{O})}$ is a submodule of $E_{K(B)}$, $\mathrm{Hom}_M(wE_B, E_{K(B)})$ being zero for $w \notin W(M,\mathcal{O})$ implies that $\mathrm{Hom}_M(wE_{B_\mathcal{O}}, (E_{B_\mathcal{O}})_{K(B_\mathcal{O})})$ is zero.

Suppose for now $w \in W(M,\mathcal{O})$. We have $wE_{B_\mathcal{O}} = \mathrm{ind}_{M^1}^M wE_1$. As $W(M,\mathcal{O})$ permutes the irreducible components of $E_{|M^1}$, we find that $\mathrm{ind}_{M^1}^M wE_1 = E_{B_\mathcal{O}}$, hence, according to **4.2**, $\mathrm{Hom}_M(wE_{B_\mathcal{O}}, E_{K(B_\mathcal{O})}) = K(B_\mathcal{O})$. □

## 5

In this section, we will define for every $w \in W_\mathcal{O}$ an operator $T_w \in \mathrm{End}_G(i_P^G E_{B_\mathcal{O}})$ and show that the operators $J_r T_w$, $r \in R(\mathcal{O})$ and $w \in W_\mathcal{O}$, form a basis of $B_\mathcal{O}$-module $\mathrm{End}_G(i_P^G E_{B_\mathcal{O}})$.

### 5.1

**Lemma** *Let $\alpha \in \Delta_\mathcal{O}$, $s = s_\alpha$. Set $c_s = c_s'^{-1} c_s''$ (where $c_s'$ and $c_s''$ are complex numbers defined respectively in **1.5** and in **3.4**).*

(i) *The operator $sp_1(X_\alpha - 1)J_s$ is scalar on $i_{P\cap K}^K E$ and*

$$(sp_1(X_\alpha - 1)J_s)^2 = c_s \frac{(1-q^{-a_s})^2(1+q^{-b_s})^2}{4}.$$

*Moreover, if $\chi \in \mathfrak{X}^{\mathrm{nr}}(M)$ satisfies $X_\alpha(\chi) = 1$, thus, for $v \in i_{P\cap K}^K E_{B_\mathcal{O}}$,*

$$sp_\chi(X_\alpha - 1)J_s v = (sp_1(X_\alpha - 1)J_s) sp_\chi v.$$

(ii) *Suppose that $\mu^{M_\alpha}$ cancels out in $X_\alpha = -1$. Fix $\chi_{-1} \in \mathfrak{X}^{\mathrm{nr}}(M)$ such that $X_\alpha(\chi_{-1}) = -1$. The operator $sp_{\chi_{-1}}(X_\alpha + 1)J_s$ is scalar on $i_{P\cap K}^K E$ and*

$$(sp_{\chi_{-1}}(X_\alpha + 1)J_s)^2 = c_s \frac{(1+q^{-a_s})^2(1-q^{-b_s})^2}{4}.$$

*Moreover, if* $\chi \in \mathfrak{X}^{\mathrm{nr}}(M)$ *satisfies* $X_\alpha(\chi) = -1$*, thus, for* $v \in i_{P\cap K}^K E_{B_{\mathcal{O}}}$*,*

$$sp_\chi (X_\alpha + 1) J_s v = (sp_{\chi^{-1}}(X_\alpha - 1)J_s) sp_\chi v.$$

*Proof.* (i) According to Lemma **1.8** and our choices, the representation $i_{P\cap M_\alpha}^{M_\alpha}\sigma$ is irreducible, and the restriction of $J_s$ to $i_{P\cap K}^K E$ has a simple pole at $Y_\alpha = 1$. The operator $sp_1(X_\alpha - 1)J_s$ is defined from an intertwining of this representation with itself. So it must be a scalar. We have

$$(sp_1(X_\alpha - 1)J_s)^2 = sp_1((X_\alpha - 1)J_s(X_\alpha - 1)J_s) = sp_1((X_\alpha - 1)(X_\alpha^{-1} - 1)J_s^2).$$

As $(J_s)^2 = c_s''(\mu^{M_\alpha})^{-1}$, we indeed find the expression in the statement, after having replaced $\mu^{M_\alpha}$ by the expression given in the proposition **1.6**.

Suppose for now $\chi \in \mathfrak{X}^{\mathrm{nr}}(M)$ such that $X_\alpha(\chi) = 1$. There thus exist $\eta \in \mathrm{Stab}(\mathcal{O})$ and $\chi_1 \in \mathfrak{X}^{\mathrm{nr}}(M_\alpha)$ such that $\chi = \eta\chi_1$. We have ${}^s\chi_1 = \chi_1$ and $\chi_1(m_s) = 1$. Note moreover that as a vector space, $i_{P\cap K}^K E_{B_{\mathcal{O}}} = \bigoplus_{m\in M/M^\sigma} i_{P\cap K}^K \sigma(m)E_1 \otimes b_m$. Let $m \in M$ and $v_m \in i_{P\cap K}^K \sigma(m)E_1$. We have, as rational functions of $\chi' \in \mathfrak{X}^{\mathrm{nr}}(M)$,

$$\begin{aligned}
&sp_{\chi'\chi}(X_\alpha - 1)J_s v_m \otimes b_m\\
&= (X_\alpha(\chi\chi') - 1)(\chi\chi')(m_s)\rho_{\sigma,s}\lambda(s)J_{sP|P}(\sigma \otimes {}^s(\chi'\chi))^s(\chi'\chi)(m)v_m\\
&= \chi_1(m)(X_\alpha(\chi') - 1)(\eta\chi')(m_s)\rho_{\sigma,s}\lambda(s)J_{sP|P}(\sigma \otimes {}^s(\chi'\eta))^s\chi'(m)^s\eta(m)v_m\\
&= \chi_1(m)(X_\alpha(\chi') - 1)(\eta\chi')(m_s)\rho_{\sigma,s}\lambda(s)J_{sP|P}(\sigma \otimes {}^s(\chi'\eta))^s\chi'(m)\phi_{\sigma,s\eta}v_m\\
&= \chi_1(m)(X_\alpha(\chi') - 1)(\eta\chi')(m_s)\rho_{\sigma,s}\phi_{\sigma,s\eta}\lambda(s)J_{sP|P}(\sigma \otimes {}^s\chi')^s\chi'(m)v_m.
\end{aligned}$$

Specializing at $\chi' = 1$, we find

$$sp_\chi(X_\alpha - 1)J_s v_m \otimes b_m = \chi_1(m)\eta(m_s)\rho_{\sigma,s}\phi_{\sigma,s\eta}\rho_{\sigma,s}^{-1} sp_1(X_\alpha - 1)J_s v_m \otimes b_m.$$

We have just seen that $sp_1(X_\alpha - 1)J_s v_m \otimes b_m$ is an element of $i_{P\cap K}^K\sigma(m)E_1$. Let $e_1 \in E_1$. Note that $\rho_{\sigma,s}^{-1}e_1$ is an element of $(s\sigma)(m_s^{-1})E_1$. Therefore

$$\begin{aligned}
\rho_{\sigma,s}\phi_{\sigma,s\eta}\rho_{\sigma,s}^{-1}\sigma(m)e_1 &= \rho_{\sigma,s}\phi_{\sigma,s\eta}(s\sigma)(m)\rho_{\sigma,s}^{-1}e_1\\
&= \rho_{\sigma,s}\eta(m_s^{-1})\eta(m)(s\sigma)(m)\rho_{\sigma,s}^{-1}e_1\\
&= \eta(m_s^{-1})\eta(m)\sigma(m)e_1,
\end{aligned}$$

i.e. $(\rho_{\sigma,s}\phi_{\sigma,s\eta}\rho_{\sigma,s}^{-1})_{|i_{P\cap K}^K\sigma(m)E_1} = \eta(m_s^{-1})\eta(m)$. The last assertion of (i) of Lemma follows.

(ii) Observe first that ${}^s\chi_{-1}\chi_{-1}^{-1} \in \mathrm{Stab}(\mathcal{O})$: as $\mu^{M_\alpha}$ cancels out at $X_\alpha = -1$, we find by results of Harish-Chandra that $\sigma \otimes \chi_{-1} \simeq s(\sigma \otimes \chi_{-1}) \simeq s\sigma \otimes {}^s\chi_{-1} \simeq \sigma \otimes {}^s\chi_{-1}$, hence $\sigma \simeq \sigma \otimes {}^s\chi_{-1}\chi_{-1}^{-1}$. Set $\eta_{-1} = {}^s\chi_{-1}\chi_{-1}^{-1}$. Let $v_m \in i_{P\cap K}^K\sigma(m)E_1$.

We thus have, as rational function of $\chi' \in \mathfrak{X}^{\mathrm{nr}}(M)$,

$$\begin{aligned}
&sp_{\chi'\chi_{-1}}(X_\alpha+1)J_s v_m \otimes b_m \\
&= (X_\alpha(\chi'\chi_{-1})+1))(\chi'\chi_{-1})(m_s)\rho_{\sigma,s}\lambda(s)J_{sP|P}(\sigma\otimes {}^s(\chi_{-1}\chi')^s(\chi_{-1}\chi')(m)v_m \\
&= (X_\alpha(\chi'\chi_{-1})+1))(\chi'\chi_{-1})(m_s)\rho_{\sigma,s}\lambda(s)J_{sP|P}(\sigma\otimes\eta_{-1}\chi_{-1}{}^s\chi')(\chi_{-1}{}^s\chi')(m)\times\phi_{\sigma,\eta_{-1}}v_m \\
&= (X_\alpha(\chi'\chi_{-1})+1))(\chi'\chi_{-1})(m_s)\rho_{\sigma,s}\phi_{\sigma,\eta_{-1}}\lambda(s)J_{sP|P}(\sigma\otimes\chi_{-1}{}^s\chi')(\chi_{-1}{}^s\chi')(m)v_m \\
&= \chi_{-1}(m)\chi_{-1}(m_s)\rho_{\sigma,s}\phi_{\sigma,\eta_{-1}}\rho^{-1}_{\sigma\otimes\chi_{-1},s}(-X_\alpha(\chi')+1)\chi'(m_s)\rho_{\sigma\otimes\chi_{-1},s}\lambda(s)\times J_{sP|P}(\sigma\otimes\chi_{-1}{}^s\chi')\chi'(m)v_m
\end{aligned}$$

Denote by $J_{s,\chi_{-1}}$ the operator in $\mathrm{End}_G((E_{\chi_{-1}})_{B_\mathcal{O}})_{K(B_\mathcal{O})}$ defined like $J_s$, by replacing $\sigma$ by $\sigma\otimes\chi_{-1}$. Thus, according to (i), we have for $\chi\in\mathfrak{X}^{\mathrm{nr}}(M)$ with $X_\alpha(\chi)=1$,

$$\begin{aligned}
&sp_{\chi_{-1}\chi}((X_\alpha+1)J_s)v_m\otimes b_m \\
&= \chi_{-1}(m)\chi_{-1}(m_s)\rho_{\sigma,s}\phi_{\sigma,\eta_{-1}}\rho^{-1}_{\sigma\otimes\chi_{-1},s}sp_\chi(-X_\alpha+1)J_{s,\chi_{-1}}v_m\otimes b_m \\
&= \chi_{-1}(m_s)\rho_{\sigma,s}\phi_{\sigma,\eta_{-1}}\rho^{-1}_{\sigma\otimes\chi_{-1},s}(sp_1(-X_\alpha+1)J_{s,\chi_{-1}})sp_{\chi_{-1}\chi}v_m\otimes b_m.
\end{aligned}$$

The operator $\rho_{\sigma,s}\phi_{\sigma,\eta_{-1}}\rho^{-1}_{\sigma\otimes\chi_{-1},s}$ comes via functoriality from an automorphism of an irreducible representation in $E$, since $s(\sigma\otimes\chi_{-1}) = s\sigma\otimes\chi_{-1}\eta_{-1}$. This is thus a scalar. By setting above $\chi=1$ and by using the part (i) of lemma with respect to the operator $J_{s,\chi_{-1}}$, we find that the operator $sp_{\chi_{-1}}(X_\alpha+1)J_s$ is scalar on $i^K_{P\cap K}E$. The equality above thus also implies the last assertion of part (ii) of the lemma.

Concerning the second assertion, it suffices to note that

$$(sp_{\chi_{-1}}(X_\alpha+1)J_s)^2 = sp_{\chi_{-1}}(X_\alpha+1)(X_\alpha^{-1}+1)J_s^2,$$

and we conclude as for (i). □

5.2

**Definition** Let $\alpha\in\Delta_\mathcal{O}$, $s=s_\alpha$. Fix a square root $c_s^{1/2}$ of $c_s$. According to lemma **5.1**, there exists $\epsilon_1,\epsilon_{-1}\in\{\pm1\}$ such that

$$sp_1(X_\alpha-1)J_s = \epsilon_1 c_s^{1/2}\frac{(1-q^{-a_s})(1+q^{-b_s})}{2}sp_1$$

(et) $$sp_{\chi_{-1}}(X_\alpha+1)J_s = \epsilon_{-1}c_s^{1/2}\frac{(1+q^{-a_s})(1-q^{-b_s})}{2}sp_{\chi_{-1}}.$$

Set

$$R_s = \begin{cases} -\epsilon_1 q^{a_s+b_s}c_s^{-1/2}J_s, & \text{if } \epsilon_1 b_s\neq\epsilon_{-1}b_s; \\ -\epsilon_1 q^{a_s+b_s}c_s^{-1/2}X_\alpha J_s, & \text{if } \epsilon_1 b_s=\epsilon_{-1}b_s \end{cases}.$$

(et) $$T_s = R_s + (q^{a_s+b_s}-1)\frac{X_\alpha\Big(X_\alpha - \frac{q^{b_s}-q^{a_s}}{q^{a_s+b_s}-1}\Big)}{X_\alpha^2-1}.$$

5.3

**Lemma** *An element* $A \in \mathrm{Hom}_G(i_P^G E_{B_{\mathcal{O}}}, i_P^G E_{K(B_{\mathcal{O}})})$ *is in* $\mathrm{End}_G(i_P^G E_{B_{\mathcal{O}}})$, *if and only if, for every* $v \in i_P^G E_{B_{\mathcal{O}}}$ *and for every* $v^\vee \in i_{P\cap K}^K E^\vee$, *the map* $\mathfrak{X}^{\mathrm{nr}}(M) \to \mathbb{C}$, $\chi \mapsto \langle sp_\chi A v, v^\vee\rangle$, *is regular.*

*Proof.* The condition is necessary, since, if $A(v) \in i_P^G E_{B_{\mathcal{O}}}$, thus $sp_\chi(A(v)) \in i_P^G E_\chi$ for every $\chi \in \mathfrak{X}^{\mathrm{nr}}(M)$.

The condition is sufficient: let $v \in i_P^G E_{B_{\mathcal{O}}}$. By the inclusion $i_P^G E_{B_{\mathcal{O}}} \subseteq i_{P\cap K}^K E \otimes K(B)$, we can write $Av = \sum_{i\in I} v_i \otimes b_i$ with $b_i \in K(B)$ and the $v_i$ are $\mathbb{C}$-linearly-independent in $i_{P\cap K}^K E$.

Let $v_i^\vee \in i_{P\cap K}^K E^\vee$, $i \in I$, the dual elements of $v_i$. Thus $\langle sp_\chi(A(v)), v_i^\vee\rangle = b_i$. It thus follows that $b_i \in B$ for every $i \in I$, i.e. $Av \in i_P^G E_B \cap (i_P^G E_{B_{\mathcal{O}}})_{K(B_{\mathcal{O}})} = i_P^G E_{B_{\mathcal{O}}}$. □

5.4

**Proposition** *The operators* $T_{s_\alpha}$, $\alpha \in \Delta_{\mathcal{O}}$, *defined in* ***5.2***, *belong to* $\mathrm{End}_G(i_P^G E_{B_{\mathcal{O}}})$.

*Proof.* Set $s = s_\alpha$. It is necessary to show that $T_s$ is regular at every $\chi \in \mathfrak{X}^{\mathrm{nr}}(M)$. Write $T_s = p_s J_s + r_s$ with $p_s \in B^\times$ and $r_s \in K(B_{\mathcal{O}})$. Let $v \in i_{P\cap K}^K E_{\mathcal{O}}$ and $v^\vee \in E^\vee$. The map $\chi \mapsto \langle sp_\chi T_s v, v^\vee\rangle$ is regular at $\chi$, unless maybe if $X_\alpha(\chi) = 1$ or if $X_\alpha(\chi) = -1$ and $b_s \neq 0$.

First note that

$$\mathrm{Res}_{X_\alpha = \pm 1}\, r_s = \frac{\pm 1}{2}(q^{a_s+b_s} - 1 - (\pm 1)(q^{b_s} - q^{a_s})).$$

The regularity at $\chi \in \mathfrak{X}^{\mathrm{nr}}(M)$ such that $X_\alpha(\chi) = 1$ thus follows from the equality $sp_\chi(X_\alpha - 1)p_s J_s = -q^{a_s+b_s}\frac{(1-q^{-a_s})(1+q^{-b_s})}{2} sp_\chi$, and that at $\chi \in \mathfrak{X}^{\mathrm{nr}}(M)$ such that $X_\alpha(\chi) = -1$ follows from the equality $sp_\chi(X_\alpha+1)p_s J_s = q^{a_s+b_s}\frac{(1+q^{-a_s})(1-q^{-b_s})}{2} sp_\chi$, if $b_s \neq 0$. □

5.5

**Proposition** *Let* $\alpha \in \Delta_{\mathcal{O}}$, $s = s_\alpha$. *Thus* $(T_s + 1)(T_s - q^{a_s+b_s}) = 0$.

*Proof.* Let $T_s$ be an operator of the form $p_s J_s + r_s$ with $p_s \in \mathbb{C}^\times$ and $r_s \in K(B)$. Thus

$$\begin{aligned}(T_s + 1)(T_s - q^{a_s+b_s}) =& T_s^2 + (1 - q^{a_s+b_s})T_s - q^{a_s+b_s} \\ =& (r_s + {}^s r_s + 1 - q^{a_s+b_s})p_s J_s \\ & + p_s^2 J_s^2 + r_s^2 + (1 - q^{a_s+b_s})r_s - q^{a_s+b_s}.\end{aligned}$$

As 1 and $J_s$ are $K(B)$-linearly independent, the equation $(T_s+1)(T_s-q^{a_s+b_s}) = 0$ equals to

$$r_s + {}^s r_s + 1 - q^{a_s+b_s} = 0 \quad \text{and} \quad p_s^2 J_s^2 + r_s^2 + (1-q^{a_s+b_s})r_s - q^{a_s+b_s} = 0.$$

For $T_s$ defined as above in the case $\epsilon_1 b_s \neq \epsilon_{-1} b_s$, the verification of these identities is done by an elementary calculation. The case $\epsilon_1 b_s = \epsilon_{-1} b_s$ as a result, since $(X_\alpha J_s)(X_\alpha J_s) = J_s^2$. □

5.6

**Definition** Fix for every $w \in W_{\mathcal{O}}$ a decomposition into simple symmetries $w = s_1 \cdots s_\ell$ and define $T_w = T_{s_1} \cdots T_{s_\ell}$. (In particular, $T_1 = \mathrm{Id}$.)

5.7

**Proposition** *Let $w, w' \in W_{\mathcal{O}}$.*

(i) *If $\ell_{\mathcal{O}}(ww') = \ell_{\mathcal{O}}(w) + \ell_{\mathcal{O}}(w')$, then $T_w T_{w'}$ is of teh form $b_{w,w'} J_{ww'} + \sum_{w''} f_{w,w',w''} J_{w''}$, with $b_{w,w'} \in B_{\mathcal{O}}^\times$ and $f_{w,w',w''} \in K(B_{\mathcal{O}})$, $f_{w,w',w''} = 0$ when $\ell_{\mathcal{O}}(w'') \geq \ell_{\mathcal{O}}(ww')$.*

(ii) *If $\ell_{\mathcal{O}}(ww') < \ell_{\mathcal{O}}(w) + \ell_{\mathcal{O}}(w')$, then $T_w T_{w'}$ is of the form $\sum_{w''} f_{w,w',w''} J_{w''}$, with $f_{w,w',w''} \in K(B_{\mathcal{O}})$, $f_{w,w',w''} = 0$ when $\ell_{\mathcal{O}}(w'') \geq \ell_{\mathcal{O}}(w) + \ell_{\mathcal{O}}(w')$.*

*Proof.* This follows directly from the definition of $T_w$ and from **4.7** (i), (ii). □

5.8

**Lemma** *Let $w' \in W_{\mathcal{O}}$ and denote by $w_0$ the longest element in $W_{\mathcal{O}}$. Suppose for every $w \in W_{\mathcal{O}}$ we're given an operator $T_{w,w'}$ in $\mathrm{End}_G(i_P^G E_{B_{\mathcal{O}}})$ of the form $\sum_{w''} f_{w,w''} J_{w''}$, with $f_{w,w''} \in K(B_{\mathcal{O}})$, $f_{w,w_0} = 0$ if $w \neq w'$ and $f_{w',w_0} \in B_{\mathcal{O}}^\times$. Let $\chi \in \mathfrak{X}^{\mathrm{nr}}(M)$ and let $c_{r,w}$, $r \in R(\mathcal{O})$, $w \in W_{\mathcal{O}}$ be complex numbers such that $\sum_{r,w} c_{r,w} sp_\chi J_r T_{w,w'} = 0$. Thus $c_{1,w'} = 0$.*

*Proof.* Fix a subgroup $H$ of $G$ satisfying the indicated properties in [Hei01, 5.] relative to $P$ and $\mathcal{O}$. Fix an element $v \neq 0$ of $(i_{P\cap K}^K E)^H$ with support contained in $(P \cap K)H$ and with values in $E_1$. It is thus identified with the element of $i_P^G E_{B_{\mathcal{O}}}$ defined for $m \in M$, $u \in U$ and $k \in K$ by $muk \mapsto \delta_P^{1/2} \sigma_B(m) v(k)$. Let $v^\vee$ be an element of $(i_{\overline{P}\cap K}^K E^\vee)^H$ with support contained in $(\overline{P} \cap K)H$. We suppose that $v^\vee$ was chosen such that $\langle v(1), v^\vee(1)\rangle \neq 0$. The rational operator $J_{P|\overline{w_0 P}} \lambda(w_0)$ is regular on $\mathcal{O}^\vee$ by **1.10** (where $w_0$ designates the element of maximal length of $W_{\mathcal{O}}$). Set $A = J_{P|\overline{w_0 P}}(\sigma \otimes \cdot)\lambda(w_0)$.

We deduce from [Hei01, 5.3] (where one has to reverse the roles of $v$ and $v^\vee$) that, for $w \in W(M,\mathcal{O})$, $\langle J_{P|wP}(\sigma\otimes\chi')\lambda(w)\rho_w v, A(\chi')\rho_{w_0}^\vee v^\vee\rangle$ is a constant function in $\chi'$ which can only be $\neq 0$ if $w = w_0$. If $w = w_0$, then its value is of the form $c_w\langle v(1), v^\vee(1)\rangle$, where $c_w$ is a complex number $\neq 0$ which does not depend on $v$ or $v^\vee$.

Suppose $\sum_{r,w} c_{r,w} sp_\chi J_r T_{w,w'} = 0$. Therefore,

$$\begin{aligned} 0 &= \left\langle \sum_{r,w} c_{r,w} J_r T_{w,w'} v, A\rho_{w_0}^\vee v^\vee \right\rangle \\ &= \sum_{r,w}\sum_{w''} \langle c_{r,w} J_r f_{w,w''} J_{w''} v, A\rho_{w_0}^\vee v^\vee\rangle \\ &= \sum_{r,w}\sum_{w''} c_{r,w} f_{w,w''}\langle J_r J_{w''} v, A\rho_{w_0}^\vee v^\vee\rangle \end{aligned}$$

Observe for now that, according to lemma **3.2** and the proof of **3.3**, for $w \in W_\mathcal{O}$, $r \in R(\mathcal{O})$ and $\chi' \in \mathfrak{X}^{\mathrm{nr}}(M)$ a regular point, and considering the definition of $J_r J_w$, $sp_{\chi'} J_r J_w v$ is the product of a nonzero complex number with

$$\begin{aligned} &\rho_{\sigma,rw}\lambda(rw) J_{w^{-1}r^{-1}P|P}(\sigma\otimes(w^{-1}r^{-1}\chi')) sp_{w^{-1}r^{-1}\chi'} v \\ &= \rho_{\sigma,rw} J_{P|rwP}(rw\sigma\otimes\chi')\lambda(rw) v \\ &= J_{P|rwP}(\sigma\otimes\chi')\lambda(rw)\rho_{\sigma,rw} v, \end{aligned}$$

the commutation properties for $\rho_{\sigma,rw}$ compared to $J_{P|rwP}(\cdot)\lambda(rw)$ resulting from the fact that this one is defined by functoriality from an isomorphism between the representations in $E$ (cf. **2.5**).

According to the above, we see that the product $\langle J_{P|rwP}(\sigma\otimes\chi')\lambda(rw)\rho_{\sigma,rw} v, A(\chi')\rho_{w_0}^\vee v^\vee\rangle$ is zero, unless $r = 1$ and $w = w_0$, and thus its value is by choice of nonzero $v^\vee$ equal to $c_{w_0}\langle(v(1), v^\vee(1))\rangle =: c$. It follows that

$$\begin{aligned} 0 &= \left\langle \sum_{r,w} sp_\chi c_{r,w} J_r T_{w,w'} v, A(\chi)\rho_{w_0}^\vee v^\vee \right\rangle \\ &= c_{1,w'} f_{w',w_0}(\chi) c. \end{aligned}$$

As $f_{w',w_0} \in B^\times$, we have $f_{w',w_0}(\chi) \neq 0$, and it follows that $c_{1,w'} c = 0$, hence $c_{1,w'} = 0$. □

5.9

**Proposition** *For every $\chi \in \mathfrak{X}^{\mathrm{nr}}(M)$, the operators $sp_\chi J_r T_w$, $r \in R(\mathcal{O})$, $w \in W_\mathcal{O}$, are linearly independent.*

*Proof.* Let $\chi \in \mathfrak{X}^{\mathrm{nr}}(M)$. We will first prove by decreasing induction on $\ell$, $0 \leq \ell \leq \ell_\mathcal{O}(w_0)$, that, when $c_{r,w}$, $r \in R(\mathcal{O})$, $w \in W_\mathcal{O}$, are complex numbers such that $\sum_{r,w} c_{r,w} sp_\chi J_r T_w = 0$, thus $c_{1,w} = 0$ for every $w$ of top length or equal to $\ell$.

It follows from the previous lemma and from **5.7** (taking $w' = 1$) that this affirmation is true for $\ell = \ell(w_0)$. So $0 \leq \ell < \ell(w_0)$ and suppose it's true for $\ell+1$. Let $c_{r,w}$ be complex numbers such that $\sum_{r,w} c_{r,w} sp_\chi J_r T_w = 0$. We have thus $c_{1,w} = 0$ when $\ell(w) > \ell$. Let $w_1 \in W_{\mathcal{O}}$ such that $\ell(w_1) = \ell$. Set $w' = w_1^{-1} w_0$ and $A := \sum_{r,w} c_{r,w} J_r T_w$. By hypothesis, $sp_\chi A = 0$ and thus $sp_\chi A T_{w'} = 0$.

On the other hand,

$$AT_{w'} = \sum_{r,w} c_{r,w} J_r (T_w T_{w'}).$$

We will deduce from Lemma **5.8** that $c_{1,w_1} = 0$. If $\ell(w) > \ell$, then $c_{1,w} = 0$ by recurrence hypothesis. If $\ell(w) = \ell$, thus $\ell(ww') \leq \ell(w) + \ell(w') = \ell(w_0)$, since by [Blo12, VI, paragraph 1, Cor. 3], $\ell(w') = \ell(w_0) - \ell$. By the proposition **5.7**, $T_w T_{w'}$ is, for $\ell(w) \leq \ell$ and $w \neq w_1$, $K(B_{\mathcal{O}})$-linear combination of operators $J_{w''}$, $\ell(w'') < \ell(w_0)$, while $T_{w_1} T_{w'}$ is of the form $f_{w_1,w_0} J_{w_0} + \sum_{w'' \neq w_0} f_{w_1,w''} J_{w''}$ with $f_{w_1,w_0} \in B_{\mathcal{O}}^\times$ and $f_{w_1,w''} \in K(B_{\mathcal{O}})$. The hypothesis of Lemma **5.8** is thus indeed verified, and, consequently, $c_{1,w_1} = 0$.

Let us now deduce that $c_{r',w} = 0$ for every $r' \in R(\mathcal{O})$ and $w \in W_{\mathcal{O}}$. Let $r' \in R(\mathcal{O})$. Remark that $T_w J_{r'^{-1}} = J_{r'^{-1}} T_{r'wr'^{-1}}$ by **3.5**.

As $0 = sp_\chi A J_{r'^{-1}}$ and, on the other hand,

$$AJ_{r'^{-1}} = \sum_{r,w} c_{r,w} J_r J_{r'^{-1}} T_{r'wr'^{-1}},$$

we find, by using **4.7** (iii),

$$0 = sp_\chi A J_{r'^{-1}} = \sum_{r,w} c_{r,r'^{-1}wr'} sp_\chi J_{rr'^{-1}} T_w.$$

By what was shown in the first part, it follows that $c_{r',w} = 0$ for every $w \in W_{\mathcal{O}}$. □

5.10

**Theorem**

$$\mathrm{End}_G(i_P^G E_{B_{\mathcal{O}}}) = \bigoplus_{r \in R(\mathcal{O}), w \in W_{\mathcal{O}}} B_{\mathcal{O}} J_r T_w.$$

*Proof.* It is clear by definition and **5.4** that the space on the right is included in the space on the left.

On the other hand, by theorems **4.9** and **5.7**, the space on the left is included in $\bigoplus_{r,w} K(B_{\mathcal{O}})J_rT_w$. It suffices thus to prove that a linear combination of the form $A := \sum_{r,w} b^{-1}b_{r,w}J_rT_w$ with $b_{r,w}$ and $b$ in $B_{\mathcal{O}}$ is in $\mathrm{End}_G(i_P^G E_{B_{\mathcal{O}}})$, if and only if $b$ divides all the $b_{r,w}$ in $B_{\mathcal{O}}$.

We are brought back to the case where $b$ and $b_{r,w}$ are prime to each other in $B_{\mathcal{O}}$. It is then necessary to prove that $b$ is a unit. Suppose by contradiction that $b$ is not a unit. The ring $B_{\mathcal{O}}$ being factorial, it suffices to consider the case where $b$ is irreducible. Thus there exists at least an element $\chi \in \mathfrak{X}^{\mathrm{nr}}(M)$ such that $b(\chi) = 0$, thus $b_{r,w}(\chi) \neq 0$ for at least a pair $(r, w)$. As $0 = sp_\chi bA = \sum_{r,w} b_{r,w}(\chi)sp_\chi J_rT_w$, this contradicts the proposition **5.9** which says that the operators $sp_\chi J_rT_w$ are $\mathbb{C}$-linearly independent for every $\chi \in \mathfrak{X}^{\mathrm{nr}}(M)$. □

5.11

**Corollary** *The center* $\mathcal{Z}_{\mathcal{O}}$ *of* $\mathrm{End}_G(i_P^G E_{B_{\mathcal{O}}})$ *is formed of* $W(M, \mathcal{O})$*-invariant elements in* $B_{\mathcal{O}}$*.*

*Proof.* Indeed, suppose $\sum_{r,w} b_{r,w}J_rJ_w$ in the center. Thus, by the commutation relations (cf.**4.8**), we have for every $r \in R(\mathcal{O})$, $w \in W_{\mathcal{O}}$, $b \in B_{\mathcal{O}}$, $bb_{r,w} = {}^{rw}bb_{r,w}$. Consequently, by the linear independent of operators $J_rJ_w$ (cf.**4.8**), $b_{r,w} = 0$ if $rw \neq 1$. The center of $\mathrm{End}_G(i_P^G E_{B_{\mathcal{O}}})$ is thus contained in $B_{\mathcal{O}}$, and it is thus clear that this is $B_{\mathcal{O}}^{W(M,\mathcal{O})}$. □

5.12

Denote by $K(\mathcal{Z}_{\mathcal{O}})$ the fraction field of $\mathcal{Z}_{\mathcal{O}}$ and set $\mathrm{End}_G(i_P^G E_{B_{\mathcal{O}}})_{K(\mathcal{Z}_{\mathcal{O}})} = \mathrm{End}_G(i_P^G E_{B_{\mathcal{O}}}) \otimes_{\mathcal{Z}_{\mathcal{O}}} K(\mathcal{Z}_{\mathcal{O}})$.

**Corollary** *The algebra* $\mathrm{End}_G(i_P^G E_{B_{\mathcal{O}}})_{K(\mathcal{Z}_{\mathcal{O}})}$ *is a* $K(B_{\mathcal{O}})$*-module which is canonically isomorphic to* $\mathrm{Hom}_G(i_P^G E_{B_{\mathcal{O}}}, i_P^G E_{K(B_{\mathcal{O}})})$*.*

*Proof.* First note that $B_{\mathcal{O}} \otimes_{\mathcal{Z}_{\mathcal{O}}} K(\mathcal{Z}_{\mathcal{O}})$ is canonically isomorphic to $K(B_{\mathcal{O}})$ (cf.[Lus89]). It follows by **5.10** that $\mathrm{End}_G(i_P^G E_{B_{\mathcal{O}}})_{K(\mathcal{Z}_{\mathcal{O}})}$ is isomorphic to $\bigoplus_{r,w} K(B_{\mathcal{O}})J_rT_w$, hence the corollary by **4.9** and **5.7**. □

**6**

Recall that we have defined in **1.5** for every $\alpha \in \Delta_{\sigma,\mu}$ an element $h_\alpha$ in $M \cap M_\alpha^1$ and a real number $t_\alpha$. Set $\widetilde{\alpha} := H_M(h_\alpha^{t_\alpha})$. Denote by $\Lambda_{\mathcal{O}}$ the free $\mathbb{Z}$-module in $a_M$ equal to the image of $M^\sigma/M^1$ by $H_M$, by $\Lambda_{\mathcal{O}}^\vee$ the free $\mathbb{Z}$-module included in $a_M^*$ which is by $\langle\cdot,\cdot\rangle$ in perfect duality with $\Lambda_{\mathcal{O}}$, by $\Sigma_{\mathcal{O}}$ the set of $\widetilde{\alpha}$, $\alpha \in \Sigma_{\mathcal{O},\mu}$, and, finally, by $\Sigma_{\mathcal{O}}^{\vee}$ the set of multiples $\alpha^*$ of $\alpha$, $\alpha \in \Sigma_{\mathcal{O},\mu}$, such that $\langle \alpha^*, \widetilde{\alpha} \rangle = 2$.

6.1

**Proposition** *The quadruple $(\Lambda_{\mathcal{O}}, \Sigma_{\mathcal{O}}, \Lambda_{\mathcal{O}}^{\vee}, \Sigma_{\mathcal{O}}^{\vee})$ is a given radical. The root system underlying it is reduced. The Weyl group of $\Sigma_{\mathcal{O}}$ is canonically isomrophic to $W_{\mathcal{O}}$, and the set $\Delta_{\mathcal{O}} = \{H_M(h_\alpha^{t_\alpha}), \alpha \in \Delta_{\mathcal{O},\mu}\}$ form a basis of $\Sigma_{\mathcal{O}}$.*

*Moreover, if $\Sigma'_\mu$ is an irreducible component of $\Sigma_{\mathcal{O},\mu}$, and if $\Sigma'_\mu$ is not of type $C_n$, then the set of $\widetilde{\alpha}$, $\alpha \in \Sigma'_\mu$, is an irreducible component of $\Sigma_{\mathcal{O}}$ which is of the same type as $\Sigma'$. In the contrary case, the corresponding irreducible component of $\Sigma_{\mathcal{O}}$ is of type $B_n$.*

First a lemma:

6.2

**Lemma** *Let $\alpha \in \Sigma_{\mathcal{O},\mu}$. Denote by $m_\alpha$ the largest number $> 0$, such that $\chi_{\frac{2\pi i}{m_\alpha \log q}\alpha} \in \mathfrak{X}^{\mathrm{nr}}(M_\alpha)$. Thus $H_M(h_\alpha) = \frac{m_\alpha}{2}\alpha^\vee$.*

*Proof.* Write $H_M(h_\alpha) = \frac{m}{2}\alpha^\vee$ with $m > 0$. Thus

$$1 = \chi_{\frac{2\pi i}{m_\alpha \log q}\alpha}(h_\alpha) = q^{-\langle \frac{2\pi i}{m_\alpha \log q}\alpha, H_M(h_\alpha)\rangle} = q^{-\frac{2\pi i}{\log q}\frac{m}{m_\alpha}}.$$

It follows that $m \in m_\alpha\mathbb{Z}$. Thus thanks to the maximality of $m_\alpha$, it is necessary that $m_\alpha = m$. □

*Proof.* (of the proposition) Fix a component $\Sigma_{\mathcal{O},\mu,i}$ of $\Sigma_{\mathcal{O},\mu}$, and denote by $\{\alpha_1, \cdots, \alpha_d\}$ its basis. Recall that $d = d_i$, unless $\Sigma_{\mathcal{O},\mu,i}$ is of type $A_d$. Thus $d = d_i - 1$. An element of $\mathfrak{X}^{\mathrm{nr}}(M)$ is of the form $|\det_{m_1}|^{s_{1,1}} \otimes \cdots \otimes |\det_{m_1}|^{s_{1,d_1}} \otimes |\det_{m_2}|^{s_{2,1}} \otimes \cdots \otimes |\det_{m_2}|^{s_{2,d_2}} \otimes \cdots \otimes |\det_{m_r}|^{s_{r,1}} \otimes \cdots \otimes |\det_{m_r}|^{s_{r,d_r}} \otimes 1$, where $s_{i,j}$ are complex numbers. For $j = 1, \cdots, d_i - 1$, $\chi_{\alpha_j}$ is given by $s_{i,j} = 1/m_i$, $s_{i,j+1} = -1/m_i$, the other exponents $s_{i',j'}$ being zero. If $d = d_i$, thus $\chi_{\alpha_d}$ is given by $s_{i,d} = 1/m_i$ and $s_{i,d-1} = 0$ (resp. $s_{i,d-1} = 1/m_i$), if $\Sigma_{\mathcal{O},\mu,i}$ is of type $B_d$ (resp. $D_d$), and by $s_{i,d} = 2/m_i$ if $\Sigma_{\mathcal{O},\mu,i}$ is of type $C_d$, the other exponents being zero.

Fix a uniformizer $\widetilde{\omega}$ of $F$. We can choose $h_{\alpha_j} = \mathrm{diag}(1, 1, \cdots, 1, \widetilde{\omega}, \widetilde{\omega}^{-1}, 1, \cdots, 1)$ for $j = 1, \cdots, d_i - 1$, and, if $d = d_i$, $h_{\alpha_d} = \mathrm{diag}(1, \cdots, 1, \widetilde{\omega}, 1, \cdots, 1)$ if $\Sigma_{\mathcal{O},\mu,i}$ is of type $B_d$ or $C_d$, and $h_{\alpha_d} = \mathrm{diag}(1, 1, \cdots, 1, \widetilde{\omega}, \widetilde{\omega}, 1, \cdots, 1)$ if $\Sigma_{\mathcal{O},\mu,i}$ is of type $D_d$, $\widetilde{\omega}$ each time being in the last place on the diagonal of the $j$-th (or $d$-th) copy of $\mathrm{GL}_{m_i}$.

Set $\epsilon = 2$, if $\Sigma_{\mathcal{O},\mu,i}$ is of type $B_d$, $\epsilon = 1$ if not. Thus, we can deduce $m_{\alpha_j} = 2m_i^{-1}$ for $j = 1, \cdots, d_i - 1$, and, if $d = d_i$, $m_{\alpha_d} = \frac{2}{\epsilon}m_i^{-1}$. If $t_i$ denotes the order of stabilizer of $\sigma_i$, we have moreover $t_{\alpha_j} = t_i$ for $j = 1, \cdots, d$. As a result, thanks to lemma **6.2**,

$$\Delta_{\mathcal{O},i} = \left\{ \frac{t_i}{m_i}\alpha_1^\vee, \cdots, \frac{t_i}{m_i}\alpha_{d_i-1}^\vee, \frac{t_i}{\epsilon m_i}\alpha_d^\vee \right\}.$$

These are elements of $\Lambda_{\mathcal{O}} = H_M(M^\sigma/M^1)$, since by construction these elements are of the form $H_M(h_\alpha^{t_\alpha})$, $\alpha \in \Delta_{\mathcal{O},\mu}$, thus $h_\alpha^{t_\alpha} \in M^\sigma$ according to **1.13**, by remarking that the group $\mathrm{Stab}(\mathcal{O})$ acts, according definition **1.5**, trivially on $h_\alpha^{t_\alpha}$.

We set $\Delta_{\mathcal{O},i}^\vee = \left\{ \frac{m_i}{t_i}\alpha_1, \cdots, \frac{m_i}{t_i}\alpha_{d_i-1}, \frac{\epsilon m_i}{t_i}\alpha_{d_i} \right\}$. We have $\Delta_{\mathcal{O},i}^\vee \subseteq \Lambda_{\mathcal{O}}^\vee$: first observe that the groups $\mathrm{GL}_{m_i}/(\mathrm{GL}_{m_i})^1$ are cyclic. The group $(\mathrm{GL}_{m_i})^{\sigma_i}/(\mathrm{GL}_{m_i})^1$ is a subgroup of index $t_i$ (because $(\mathrm{GL}_{m_i}/(\mathrm{GL}_{m_i})^1)^{t_i} = (\mathrm{GL}_{m_i})^\sigma/(\mathrm{GL}_{m_i})^1$), and the group $M^\sigma/M^1$ is a product of such groups. It is clear that the rational characters $m_i\alpha_j$, considered as elements of $a_M^*$, send the image of $M$ by $H_M$ to $\mathbb{Z}$. As elements of $\Delta_{\mathcal{O},i}^\vee$ are trivial on the factors $\mathrm{GL}_{m_j}$, $j \neq i$, it follows that $\frac{m_i}{t_i}\alpha_j$ sends $M^\sigma/M^1$ to $\mathbb{Z}$, i.e. these are elements of $\Lambda_{\mathcal{O}}^\vee$.

We see that the root system obtained is of the same type as $\Sigma_{\mathcal{O},\mu,i}$ except in the case where it is of type $C_d$. Thus, we find a system equivalent to the system dual to $\Sigma_{\mathcal{O},\mu,i}$. It is thus of type $B_d$. □

## 7

We will now make the link with Hecke algebra with parameters.

### 7.1

Recall first the notion of Hecke algebra with parameters defined in [Lus89].

Let $(\Lambda, \Lambda^\vee, \Sigma, \Sigma^\vee, \Delta)$ a quintuplet, where $\Lambda$ and $\Lambda^\vee$ are free abelian groups of finite type with duality by a $\mathbb{Z}$-bilinear map $\langle\cdot,\cdot\rangle : \Lambda \times \Lambda^\vee \to \mathbb{Z}$, $\Sigma \subseteq \Lambda$ a root system, $\Delta$ a basis of $\Sigma$ and $\Sigma^\vee \subseteq \Lambda^\vee$ the dual root system of $\Sigma$, the duality being given by $\langle\cdot,\cdot\rangle$.

Denote by $W(\Sigma)$ the Weyl group of $\Sigma$. Denote by $\Sigma_1, \cdots, \Sigma_r$ the irreducible components of $\Sigma$ and $\Delta_i := \Delta \cap \Sigma_i$. For every $i \in \{1, \cdots, r\}$, suppose we're given a set $\{q_\alpha, \alpha \in \Delta_i\}$ of real numbers $> 1$ such that $q_\alpha = q_\beta$, if $\alpha$ and $\beta$ are conjugate by an element of Weyl group of $\Sigma_i$. If $\Sigma_i$ is of type $B_n$, we moreover give ourselves a real number $q_i > 1$.

When $\alpha$ is in $\Sigma$, denote by $s_\alpha$ the elementary symmetry of $\Lambda$ associated to $\alpha$, and, for $\alpha, \beta \in \Sigma$, by $m(\alpha,\beta)$ the order of $s_\alpha s_\beta$.

Consider the group $B_0(\Sigma)$ of generators $U_{s_\alpha}$, $\alpha \in \Delta$, with the braid relations $U_{s_\alpha}U_{s_\beta}U_{s_\alpha}\cdots = U_{s_\beta}U_{s_\alpha}U_{s_\beta}\cdots$ ($m(\alpha,\beta)$ factors) for every $\alpha, \beta \in \Delta$. Thus note that $U_{s_{\alpha_1}}\cdots U_{s_{\alpha_r}} = U_{s_{\beta_1}}\cdots U_{s_{\beta_r}}$, when $s_{\alpha_1}\cdots s_{\alpha_r}$ and $s_{\beta_1}\cdots s_{\beta_r}$ are reduced decompositions of the same element $w$ of $W(\Sigma)$ [Spr98, 8.3.3]. We can thus set $U_w := U_{s_{\alpha_1}}\cdots U_{s_{\alpha_r}}$. Denote by $\mathcal{H}_0(\Sigma)$ the quotient of the group algebra of $B_0(\Sigma)$ by the two-sided ideal generated by the elements $(U_{s_\alpha}+1)(U_{s_\alpha}-q_\alpha)$, $\alpha\in\Delta$.

Denote by $C=\mathbb{C}[\Lambda]$ the group algebra of $\Lambda$ and $Z_\lambda$ the element of $C$ associated to $\lambda\in\Lambda$.

We will call *Hecke algebra with parameters* $\{q_\alpha\}$ *and* $\{q_i\}$, and we denote by $\mathcal{H}:=\mathcal{H}(\Sigma,\{q_\alpha\},\{q_i\})$ the $\mathbb{C}$-algebra which, as $\mathbb{C}$-vector-space is generated by $U_wZ_\lambda$, $w$ and $\lambda$ running through respectively $W(\Sigma)$ and $\Lambda$. The multiplication in $\mathcal{H}(\Sigma,\{q_\alpha\},\{q_i\})$ is deduced from that in $\mathcal{H}_0(\Sigma)$ and from that in $C$ with the commutation rule

$$Z_\lambda U_{s_\alpha}-U_{s_\alpha}Z_{s_\alpha\lambda}=\begin{cases}(q_\alpha-1)\frac{Z_\lambda-Z_{s_\alpha(\lambda)}}{1-Z_{-\alpha}}, & \text{if } \alpha^\vee\notin 2\Lambda^\vee,\\ (q_\alpha-1+Z_{-\alpha}((q_\alpha q_i)^{1/2}-(q_\alpha q_i^{-1})^{1/2}))\frac{Z_\lambda-Z_{s_\alpha(\lambda)}}{1-Z_{-2\alpha}}, & \text{if not;}\end{cases}$$

for $\alpha\in\Delta$ and $\lambda\in\Lambda$.

This is well-defined because $\alpha^\vee\in 2\Lambda^\vee$ is equivalent to saying that the irreducible component $\Sigma_i$ of $\Sigma$ to which $\alpha$ belongs is of type $B_n$ and that $\alpha$ is the short root in $\Delta\cap\Sigma_i$.

Denote by $\mathcal{Z}$ the center of $\mathcal{H}$, $K(\mathcal{Z})$ its fraction field, and set $\mathcal{H}_{K(\mathcal{Z})}=\mathcal{H}\otimes_{\mathcal{Z}}K(\mathcal{Z})$. Following [Lus89, 3.12], we have $\mathcal{H}_{K(\mathcal{Z})}=\bigoplus_{w\in W(\Sigma)}K(C)U_w$, where $K(C)$ denotes the fraction field of $C$.

7.2

We will now come back to our algebra $\mathrm{End}_G(i_P^GE_{B_\mathcal{O}})$ and establish the properties of operators $T_w$ to show that it is indeed a Hecke algebra with parameters.

**Lemma** *Let $\alpha,\alpha'\in\Delta_\mathcal{O}$, $s=s_\alpha$, $s'=s_{\alpha'}$, and $m\in\mathbb{Z}$ such that $(ss')^m=1$. Suppose that $b_s=b_{s'}=0$, or that $s$ and $s'$ commute, or that $\alpha$ and $\alpha'$ generate a root system of type $B_2$, that $\alpha'$ is the short root and that $b_s=0$. Thus $T_sT_{s'}T_s\cdots=T_{s'}T_sT_{s'}\cdots$, the number of factors of each side being $m$.*

*Remark* It seems that the conclusions of the above lemma will be wrong, if we omit one of the hypotheses.

*Proof.* The root system $\Sigma_1$ generated by $\alpha$ and $\alpha'$ is reduced of rank 2. If it is reducible, it is of type $A_1\times A_1$, and $J_s$ and $J_{s'}$ commute. But then $T_s$ and $T_{s'}$ commute also, and the assertion is trivial.

If not, $\Sigma_1$ is of type $A_2$, $B_2$ or $G_2$. Without loss of generality, we can suppose that the root $\alpha'$ is at most as long as $\alpha$. Denote by $m$ the order of $s_\alpha s_{\alpha'}$. We set $q_\alpha=q^{a_s}$, $q_{\alpha'}=q^{a_{s'}+b_{s'}}$ and $q_0=q^{a_{s'}-b_{s'}}$. Consider the quintuplet $(\Lambda_\mathcal{O},\Lambda_\mathcal{O}^\vee,\Sigma_1,\Sigma_1^\vee,\{\alpha,\alpha'\})$ and the associated Hecke algebra with parameters $\mathcal{H} = \mathcal{H}(\Sigma_1, \{q_\alpha, q_{\alpha'}\}, \{q_0\})$. Remark that $\mathbb{C}[\Lambda_\mathcal{O}]$ is identified with $B_\mathcal{O}$. Denote by $\mathcal{Z}$ the center of $\mathcal{H}$. The algebra $\mathcal{H}_{K(\mathcal{Z})}$ is thus a $K(B_\mathcal{O})$-module.

With the notations of **7.1**, set

$$j_\alpha = \frac{q_\alpha - Z_\alpha}{1 - Z_\alpha} \qquad \text{and} \qquad j_{\alpha'} = \frac{(q_{\alpha'}^{1/2} q_0^{1/2} - Z_{\alpha'})(Z_{\alpha'} + q_{\alpha'}^{1/2} q_0^{-1/2})}{1 - Z_{\alpha'}^2}.$$

For $P \in C$, $\beta \in \{\alpha, \alpha'\}$ and $t = s_\beta$, the commutation law implies that

$$(*) \qquad P(U_\beta - q_\beta) - (U_\beta - q_\beta)^t P = ({}^t P - P) j_\beta.$$

If $\beta \in \Sigma_1$ is conjugate to $\alpha$ (resp. $\alpha'$), put $j_\beta = j_\alpha(Y_\beta)$ (resp. $j_\beta = j_{\alpha'}(Y_\beta)$).

In $\mathcal{H}_{K(\mathcal{Z})}$, set $S_s = U_s - q_\alpha + j_\alpha$, $S_{s'} = U_{s'} - q_{\alpha'} + j_{\alpha'}$, and, when $w \in W(\Sigma_1)$ and that $w = s_{\alpha_1} \cdots s_{\alpha_r}$ is a reduced decomposition of $w$, $S_w = S_{\alpha_1} \cdots S_{\alpha_r}$. This does not depend on the reduced decomposition chosen [Ree97, 4.3]. Thus, for every $w, w' \in W(\Sigma_1)$,

$$S_w S_{w'} = \Big(\prod_\beta j_\beta j_{-\beta}\Big) S_{ww'},$$

the product taken over the set of $\beta \in \Sigma_1^+$ such that $w^{-1}\beta < 0$ and $w'^{-1} w^{-1} \beta > 0$ [Ree97, 4.3].

We will show that there is a unique homomorphism of algebras $\mathcal{H}_{K(\mathcal{Z})} \to \mathrm{End}_G(i_P^G E_{B_\mathcal{O}})_{K(\mathcal{Z}_\mathcal{O})}$ which sends $S_s$ to $R_s$, $S_{s'}$ to $R_{s'}$ and which induces the identity on $K(B_\mathcal{O})$. It will follow the relation $T_s T_{s'} T_s \cdots = T_{s'} T_s T_{s'} \cdots$, the number of factors of each side being equal to the order of $ss'$, since $U_s \mapsto T_s$ and that $U_{s'} \mapsto T_{s'}$.

Since the operators $S_s$ and $S_{s'}$ generate together with $K(B_\mathcal{O})$ the algebra $\mathcal{H}_{K(\mathcal{Z})}$, the uniqueness is evident. It remains to prove the existence, i.e. that $R_s$, $R_{s'}$ satisfy the same relations as $S_s$, $S_{s'}$.

The commutation relations $(*)$ for the operatos $R_s$ and $R_{s'}$ follow directly from their definition **5.2** with the lemma **4.8**. It is the same with identities $R_s^2 = j_\alpha j_{-\alpha}$ and $R_{s'}^2 = j_{\alpha'} j_{-\alpha'}$ with the lemma **4.7**. It remains to verify that, when $w \in W(\Sigma_1)$ and that $w = s_{\alpha_1} \cdots s_{\alpha_r}$ is a decomposition into simple symmetries of $w$, then $R_{s_{\alpha_1}} \cdots R_{s_{\alpha_r}}$ is independent of choice of this decomposition. For this, we will study the different case. Indeed, if $m(s, s')$ denotes the order of $ss'$, it suffices to show that $R_s R_{s'} R_s \cdots = R_{s'} R_s R_{s'} \cdots$, the number of factors of each side being $m(s, s')$ [Spr98, 8.3.3].

If $\Sigma_1$ is of type $A_2$ or $G_2$, then $b_{s'} = 0$ by hypothesis and $c_{s'} = c_s$, since $s$ and $s'$ are conjugate. The relations above are thus a direct consequence of corollary **3.4**.

Suppose for now $\Sigma_1$ is of type $B_2$. It is necessary that $R_{s_1} R_{s_2} R_{s_1} R_{s_2} = R_{s_2} R_{s_1} R_{s_2} R_{s_1}$. Since each factor $R_{s_i}$ appears with the same multiplicity, the scalar factors are equal. We are therefore brought back to the equality

$$X_\alpha J_{s_\alpha} X_{\alpha'} J_{s_{\alpha'}} X_\alpha J_{s_\alpha} X_{\alpha'} J_{s_{\alpha'}} = X_{\alpha'} J_{s_{\alpha'}} X_\alpha J_{s_\alpha} X_{\alpha'} J_{s_{\alpha'}} X_\alpha J_{s_\alpha}$$

which follows from lemma **2.7**. □

7.3

**Proposition** *When $\Sigma_i$ is an irreducible component of $\Sigma_{\mathcal{O}}$ and that $\alpha \in \Delta \cap \Sigma_i$, then $b_{s_\alpha} \neq 0$ implies that $\Sigma_i$ is of type $B_n$ and that $\alpha$ corresponds to the short root of $B_n$.*

*Proof.* It follows from the work of Bernstein-Zelevinsky [BZ77, 4.2] and, in the case of a simple algebra, of [Tad90] that $b_{s_\alpha} = 0$ if the function $\mu^{M_\alpha}$ is equal to that for a general linear group or the multiplicative group of a simple algebra. This is also true for $\mathrm{SL}_2(F)$ and $\mathrm{PGL}_2(F)$. Now, in the other case, $\Sigma_i$ is, according to **1.13** and **6.1**, of type $B_n$ and $\alpha$ is the short root. □

7.4

When $\Sigma_i$ is an irreducible component of $\Sigma_{\mathcal{O}}$, denote by $\alpha_i$ the unique root in $\Delta \cap \Sigma_i$ such that $\alpha_i^\vee \in 2\Lambda_{\mathcal{O}}^\vee$, if such a root exists, and set $a_i = a_{s_{\alpha_i}}$ and $b_i = b_{s_{\alpha_i}}$.

**Proposition** *The algebra $\bigoplus_{w \in W_{\mathcal{O}}} B_{\mathcal{O}} T_w$ is a Hecke algebra with parameters $\{q^{a_{s_\alpha}+b_{s_\alpha}}\}$ and $\{q^{a_i - b_i}\}$ and given basic radical $(\Lambda_{\mathcal{O}}, \Sigma_{\mathcal{O}}, \Lambda_{\mathcal{O}}^\vee, \Sigma_{\mathcal{O}}^\vee, \Delta_{\mathcal{O}})$, $i$ correspond to different irreducible components of $\Sigma_{\mathcal{O}}$.*

*Proof.* The only property which remains to verify is the commutation law. This is an elementary calculation. □

7.5

Note that the possible values for the parameters $a_s$ and $b_s$ are known for the groups considered here via the Langlands correspondence: the order of group of unramified characters which stabilize the isomorphism class of a cuspidal irreducible representation of $\mathrm{GL}_n(F)$ is a divisor of $n$. (This follows from the local Langlands correspondance [HT01].) We have a similar result in the case of simple algebras by the Jacquet-Langlands correspondance [DKV84]. In the case where $G$ is equal to $\mathrm{GL}_n(F)$ or the multiplicative group of a simple algebra, we deduce the values of $a_s$ and $b_s$ directly from works cited in the proof of **7.3**. For $H_k$ a symplectic or orthogonal group and $\sigma \otimes \tau$ a cuspidal irreducible representation of $\mathrm{GL}_n(F) \times H_k$ such that $\mu(\sigma \otimes \tau) = 0$, C. Moeglin [Mg11] recently deducted from work of J. Arthur that the number $a > 0$ such that $\sigma |\det_n|_F^a \otimes \tau$ is a pole of $\mu$ is a half-integer. Moreover, if $\chi$ is an unramified unitary character such that the function $\mu$ is also cancelled at $\sigma\chi \otimes \tau$, then the real number $b > 0$ such that $(\sigma\chi)|\det|_F^b \otimes \tau$ is a pole of $\mu$ multiplied by the order of stabilizer $t$ of $\sigma$ is an integer if and only if $ta$ is. Other than that,

unless $\sigma$ and $\sigma\chi$ are isomorphic, the value of $b$ is in general not conditional on that of $a$. The reader can consult [Hei10] where this is made completely explicit and applied to Hecke algebras with parameters obtained here.

Determining the values of the parameters $a_s$ and $b_s$ directly on $\sigma$, for example by the theory of types, seems to be a very hard and difficult problem (see for example [Blo12] for a partial result).

7.6

**Proposition** *Let $w \in W_{\mathcal{O}}$ and $r \in R(\mathcal{O})$. Thus $r^{-1}wr \in W_{\mathcal{O}}$ and $T_w J_r = J_r T_{r^{-1}wr}$.*

*Proof.* It suffices to consider the case where $w$ is a simple symmetry $s_\alpha$ and to show that $r^{-1}s_\alpha r$ is also a simple symmetry. Now, since by definition $r$ leaves $\Sigma_{\mathcal{O}} \cap \Sigma(P)$ invariant, it is the same for $\Delta_{\mathcal{O}}$, and, as a result, $r^{-1}s_\alpha r = s_{r\alpha}$ is a simple symmetry. □

7.7

**Theorem** *The algebra $\mathrm{End}_G(i_P^G E_{B_{\mathcal{O}}})$ is isomorphic to semi-direct product $\mathbb{C}[R(\mathcal{O})] \ltimes \mathcal{H}(\Sigma_{\mathcal{O}}, \{q^{a_{s_\alpha}+b_{s_\alpha}}\}, \{q^{a_i-b_i}\})$ of the group algebra of $R(\mathcal{O})$ with the Hecke algebra with parameters $\{q^{a_{s_\alpha}+b_{s_\alpha}}\}$ and $\{q^{a_i-b_i}\}$.*

*Proof.* This is an immediate consequence of **5.10, 7.4** and of **7.6**. □

7.8

**Corollary** *The category $\mathrm{Rep}_G({}^W\mathcal{O})$ is isomorphic to the category of right modules over the algebra $\mathbb{C}[R(\mathcal{O})] \ltimes \mathcal{H}(\Sigma_{\mathcal{O}}, \{q^{a_{s_\alpha}+b_{s_\alpha}}\}, \{q^{a_i-b_i}\})$.*

*Proof.* This is an immediate consequence of **7.7** combined with the theorem of Bernstein [Roc02, 1.6] mentioned in the introduction. □

7.9

*Remark* The isomorphism of categories of **7.8** is compatible with the parabolic induction in the following sense: remark first that the theorem **7.7** and the corollary **7.8** generalize easily to subgroups of Levi of groups considered here. Indeed, such a subgroup of Levi is a direct product of groups of this type, and the algebra considered are thus the product of algebras associated to the factors.

Then let $P' = M'U'$ be a parabolic subgroup of $G$ such that $P' \supseteq P$ and $M' \supseteq M$. For $\alpha \in \Delta_{\mathcal{O}}^{M'}$ and $r \in R^{M'}(\mathcal{O})$, denote by $T_{s_\alpha}^{M'}$ and $J_r^{M'}$ the operators $T_{s_\alpha}$ and $J_r$ defined relative to $M'$. We have a canonical injection of $\mathrm{End}_{M'}(i^{M'}_{P\cap M'}E_{B_\mathcal{O}})$ into $\mathrm{End}_G(i^G_P E_{B_\mathcal{O}})$ defined by the functor of parabolic induction $i^G_{P'}$. We thus observe that this inclusion corresponds to the natural inclusion of respective Hecke algebras with parameters (also the algebras of finite group). Indeed, this inclusion of $\mathrm{End}_{M'}(i^{M'}_{P\cap M'}E_{B_\mathcal{O}})$ into $\mathrm{End}_G(i^G_P E_{B_\mathcal{O}})$ sends an operator $T^{M'}_{s_\alpha}$, $\alpha \in \Delta^{M'}_\mathcal{O}$, to $T_{s_\alpha}$ and an operator $J^{M'}_r$, $r \in R^{M'}(\mathcal{O})$, to $J_r$, the construction of these operators being compatible with the parabolic induction (for the intertwining operators see [Wal03, IV.1(14)]) and being done in relation to a Levi subgroup contained in $M'$.

Recall, as already noted in the introduction, that the isomorphism of categories between $\mathrm{Rep}_G({}^W\mathcal{O})$ and the category of right $\mathrm{End}_G(i^G_P E_{B_\mathcal{O}})$-modules is itself also compatible with the parabolic induction [Roc02, 5.3]. The isomorphism of category **7.8** is thus compatible with the parabolic induction.

We see likewise that it is compatibel with the Jacquet functor: indeed, the Jacquet functor $r^G_{\overline{P'}}$ defines a functor of the category $\mathrm{Rep}_G({}^W\mathcal{O})$ to the category $\mathrm{Rep}_{M'}({}^{W^{M'}}\mathcal{O})$. On the other hand, by the inclusion of $\mathrm{End}_{M'}(i^{M'}_{P\cap M'}E_{B_\mathcal{O}})$ into $\mathrm{End}_G(i^G_P E_{B_\mathcal{O}})$, every right $\mathrm{End}_G(i^G_P E_{B_\mathcal{O}})$-module has a structure of right $\mathrm{End}_{M'}(i^{M'}_{P\cap M'}E_{B_\mathcal{O}})$-module, which induces a functor between these two categories of modules. By [Roc02, 5.3], this commutes with the Jacquet functor.

Author: V.H., current address: Université d'Aix-Marseille, I2M - CNRS (UMR 7373), Campus de Luminy, Luminy Case 907, 13288 Marseille Cédex 9, France

*Email address*: volker.heiermann@univ-amu.fr

(Translator: Y.X., current address: Columbia University, 2990 Broadway, New York, NY 10027, USA

*Email address*: xu.yujie@columbia.edu)